
  \input amssym
  \input miniltx
  \input pictex
    \input graphicx.sty


  \font \bbfive = bbm5
  \font \bbeight = bbm8
  \font \bbten = bbm10
  \font \rs = rsfs10 \font \rssmall = rsfs10 scaled 833 \font \rstiny = rsfs10 scaled 600  
  \font \eightbf = cmbx8
  \font \eighti = cmmi8 \skewchar \eighti = '177
  \font \fouri = cmmi5 scaled 800 
  \font \eightit = cmti8
  \font \eightrm = cmr8
  \font \eightsl = cmsl8
  \font \eightsy = cmsy8 \skewchar \eightsy = '60
  \font \eighttt = cmtt8 \hyphenchar \eighttt = -1

  \font \sixi = cmmi6 \skewchar \sixi = '177
  \font \sixrm = cmr6
  \font \sixsy = cmsy6 \skewchar \sixsy = '60
  \font \tensc = cmcsc10

  \scriptfont \bffam = \bbeight
  \scriptscriptfont \bffam = \bbfive
  \textfont \bffam = \bbten

  \newskip \ttglue

  \def \eightpoint {\def \rm {\fam 0 \eightrm }\relax
  \textfont 0 = \eightrm \scriptfont 0 = \sixrm \scriptscriptfont 0 = \fiverm
  \textfont 1 = \eighti \scriptfont 1 = \sixi \scriptscriptfont 1 = \fouri
  \textfont 2 = \eightsy \scriptfont 2 = \sixsy \scriptscriptfont 2 = \fivesy
  \textfont 3 = \tenex \scriptfont 3 = \tenex \scriptscriptfont 3 = \tenex
  \def \it {\fam \itfam \eightit }\relax
  \textfont \itfam = \eightit
  \def \sl {\fam \slfam \eightsl }\relax
  \textfont \slfam = \eightsl
  \def \bf {\fam \bffam \eightbf }\relax
  \textfont \bffam = \bbeight \scriptfont \bffam = \bbfive \scriptscriptfont \bffam = \bbfive
  \def \tt {\fam \ttfam \eighttt }\relax
  \textfont \ttfam = \eighttt
  \tt \ttglue = .5em plus.25em minus.15em
  \normalbaselineskip = 9pt
  \def \MF {{\manual opqr}\-{\manual stuq}}\relax
  \let \sc = \sixrm
  \let \big = \eightbig
  \let \rs = \rssmall
  \setbox \strutbox = \hbox {\vrule height7pt depth2pt width0pt}\relax
  \normalbaselines \rm }

  \def \setfont #1{\font \auxfont =#1 \auxfont }
  \def \withfont #1#2{{\setfont {#1}#2}}


  \def \text #1{\mathchoice {\hbox {#1}} {\hbox {#1}} {\hbox {\eightrm #1}} {\hbox {\sixrm #1}}}
  \def \varbox #1{\setbox 0\hbox {$#1$}\setbox 1\hbox {$I$}{\ifdim \ht 0< \ht 1 \scriptstyle #1 \else \scriptscriptstyle #1 \fi }}

  \def \rsbox #1{{\mathchoice {\hbox {\rs #1}} {\hbox {\rs #1}} {\hbox {\rssmall #1}} {\hbox {\rstiny #1}}}}
  \def \mathscr #1{\rsbox {#1}}

  \def \paren #1{\setbox 0\hbox {$#1$} \htzero =\ht 0
    \setbox 1\hbox {$($}
      \ifdim \htzero< \ht 1
        (#1)
      \else
        \setbox 1\hbox {$\big($}
          \ifdim \htzero< \ht 1
            \big(#1\big)
          \else
            \left(#1\right)
          \fi
      \fi}


  \def \VarItem #1#2{\smallskip \par \noindent \kern #1 \hangindent #1 \llap {#2\enspace }\ignorespaces }
  \def \item #1{\VarItem {1.1truecm}{\rm #1}}


  \def \ifundef #1{\expandafter \ifx \csname #1\endcsname \relax }

  \def \section #1 #2 \par {%
    \goodbreak \bigbreak
    \null\vskip 0.5cm
    \centerline {\tensc #1.\enspace #2}
    \null\vskip 0.3cm
    \edef \rightRunningHead {\ifundef {authorRunninhHead}#1\else \authorRunninhHead \fi }
    \Headlines {#2}{\rightRunningHead }\par}

  \def \firstLine {\medskip \noindent }

  \def \state #1 #2. #3\par {\edef \a {Definition} \edef \b {#2}
    \medbreak \noindent {\bf #1\enspace #2.\enspace \ifx \a \b \rm \else \sl \fi
    #3\par } \medbreak }

  \def \Proof {\global \def \iP {Y}
    \medbreak \noindent {\it Proof.\enspace }}

  \def \endProof {\if \iP Y \hfill $\endproofmarker $ \looseness = -1 \fi
    \medbreak}

  \def \closeProof {\eqno \endproofmarker
    \global \def \iP {N}}

  \def \possundef #1{\csname #1\endcsname}


  \def \Bibitem #1 #2; #3; #4 \par {\VarItem {0.5truecm}{[\possundef {#1}]} #2, {``#3''}, #4.\par}

  \def \jrn #1, #2 (#3), #4-#5;{{\sl #1}, {\bf #2} (#3), #4--#5}
  \def \Article #1 #2; #3; #4 \par {\Bibitem #1 #2; #3; \jrn #4; \par }

  \def \references {\begingroup \bigbreak \eightpoint \centerline {\tensc References} \nobreak \medskip \frenchspacing }



  %
  \catcode `@ = 11
  \def \@Nnil {\@Nil }%
  \def \@Fornoop #1\@@ #2#3{}%
  \def \For #1:=#2\do #3{%
     \edef \@Fortmp {#2}%
     \ifx \@Fortmp \empty \else
        \expandafter \@Forloop #2,\@Nil ,\@Nil \@@ #1{#3}%
     \fi
  }%
  \def \@Forloop #1,#2,#3\@@ #4#5{\@Fordef #1\@@ #4\ifx #4\@Nnil \else
         #5\@Fordef #2\@@ #4\ifx #4\@Nnil \else #5\@iForloop #3\@@ #4{#5}\fi \fi
  }%
  \def \@iForloop #1,#2\@@ #3#4{\@Fordef #1\@@ #3\ifx #3\@Nnil
         \let \@Nextwhile =\@Fornoop \else
        #4\relax \let \@Nextwhile =\@iForloop \fi \@Nextwhile #2\@@ #3{#4}%
  }%
  \def \@Forspc { }%
  \def \@Fordef {\futurelet \@Fortmp \@@Fordef }
  \def \@@Fordef {%
    \expandafter \ifx \@Forspc \@Fortmp 
      \expandafter \@Fortrim
    \else
      \expandafter \@@@Fordef
    \fi
  }%
  \expandafter \def \expandafter \@Fortrim \@Forspc #1\@@ {\@Fordef #1\@@ }%
  \def \@@@Fordef #1\@@ #2{\def #2{#1}}%
  %
  %
  \def \citetrk #1{{{\def \a {}\For \name :=#1\do {\a {\bf \possundef {\name }}\def \a {, }}}}} 
  \def \c@ite #1{{\rm [\citetrk {#1}]}}
  \def \sc@ite [#1]#2{{\rm [\citetrk {#2}\hskip 0.7pt:\hskip 2pt #1]}}
  \def \du@lcite {\if \pe@k [\expandafter \sc@ite \else \expandafter \c@ite \fi }
  \def \cite {\futurelet \pe@k \du@lcite }
  \catcode `\@ =12
  %


  \def \Headlines #1#2{\nopagenumbers
    \headline {\ifnum \pageno = 1 \hfil
    \else \ifodd \pageno \tensc \hfil \lcase {#1} \hfil \folio
    \else \tensc \folio \hfil \lcase {#2} \hfil
    \fi \fi }}

  \font \titlefont = cmbx12
  \long \def \title #1{\begingroup
    \titlefont
    \parindent = 0pt
    \baselineskip = 16pt
    \leftskip = 35truept plus 1fill
    \rightskip = \leftskip
    #1\par \endgroup }

  \long \def \Quote #1\endQuote {\begingroup \leftskip 35truept \rightskip 35truept \parindent 17pt \eightpoint #1\par \endgroup }
  \long \def \Abstract #1\endAbstract {\vskip 1cm \Quote \noindent #1\endQuote }
  \def \Address #1#2{\bigskip {\tensc #1 \par \it E-mail address: \tt #2}}
  \def \Authors #1{\bigskip \centerline {\tensc #1}}
  \def \Note #1{\footnote {}{\eightpoint #1}}
  \def \Date #1 {\Note {\it Date: #1.}}


  \def \Case #1:{\medskip \noindent {\tensc Case #1:}}

  \def \fix {\smallskip \noindent $\blacktriangleright $\kern 12pt}

  \def \lcase #1{\edef \auxvar {\lowercase {#1}}\auxvar }
  \def \emph #1{{\it #1}\/}
  \def \explica #1#2{\mathrel {\buildrel \hbox {\sixrm #1} \over #2}}

  \newcount \fnctr \fnctr = 0
  \def \fn #1{\global \advance \fnctr by 1
    \edef \footnumb {$^{\number \fnctr }$}\relax
    \footnote {\footnumb }{\eightpoint #1\par \vskip -10pt}}


  \def \mathbb #1{{\bf #1}}
  \def \frac #1#2{{#1\over #2}}
  \def \<{\left \langle \vrule width 0pt depth 0pt height 8pt }
  \def \>{\right \rangle }
  \def \({\left (}
  \def \){\right )}
  
  \def \and {\mathchoice {\hbox {\quad and \quad }} {\hbox { and }} {\hbox { and }} {\hbox { and }}}

  \def \imply {\mathrel {\Rightarrow }}
  \def \IFF {\kern 7pt\Leftrightarrow \kern 7pt}
  \def \IMPLY {\kern 7pt \Rightarrow \kern 7pt}
  \def \for #1{\mathchoice {\quad \forall \,#1} {\hbox { for all } #1} {\forall #1}{\forall #1}}
  \def \endproofmarker {\square }
  \def \"#1{{\it #1}\/} \def \umlaut #1{{\accent "7F #1}}
  
  \def \*{\otimes }
  \def \caldef #1{\global \expandafter \edef \csname #1\endcsname {{\cal #1}}}
  \def \mathcal #1{{\cal #1}}
  \def \bfdef #1{\global \expandafter \edef \csname #1\endcsname {{\bf #1}}}
  \bfdef N \bfdef Z \bfdef C \bfdef R
  \def \percent {\char "25}
  \def \exists {\mathchar "0239\kern 1pt }
  \def \labelarrow #1{\setbox 0\hbox {\ \ $#1$\ \ }\ {\buildrel \textstyle #1 \over {\hbox to \wd 0 {\rightarrowfill }}}\ }
  \def \subProof #1{\medskip \noindent #1\enspace }
  \def \itmProof (#1) {\subProof {(#1)}}
  \def \itemImply #1#2{\subProof {#1$\Rightarrow $#2}}
  \def \itmImply (#1) > (#2) {\itemImply {(#1)}{(#2)}}

  \def \close {\end }

  \def \medprod {\mathop {\mathchoice {\hbox {$\mathchar "1351$}}{\mathchar "1351}{\mathchar "1351}{\mathchar "1351}}}
  
  \def \medcup {\mathop {\mathchoice {\raise 1pt \hbox {$\mathchar "1353$}}{\mathchar "1353}{\mathchar "1353}{\mathchar "1353}}}
  \def \medcap {\mathop {\mathchoice {\raise 1pt \hbox {$\mathchar "1354$}}{\mathchar "1354}{\mathchar "1354}{\mathchar "1354}}}

  \def \clauses #1{\def\quad{\enspace}\left \{\matrix {#1}\right .}
  \def \cl #1#2{\hfill #1\hfill ,&#2\hfill \vrule width 0pt height 12pt depth 6pt\cr }

  \def \paperAll #1#2#3#4#5#6{
    \hsize #1truemm   \advance \hsize by -#3truemm  \advance \hsize by -#4truemm
    \vsize #2truemm   \advance \vsize by -#5truemm  \advance \vsize by -#6truemm
    \hoffset =-1truein \advance \hoffset by #3truemm
    \voffset =-1truein \advance \voffset by #5truemm
    }

  \def \paperHV #1#2#3#4{
    \paperAll {#1}{#2}{#3}{#3}{#4}{#4}
    }

  \def \bool #1{\, [\, {#1}\, ]}

\def\iff{\Leftrightarrow} \def\proofImply#1#2{(#1$\Rightarrow$#2)\hskip6pt} \def\Medskip{\medskip\noindent}
 \def\tight{\tau} \def\bolinha{\raise1.5pt \hbox{\kern1pt$\scriptscriptstyle\bullet$\kern1pt}} \def\bitem{
  \item {$\bolinha$}} \def\pbitem{
  \item {($\bolinha$)}} \def\tsp#1{\widehat{#1}_\tight} \def\hath{\hat h} \def\hatE{\widehat E} \def\hatF{\widehat F}

\def\germ#1#2#3{\big[#1,#2\big]_{#3}} \def\gt#1#2{\germ{#1}{#2}{\theta}} \def\ga#1#2{\germ{#1}{#2}{\alpha}}
\def\gb#1#2{\germ{#1}{#2}{\beta}} \def\GermGpd#1{G_{#1}}

\def\otg/{topological groupoid with order} \def\otgs/{topological groupoids with order} \def\Up{{\cal U}}

\def\Cpl#1{#1^{\tight}}  \def\boo{\rsbox{B}} \def\ch{\check h} \def\d{d} \def\D{D} \def\O{{\cal O}}
 \def\equiva{\approx} \def\menor{\preccurlyeq} \def\ker{\text{Ker}} \def\dom{\text{dom}}
\def\ran{\text{ran}} \def\relimply{\ \Rightarrow\ } \def\reliff{\ \Leftrightarrow\ } \def\op{^{op}}
\def\Gt#1{G_\tight(#1)} \def\Tcal{{\rsbox{T\kern3pt}}} \def\Ecal{{\rsbox{E}}} \def\pointfun{\varphi} \def\setfun{\Phi}

\def\ex#1{^{\scriptscriptstyle(#1)}} \def\sm{\kern-3pt\setminus} \def\smallsup#1{^{\raise2pt \hbox{\sixrm#1}}}
\def\hz{h\smallsup{o}} \def\hzao{H} \def\chz{\check h\smallsup{o}} \def\hathz{\hat h\smallsup{o}}

\def\gt#1#2{\big[#1,#2\big]_\theta} \def\ga#1#2{\big[#1,#2\big]_\alpha} \def\gb#1#2{\big[#1,#2\big]_\beta}

\magnification1200 \overfullrule= 0cm \paperHV{210}{297}{39}{38}

\title{CONSONANT INVERSE SEMIGROUPS}

\def\authorRunninhHead{R. Exel} \Authors{\authorRunninhHead}

\Note{Partially supported by CNPq and FAPESP.}

\Abstract We study necessary and sufficient conditions for two inverse semigroups to possess identical tight groupoids
from the point of view of their algebraic, topological, and spectral order structures. The spectral order is a partial
order relation that presents itself very naturally on the tight groupoid of an inverse semigroup and is related to the
subtle difference between tight filters and ultra-filters. Up to a small glitch, the spectral order makes tight
groupoids into ordered groupoids in Ehresmann's sense.

We introduce the notions of \emph{tight injectivity} and \emph{tight surjectivity} for inverse semigroup homomorphisms
and show that, together, they provide necessary and sufficient conditions for the induced map to be an isomorphism of
ordered topological groupoids. A homomorphism is then called a \emph{consonance} provided these conditions are met.

A consonance is not necessarily injective or surjective, so it is likely not an invertible map. Due to this fact, the
existence of a consonance between inverse semigroups does not define an equivalence relation, so we consider instead the
equivalence relation it generates. When it applies to a pair of inverse semigroups we say that they are
\emph{consonant}.

The first main result of the paper shows that two inverse semigroups $S_1$ and $S_2$ are consonant if and only if their
tight groupoids are isomorphic as ordered topological groupoids, if and only if there exists another inverse semigroup
$T$ admitting consonances from each $S_i$ to $T$.

We also prove that, given an inverse semigroup $S$, there exists a largest inverse semigroup consonant to $S$, denoted
$\Cpl S$, and called the tight envelope of $S$. In a nutshell, $\Cpl S$ is the inverse semigroup formed by the compact
up-slices (i.e.~slices that are up-sets relative to the spectral order) in the tight groupoid of $S$. Among other things
$\Cpl S$ is the only flat, distributive, inverse semigroup that is consonant to $S$.

We also introduce the notion of \emph{tight-like} groupoids and, along the lines of Lawson and Vdovina's duality
theories, we show that the correspondence $ S\mapsto\Gt S $ is a bijection between the isomorphism classes of flat
distributive inverse semigroups and the isomorphism classes of tight-like groupoids. \endAbstract

\section 1 Introduction

\firstLine The genesis of the theory of Operator Algebras, largely attributed to the series of seminal papers ``On rings
of operators'' by Murray and von Neumann \cite{VNOne, VNTwo, VNThree, VNFour}, has been inextricably linked to dynamical
systems, as these were the main ingredients used to produce some of the first nontrivial examples of what is now known
as von Neumann algebras. Not much later, through the work of Zeller-Meier \cite{ZM}, dynamical systems have also become
one of the main sources of data from where C*-algebras were built. This trend continued unabated through the years,
including the revolution-inducing class of Cuntz-Krieger algebras \cite{Cuntz, CK}, whose main ingredient, a 0\kern1pt-1
matrix, is not explicitly a dynamical system, but which nevertheless turns out to be very closely related to a dynamical
system, this time the corresponding topological Markov chain. Graph C*-algebras \cite{Raeburn} ensued, wielding a strong
influence outside the field of Functional Analysis, and spawning the rich classes of Leavitt path algebras \cite{LPA}
and Steinberg algebras \cite{Steinberg, Lisa}.

Listing all instances in which a dynamical system, be it global or partial, reversible or irreversible, smooth,
continuous or symbolic, have influenced the construction of associative algebras, will extend this introduction far
beyond its goals. However there is one point of view that cannot be left out since it unifies all guises of dynamical
systems used in the construction of C*-algebras, and this is Renault's theory of groupoid C*-algebras \cite{Renault}.

In virtually all examples where C*-algebras are constructed from dynamical systems, there is a mediating groupoid, which
can be built out of the dynamical information, and from which the corresponding C*-algebra can be built. On the other
hand, in the very common case in which the relevant groupoid is \'etale and ample, there is a closely related object,
namely the inverse semigroup of compact open bisections, already recognized by Renault \cite[pag. 105]{Renault} as an
important player.

This motivated Paterson \cite{paterson} to construct an ample groupoid from an inverse semigroup, now known as the
Paterson's universal groupoid, which in turn was one of my main motivations to introduce the tight groupoid of an
inverse semigroup \cite{actions}. This theory has found many applications, and has attracted a great deal of interest
\cite{ Boava, BiceStar, DonMil, reconstru, covertojoin, ExelStar, EGS, EP, EPdois, ExelSteinbergOne, ExelSteinberg,
LaLMil, lawsonB, StoneDuality, LawsonLenz, lawsonV, MilanStein, OrtegaPardo, StarlingTwo, StarlingThree, Steinberg}. In
particular, many important algebras that can be encoded in an ample groupoid, can also be encoded in an inverse
semigroup, arguably a simpler object than a groupoid.

Among the many roles played by inverse semigroups across diverse areas of mathematics, here we focus on their role as
ingredients in the construction of \'etale groupoids. An inevitable consequence of this point of view is that it
immediately raises the question:

\medskip\centerline{\sl When do two inverse semigroups lead to the same groupoid?}

\medskip If interpreted in terms of Paterson's universal groupoid, this question doesn't seem to lead to a very deep
theory, or at least this is what our initial undertaking presented in Section (13) indicates. However, when the tight
groupoid is involved, the question becomes a lot more interesting and deep.

The raison d'\^etre of this work is thus to consider two inverse semigroups $S$ and $T$ and to try to decide whether or
not their tight groupoids are isomorphic in the most direct way possible, of course without going through the trouble of
constructing said groupoids and checking whether they are isomorphic.

Since inverse semigroups\fn{In this work all inverse semigroups are required to contain a zero element.} are equipped
with their all important idempotent semilattices, our first step is to tackle the above question for
semilattices. Noting that the tight groupoid of a semilattice $E$ is nothing but its tight spectrum $\tsp E$, seen as a
groupoid where all elements are units, we begin by determining necessary and sufficient conditions on a semilattice
homomorphism
  $$
  h:E\to F
  $$
  for the dual map
  $$
  \hath:\varphi\in\tsp F \mapsto\varphi\circ h\in\tsp E
  $$
  to be well defined and a homeomorphism. In fact, being well-defined is already a concern since $\varphi\circ h$ may
not be tight for every tight character $\varphi$ in $\tsp F$. In Corollary (4.4), we solve this problem by showing that
the well-definedness of $\hath$ is equivalent to the following conditions: \pbitem for every $f$ in $F$, there exists a
finite set $C\subseteq E$, such that $\{fh(e): e\in C\}$ is a cover for $f$, and \pbitem if $\{e_1,\ldots,e_n\}$ is a
cover for a given $e\in E$, then $\{h(e_1),\ldots,h(e_n)\}$ is a cover for $h(e)$.

\medskip Homomorphisms satisfying these conditions are called \emph{tight}. Assuming that $h$ is a tight homomorphism,
and hence that $\hath$ is well defined, we then address the question of when is $\hath$ a homeomorphism.

At this point a rather inconspicuous feature springs to life: the tight spectrum of a semilattice is more than just a
topological space -- it is equipped with an intrinsic order relation that has seemingly been overlooked until
now. Recalling that $\tsp E$ may be viewed as the set of tight filters on $E$, we may compare the points of $\tsp E$
according to whether one filter is a subset of the other. We call this the \emph{spectral order}.

The spectral order can be used e.g.~to identify the ultra-filters, which are, of course, the maximal elements, but this
is only the tip of the iceberg.

From this perspective, the question of whether the tight spectra of two semilattices are homeomorphic falls short,
insofar as it neglects the canonical order relation. Instead, the relevant question to ask is thus:

\begingroup\narrower

\medskip{\sl\noindent Given a tight homomorphism $h$, when is the dual map $\hath$ an isomorphism of ordered topological
spaces?}

\endgroup

\Medskip The answer is based on two main results, namely, \pbitem$\hath$ is surjective if and only if $h$ is
\emph{tightly injective} (Theorem 4.9), and \pbitem$\hath$ is \emph{order-injective} if and only if $h$ is \emph{tightly
surjective} (Theorem 5.4).

\medskip To explain the undefined terms above, we must first introduce a new order relation on a semilattice $E$ which,
like the spectral order, will play a fundamental role throughout this work: given $e$ and $f$ in $E$, we say that $e$ is
\emph{tightly below} $f$, in symbols $e \menor f$, if there is no nonzero $g$ in $E$, such that
  $$
  f\perp g \leq e.
  $$
  This is closely related to a congruence studied by Lawson and Lenz in \cite[Lemma 4.8]{LawsonLenz}, and it will be
called the \emph{tight order relation} on $E$ (see 3.1).

Among other things, one has that $e\menor f$ if and only if $\varphi(e)\leq\varphi(f)$ for every tight character
$\varphi$ on $E$.

It so happens that the tight order relation is reflexive and transitive, but it is not necessarily anti-symmetric. To
account for this, we consider an equivalence relation, namely the one studied by Lawson and Lenz, by saying that $e$ and
$f$ are \emph{tightly equivalent}, in symbols $e\equiva f$, whenever $e\menor f$, and $f\menor e$.

\emph{Tight injectivity} for a semilattice homomorphism $h:E\to F$ is then defined by
  $$
  h(e_1)\equiva h(e_2) \relimply e_1\equiva e_2.
  $$

On the other hand, \emph{order-injectivity}, the term used in the statement of Theorem (5.4) reproduced above, applies,
not to the tight order on a semilattice, but to the spectral order on the tight spectrum, and it means that
  $$
  \hath(\varphi_1)\leq\hath(\varphi_2) \relimply\varphi_1\leq\varphi_2,
  $$
  for any pair of tight characters $\varphi_1$ and $\varphi_2$ in $\tsp F$. Finally, the notion of \emph{tight
surjectivity} relies of the further concept of \emph{large} subsets of a semilattice $F$, and is defined as follows: a
subset $R\subseteq F$ is called \emph{large} if, for every $f$ in $F$, there exists a finite set $C\subseteq R$, such
that $e\menor f$, for all $e\in C$, and $\{fe: e\in C\}$ is a cover for $f$. With this, one says that a semilattice
homomorphism is \emph{tightly surjective} if its range is a large subset of the codomain.

The spectral order being anti-symmetric, every order-injective map is injective. Coupling this with the fact that
$\hath$ is always order-preserving (Proposition 5.3), we have that the combined effect of Theorems (4.9) and (5.4) is to
produce the desired answer to our question above: a necessary and sufficient condition for $\hath$ to be a well defined
isomorphism of ordered topological spaces is that $h$ be both tightly injective and tightly surjective, in which case we
say that $h$ is a \emph{consonance}.

Armed with a theory of consonances for semilattices, we then set our sights on inverse semigroups. The trouble, however,
is that there does not seem to be a sensible contravariant functor sending inverse semigroups to their tight groupoids,
as one would hope for based on the association $h \to\hath$.

Fortunately, since we are really only interested in the situation when $\hath$ is invertible, we may switch tacks and
focus instead on the inverse map $\hath^{-1}$, which we denote by $\ch$. This is greatly facilitated by the fact that
this map has a nice expression, namely
  $$
  \ch(\xi) = \big\{f\in F: h(e)\menor f, \text{ for some } e\in\xi\big\},
  $$
  for every tight filter $\xi$ in $\tsp E$, as proved in (5.8). Recalling that the tight groupoid $\Gt S$ of an inverse
semigroup $S$ may be constructed as the groupoid of germs for the canonical action of $S$ on the tight spectrum of its
idempotent semilattice, and considering an inverse semigroup homomorphism
  $$
  h:S\to T,
  $$
  we instead focus on constructing a covariant map between the dynamical systems from where one builds the tight
groupoids $\Gt S$ and $\Gt T$.

Denoting by $E$ and $F$ the respective idempotent semilattices of $S$ and $T$, this will obviously require a continuous
map between the tight spectra $\tsp E$ and $\tsp F$, and the only natural candidate is then $\chz$, where
  $$
  \hz:E\to F
  $$
  is the restriction of $h$ to the respective idempotent semilattices. Needless to say, we must require $\hz$ to be a
consonance in order for $\chz$ to be well defined. Making this hypothesis is however not such a big deal since the goal
is now to understand aspects of inverse semigroups that transcend their semilattices, so it is actually convenient not
to have to worry about $\hz$.

With the assumptions made so far, it is possible to define a covariant pair, and from there a groupoid homomorphism
  $$
  \ch: \Gt S\to\Gt T,
  $$
  and the problem is then to find conditions on $h$ for $\ch$ to be an isomorphism. After extending the notions of tight
injectivity/surjectivity to inverse semigroups (see 3.10), which in turn also requires the tight order to be extended
(see 3.7), we are able to prove that, as before, \pbitem$\ch$ is surjective if and only if $h$ is tightly surjective
(Theorem 6.7), and \pbitem$\ch$ is injective if and only if $h$ is tightly injective (Theorem 6.6).

\medskip Compared to the semilattice case, the play of words confronting injectivity with surjectivity and vice-versa is
now absent, of course due to the fact that we are using $\ch$ rather than $\hath$.

This justifies calling a homomorphisms of inverse semigroups a \emph{consonance}, whenever it is tightly injective and
tightly surjective, mimicking the corresponding concept for homomorphisms of semilattices. In particular, if $h$ is a
consonance between inverse semigroups, we show that the restriction $\hz$ of $h$ to the corresponding semilattices is a
consonance in the semilattice sense (see 7.2). This implies that $\ch$ is well defined, bringing us to one of our main
results, namely Theorem (7.3), reproduced here for the reader's convenience:

\Medskip{\bf Theorem}. If $h:S\to T$ is a consonance, then the map
  $$
  \ch: \Gt{S}\to\Gt{T}
  $$
  introduced in (6.3) is well defined, and it is an isomorphism of topological groupoids. Conversely, if $\hz$ is a
consonance (so that $\ch$ is well defined), and if $\ch$ is bijective, then $h$ itself is a consonance.

\medskip

Comparing this to Theorem (5.4), one cannot avoid to notice the conspicuous absence of any order related condition\fn{To
be sure, we are assuming that $\hz$ is a consonance, so the spectral orders on the tight spectra involved are indeed
playing an important role.} on $\ch$ in the last sentence of the above Theorem. Actually this is no surprise since,
contrary to tight spectra, tight groupoids are so far devoid of any order relation. However, like tight spectra, tight
groupoids may be equipped with a hugely important spectral order, without which we would have our hands tied.

Given an inverse semigroup $S$, the \emph{spectral order} on $\Gt S$ is defined as follows: given two elements
$\gamma,\delta\in\Gt{S}$, we say that $\gamma\leq\delta$ if \pbitem$d(\gamma)\leq d(\delta)$, relative to the spectral
order on $\tsp E$, and \pbitem there exists $s$ in $S$, such that both $\gamma$ and $\delta$ lie in $\Delta_s $,

\Medskip where $d$ denotes the domain, or source map of the groupoid $\Gt{S}$, and $\Delta_s$ refers to the fundamental
slice associated to $s$, defined by
  $$
  \Delta_s = \big\{[s,x]_\theta:x\in\dom(\theta_s)\big\}.
  $$
  See (2.14) for more details.

The second condition above, namely that both $\gamma$ and $\delta$ lie in $\Delta_s $, might look like it would never
lead to a transitive relation, but in fact the spectral order is a well behaved partial order on $\Gt{S}$, even
according to Ehresmann's definition of ordered groupoids (see Proposition (12.2)).

The importance of this new ingredient cannot be overstated, manifesting itself mainly through the derived notion of
up-sets: we say that a subset $U\subseteq\Gt S$ is an \emph{up-set} if, whenever $\gamma$ and $\delta$ are elements of
$\Gt S$ such that $\gamma$ lies in $U$, and $\delta\geq\gamma$, relative to the spectral order, then $\delta\in U$.

The first indication of the relevance of up-sets is that each $\Delta_s$ is such a thing. Moreover all open up-sets of
$\Gt S$ may be written as a union of slices of the form $\Delta_s$ (Proposition 8.6). In other words, the topology on
$\Gt S$ generated\fn{This topology was used in \cite[Section 3]{bussExelMeyer}.} by the $\Delta_s$ turns out to consist
precisely of the open up-sets of $\Gt S$.\fn{Recall that the standard topology on $\Gt S$ is not generated just by the
$\Delta_s$. In the simpler special case of the unit space $\tsp E$, the standard topology is generated by the $\Delta_e$
\emph{and} their complements. This also applies to $\Gt S$ in case it is Hausdorff.}

Secondly, products and inverses of up-sets are up-sets (Proposition 8.7). Therefore, denoting by $\Gt S\op$ the well
known\fn{The notation $G\op$ was introduced by Paterson in \cite[pg. xii]{paterson}.} inverse semigroup of slices, we
see that the compact \emph{up-slices} (slices which are up-sets) form an inverse sub-semigroup of $\Gt S\op$. We denote
this sub-semigroup by $\Cpl S$ and call it the \emph{tight envelope} of $S$ for reasons that we are about to explain.

Lest we forget to finish the discussion about the order properties of $\ch$ in the tight groupoid case, we must point
out that, according to Theorem (8.5), if $h:S\to T$ is a consonance then $\ch$ is an isomorphism of ordered topological
groupoids from $\Gt{S}$ to $\Gt{T}$. Combined with Theorem (7.3), this result then completely answers the question of
which inverse semigroup homomorphisms induce an isomorphism between the corresponding tight groupoids, preserving their
algebraic, topological and order structures. The answer is consonances.

Satisfactory as this answer is, no one should be under the illusion that it solves the whole problem, because there are
examples of consonant inverse semigroups $S$ and $T$ (actually even examples of semilattices) for which there are
consonances from $S$ to $T$, but none from $T$ to $S$.

In fact, the above use of the expression \emph{consonant inverse semigroups} is not really fair, not only because it has
not yet been defined, but mainly due to the misleading suggestion that it corresponds to an equivalence relation: to
make it very clear, the existence of a consonance from $S$ to $T$ is not an equivalence relation.

Thus, unless we want to abandon the notion of consonances altogether, the only solution is to generate an equivalence
relation out of this, saying that the inverse semigroups $S$ and $T$ are consonant if there exists a zig-zag of
consonances

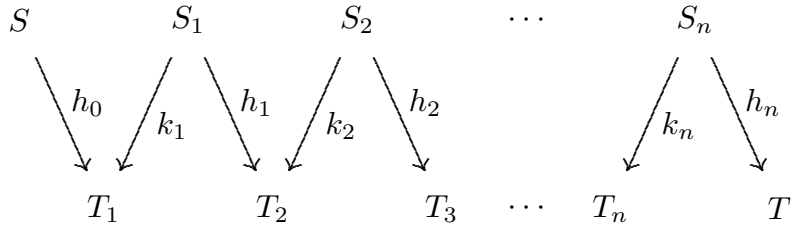
\begingroup\noindent\hfill\beginpicture\setcoordinatesystem units <0.025truecm, 0.025truecm> \setplotarea x from -30 to
430, y from 130 to -30 \put{\null} at -30 130 \put{\null} at -30 -30 \put{\null} at 430 130 \put{\null} at 430 -30
\put{$S$} at 0 100 \put{$T_1$} at 44.444 0 \arrow<0.11cm> [0.5, 1.8] from 8.935 79.896 to 35.509 20.104 \put{$h_0$} at
35.929 56.092 \put{$S_1$} at 88.889 100 \arrow<0.11cm> [0.5, 1.8] from 79.954 79.896 to 53.379 20.104 \put{$k_1$} at
80.374 43.908 \put{$T_2$} at 133.333 0 \arrow<0.11cm> [0.5, 1.8] from 97.824 79.896 to 124.398 20.104 \put{$h_1$} at
124.818 56.092 \put{$S_2$} at 177.778 100 \arrow<0.11cm> [0.5, 1.8] from 168.843 79.896 to 142.268 20.104 \put{$k_2$} at
169.263 43.908 \put{$T_3$} at 222.222 0 \arrow<0.11cm> [0.5, 1.8] from 186.713 79.896 to 213.287 20.104 \put{$h_2$} at
213.707 56.092 \put{$\cdots$} at 266.667 0 \put{$\cdots$} at 266.667 100 \put{$S_n$} at 355.556 100 \put{$T_n$} at
311.111 0 \arrow<0.11cm> [0.5, 1.8] from 346.621 79.896 to 320.046 20.104 \put{$k_n$} at 347.041 43.908 \put{$T$} at 400
0 \arrow<0.11cm> [0.5, 1.8] from 364.491 79.896 to 391.065 20.104 \put{$h_n$} at 391.485 56.092
\endpicture\hfill\null\endgroup

\noindent linking $S$ and $T$. If such a zig-zag exists we then say that $S$ and $T$ are \emph{consonant}.

This is now clearly an equivalence relation but the possibility of the length of this zig-zag going astray could make it
impractical to verify it in concrete situations.

In one of our main contributions, we show that this zig-zag is never too long. The result is Theorem (9.5), and it says:

\Medskip{\bf Theorem}. Let $S_1$ and $S_2$ be inverse semigroups. Then the following are equivalent:
  \item {(i)} $\Gt{S_1}$ and $\Gt{S_2}$ are isomorphic as ordered topological groupoids,
  \item {(ii)} $S_1$ and $S_2$ are consonant,
  \item {(iii)} there exists an inverse semigroup $T$, and consonances

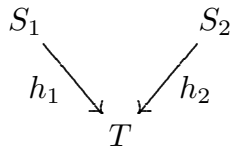
\begingroup\noindent\hfill\beginpicture\setcoordinatesystem units <0.025truecm, 0.025truecm> \setplotarea x from -30 to
130, y from 90 to -30 \put{\null} at -30 90 \put{\null} at -30 -30 \put{\null} at 130 90 \put{\null} at 130 -30
\put{$S_1$} at 0 60 \put{$T$} at 50 0 \arrow<0.11cm> [0.5, 1.8] from 9.603 48.477 to 40.397 11.523 \put{$h_1$} at 10.397
24.093 \put{$S_2$} at 100 60 \arrow<0.11cm> [0.5, 1.8] from 90.397 48.477 to 59.603 11.523 \put{$h_2$} at 89.603 24.093
\endpicture\hfill\null\endgroup

In other words, one needs never consider more than one ``zig'' and one ``zag''.

\medskip The purpose of this introductory section is not to discuss the technicalities behind proofs but the situation
demands an exception because all of the tools developed up to this point are about to converge in order to prove
(i)$\Rightarrow$(iii). The crux of the matter is the construction of the inverse semigroup $T$, and this is done as
follows: consider the diagram

\begingroup\noindent\hfill\beginpicture\setcoordinatesystem units <0.025truecm, 0.025truecm> \setplotarea x from -40 to
240, y from 130 to -40 \put{\null} at -40 130 \put{\null} at -40 -40 \put{\null} at 240 130 \put{\null} at 240 -40
\put{$S_1$} at 0 90 \put{$\Delta^1(S_1)$} at 110 90 \arrow<0.11cm> [0.5, 1.8] from 20 90 to 70 90 \put{$\Delta^1$} at 45
105 \put{$S_2$} at 0 0 \put{$\Delta^2(S_2)$} at 110 0 \arrow<0.11cm> [0.5, 1.8] from 20 0 to 70 0 \put{$\Delta^2$} at 45
-15 \put{$\Gt{S_1}\op$} at 200 90 \put{\rotatebox{0}{$\subseteq$}} at 150 90 \put{$\Gt{S_2}\op$} at 200 0
\put{\rotatebox{0}{$\subseteq$}} at 150 0 \arrow<0.11cm> [0.5, 1.8] from 200 70 to 200 20 \put{$\Phi$} at 215 45
\setdashes<1.5pt> \arrow<0.11cm> [0.5, 1.8] from 110 70 to 110 20 \put{$?$} at 125 45 \arrow<0.11cm> [0.5, 1.8] from
127.889 72.111 to 181.43 18.57 \put{$?$} at 165.266 55.947 \endpicture\hfill\null\endgroup

\noindent where $\Phi$ is the map on slices given by any isomorphism of ordered topological groupoids provided by (i),
and the $\Delta^i$ are the corresponding \emph{tight regular representations} introduced in \cite{tightrep}, namely the
maps sending each element $s$ to its fundamental slice $\Delta^i_s$.

Making each $\Delta^i$ surjective by replacing its codomain with its range $\Delta^i(S_i)$, it is not difficult to prove
that $\Delta^i$ becomes a consonance. It would be lovely should $\Phi$ send $\Delta^1(S_1)$ onto $\Delta^2(S_2)$,
because the choice of $T$ would then be obvious: any one of the isomorphic inverse semigroups $\Delta^i(S_i)$.

Needless to say, this is too much to hope for. The alternative is to map, say, $\Delta^1(S_1)$ into $\Gt{S_2}\op$ via
$\Phi$, and then consider the inverse sub-semigroup $T$ of $\Gt{S_2}\op$ generated by the union of
$\Phi\big(\Delta^1(S_1)\big)$ and $\Delta^2(S_2)$.

At this point, anyone acquainted with the hardships of working with things generated by arbitrary sets will perhaps feel
a chill run down their spine, but things can be put under control. Since each $\Delta^i_s$ is an up-set, and since
$\Phi$ preserves up-sets, we have that
  $$
  \Phi\big(\Delta^1(S_1)\big) \cup\Delta^2(S_2),
  $$
  consists of nothing but up-sets. Thanks to Proposition (8.7), which says that products and inverses of up-sets are
up-sets, we deduce that $T$ also consists of up-sets, and this is precisely the inverse semigroup that works.

The candidates for $h_1$ and $h_2$ are then clearly $h_1 = \Phi\circ\Delta^1$, and $h_2 = \Delta^2$, seen as maps into
$T$, so one must next prove that they are consonances.

For the case e.g.~of $h_2$, proving tight surjectivity essentially amounts to being able to write any $U\in T$, as
  $$
  U = \bigcup_{i\in I}\Delta_{s_i},
  $$
  for some family $\{s_i\}_{i\in I}$ of elements of $S_2$, which, as already mentioned, is precisely what (8.6) allows
us to do, as long as $U$ is an up-set, and this is why it is so important to make sure $T$ consists only of up-sets.

From this point on, things calm down considerably and it is not difficult to conclude the proof that $h_1$ and $h_2$ are
consonances.

\medskip Having already introduced the tight envelope $\Cpl S$ of a given inverse semigroup $S$, namely the
sub-semigroup of $\Gt S\op$ formed by the compact up-slices, and very much in the spirit of our discussion above, we may
see the tight regular representation of $S$ as a map
  $$
  \rho:S\to\Cpl S.
  $$
  In Proposition (8.11) we show that $\rho$ is a consonance, and for this reason we call $\rho$ the \emph{fundamental
consonance} of $S$. The above sketch of the construction of $T$ in the proof of Theorem (9.5) made it clear that $T$ is
closely related to $\Cpl S_2$, and indeed $T$ may be taken to be $\Cpl S_2$, itself, as proved in Corollary (10.2).

A further investigation of the tight envelope reveals that it has many special properties: it is flat (Proposition
8.10), meaning that its tight order coincides with the standard order, it contains the joins of all finite compatible
sets, and it is distributive in the sense that multiplication distributes over all existing finite joins (Proposition
10.7). As it turns out, according to Theorem (10.8), the tight envelope $\Cpl S$ is the only inverse semigroup with
these properties that is consonant to $S$.

\medskip Besides the tight envelope, there are other well known attempts at embedding an inverse semigroup $S$ into a
distributive one, notably Lawson's Booleanization $B(S)$ \cite{lawsonB}, and Lawson and Vdovina's tight completion
$T(S)$ \cite[Theorem 1.3] {lawsonV}. These may be defined as inverse semigroups formed by certain slices in Paterson's
universal groupoid, and in the tight groupoid of $S$, respectively. Both of these constructions produce a Boolean
inverse semigroup, and in particular their idempotent semilattices are Boolean. Incidentally a Boolean semilattice is a
Boolean algebra that has been stripped of everything but the meet operation. Due to \cite[Proposition 11]{actions},
every tight character on a Boolean semilattice is necessarily a Boolean algebra homomorphisms, which in turn implies
that every tight filter is an ultra-filter. As a result, if $E$ is a Boolean semilattice, every point of $\tsp E$ is
maximal, resulting in the the spectral order relation being trivial: it coincides with the equality relation.

Consequently, if $S$ is a given inverse semigroup with a nontrivial spectral order, then $S$ will never be consonant to
a Boolean inverse semigroup, due to the fact that their tight groupoids will never be order-isomorphic. Since the main
purpose of both the Booleanization and the tight completion is to embed an inverse semigroup into a distributive one,
the tight envelope may be seen as providing one more such embedding, with the advantage of being a consonance.

The last major endeavor we undertake is to obtain an abstract characterization of tight groupoids. Analyzing the
properties of the spectral order relation, we realize that they are surprisingly similar to the properties studied by
Ehresmann in his work on groupoids associated to inverse semigroup, except that Ehresmann's development often requires
his groupoids to be \emph{inductive}, which in turn is not a relevant condition for us.

Based on our analysis, we then go abstract and say that an ample groupoid equipped with a partial order relation is
\emph{tight-like} when it satisfies the fundamental properties found to hold on tight groupoids.

Our final main result is then to show that tight-like groupoids are indeed tight groupoids of inverse
semigroups. Precisely we show in Theorem (12.10) that the correspondence
  $$
  S\mapsto\Gt S
  $$
  is a bijection between the isomorphism classes of flat, distributive, inverse semigroups and the isomorphism classes
of tight-like groupoids. We further show that the inverse of this correspondence is given by
  $$
  G \mapsto\Up(G),
  $$
  where $\Up(G)$ is the \emph{fundamental inverse semigroup} of $G$, namely the semigroup consisting of all compact
up-slices of a given tight-like groupoid $G$.

\medskip A lot of the work we do here hinges on the delicate balance between tight filters and ultra-filters. Of course,
if all tight filters are ultra-filters, then the spectral order is trivial and all sets are up-sets, so the tight
envelope ends up coinciding with Lawson and Vdovina's tight completion. This situation, which is closely related to
Lawson's notion of compactable semilattices \cite{lawsonC}, is indeed very common, and it often presents itself when
\emph{finiteness} conditions are satisfied. For example, the canonical inverse semigroup within a Leavitt path algebra
of a finite graph has this property. Inverse semigroups occurring in the study of subshifts of finite type are also
examples. However, without any finiteness conditions the situation may be quite different, as illustrated e.g.~by
Cuntz-Krieger algebras of infinite matrices without the row-finiteness condition (see \cite[Proposition 7.9]{exelLaca}),
as well as all subshifts of infinite type, such as the even shift (see \cite{ExelSteinbergOne}, \cite{ExelSteinberg},
and the example at the end of \cite[Section 8]{DokExelShift}).

\medskip Our references for groupoids are \cite{Sims,paterson,actions,Renault}, while the reader may consult
\cite{paterson,lawsonBook} for the basic theory of inverse semigroups. Moreover our next section on preliminaries will
attempt to go over the basic topics that are most relevant to our goals.

\section 2 Preliminaries

\firstLine A \emph{semilattice} (sometimes also referred to as a \emph{meet semilattice}) is, by definition, a partially
ordered set $E$ (that is, a set equipped with a reflexive, anti-symmetric and transitive relation ``$\leq$'') such that,
for every $e$ and $f$ in $E$, the subset $\{g\in E : g\leq e, \ g\leq f\}$ admits a maximum element, usually denoted by
$e\wedge f$. This is sometimes referred to as the \emph{meet} of $e$ and $f$.

A semilattice is often profitably viewed as a multiplicative semigroup, with ``$\wedge$'' playing the role of the
multiplication operation. When the semigroup point of view is adopted, it is customary to denote $e\wedge f$ simply by
$ef$. A defining feature of this semigroup is that it is commutative and all elements are idempotent. In fact, it is
easy to prove that, if $E$ is a semigroup with these two properties, then the relation ``$\leq$'' defined on $E$ by
  $$
  e\leq f \reliff e = ef
  $$
  turns $E$ into a semilattice, and $e\wedge f = ef$, for all $e$ and $f$ in $E$.

In this work we shall take the semigroup point of view, meaning that we will adopt the notation $ef$ rather than
$e\wedge f$.

If a semilattice $E$ admits a minimum element, this element, which is necessarily unique, is often called the
\emph{zero} of $E$, and it is denoted by ``$0$''. In this case it is easy to see $0$ is also a \emph{multiplicative
zero} in the sense that $0e = e0 = 0$, for all $e$ in $E$.

\medskip{\bf In this work, the term \withfont{cmbxsl10}{semilattice} will always implicitly mean a semilattice with
zero.} \medskip

\fix From now on we fix a semilattice $E$.

\medskip If $F$ is another semilattice (with zero), we shall say that a function $h:E\to F$ is a \emph{homomorphism} if
$h(0) = 0$, and
  $$
  h(e_1e_2) = h(e_1)h(e_2),\for e_1,e_2\in E.
  $$
  The set of all homomorphisms from $E$ to $F$ will be denoted by
  $$
  \text{Hom}(E,F).
  $$

When $F = \{0,1\}$, with the obvious semilattice structure, a homomorphism from $E$ to $F$ will be called a
\emph{character}, provided it is not identically zero.

Given a character $\varphi:E\to\{0, 1\}$, the \emph{support} of $\varphi$, namely the set
  $$
  \xi= \{e\in E: \varphi(e) = 1\}
  $$
  is a \emph{filter}, in the sense that $\xi$ is a subset of $E$ satisfying:
  \item {(a)} $\xi$ is nonempty,
  \item {(b)} $0\notin\xi$,
  \item {(c)} $e,f\in\xi\relimply ef\in\xi$,
  \item {(d)} $f\geq e\in\xi\relimply f\in\xi$.

\Medskip Conversely, if $\xi\subseteq E$ is a filter, then the \emph{characteristic function} of $\xi$, namely the
function
  $$
  \varphi(e) = [e\in\xi], \for e\in E, \eqno{(2.1)}
  $$
  where the brackets stand for Boolean value, is a character. In other words, there is a natural bijective
correspondence between the set of all characters and the set of all filters.

\medskip Occasionally we will encounter subsets $\eta\subseteq E$ satisfying all filter axioms, except for (d) above. In
this case it is easy to see that
  $$
  \xi= \big\{f\in E: \exists e\in\eta,\ f\geq e\big\}
  $$
  is a genuine filter. We will then say that $\eta$ is a \emph{filter-base}, and that $\xi$ is the filter
\emph{generated} by $\eta$. An important example of a filter-base is as follows: suppose that $h$ is a homomorphism from
$E$ to $F$, such that $h(e)\neq0$, whenever $e\neq0$. If we are moreover given a filter $\xi$ in $E$, then clearly
  $$
  \eta= \{h(e): e\in\xi\}
  $$
  is a filter-base in $F$, so one might want to consider the filter generated by $\eta$.

\medskip The \emph{spectrum} of $E$, denoted $\hatE$, is the set of all characters on $E$. Observing that $\hatE$ is a
subset of $\{0,1\}^E$, we may equip $\hatE$ with the product topology (obviously viewing $\{0,1\}$ as discrete
space). This topology is sometimes also called the topology of \emph{pointwise convergence}. From now on, $\hatE$ will
always be seen as a topological space equipped with this topology.

\def\0{{\bf o}}

Recalling that the zero homomorphism, here denoted by ``$\0$'', has been explicitly excluded from being called a
character, we see that
  $$
  \hatE= \text{Hom}\big(E,\{0,1\}\big) \setminus\{\0\}.
  $$
  On the other hand, it is easy to see that $\text{Hom}\big(E,\{0,1\}\big)$ is a closed subset of $\{0,1\}^E$, hence
compact. Since $\hatE$ is obtained by removing a point from a compact set, we have that $\hatE$ is locally compact.

Given any $e$ in $E$, notice that
  $$
  V_e := \big\{\varphi\in\text{Hom}\big(E,\{0,1\}\big): \varphi(e) = 1\big\} \eqno{(2.2)}
  $$
  is clearly a clopen (closed and open) subset of $\text{Hom}\big(E,\{0,1\}\big)$, and hence compact. Moreover, since
$V_e\subseteq\hatE$ (because the condition ``$\varphi(e) = 1$'' automatically excludes $\0$) we see that $V_e$ is indeed
a compact open subset of $\hatE$.

\medskip When two elements $e$ and $f$ in $E$ satisfy $ef = 0$, we shall say that $e$ and $f$ are \emph{orthogonal}, and
we will express this fact by writing $e\perp f$. On the other hand, when $ef\neq0$, we shall say that $e$ and $f$
\emph{intercept}, expressing this as $e\Cap f$.

The following is part of \cite[Definition 11.5]{actions} and it will be of great interest to us.

\state 2.3. Definition. Given $f\in E$, and given a finite set $C\subseteq E$, we shall say that $C$ is a \emph{cover}
for $f$, if
  \item {(a)} $e\leq f$, for every $e\in C$,
  \item {(b)} whenever $g$ is a nonzero element in $E$, such that $g\leq f$, there exists some $e\in C$, such that
$e\Cap g$.

The definition of covers is inspired by the class of examples described below, which incidentally are the main examples
guiding our intuition.

\state 2.4. Definition. Let $X$ be a locally compact, totally disconnected Hausdorff space, and let $\boo$ be the
collection of all compact open subsets $X$, seen as a semilattice with respect to the order of inclusion. A
sub-semilattice $B\subseteq\boo$ (required to containing the zero element of $\boo$, namely the empty set) will be
called \emph{basic}, provided it forms a basis for the topology of $X$.

As already mentioned, basic semilattices are the inspiration for the notion of covers, and the reason is as follows: let
$B$ be a basic semilattice on $X$, pick $f$ in $B$, and let $C$ be a finite cover for $f$ in the above sense. It is then
easy to prove that
  $$
  c := \medcup_{e\in C}e \eqno{(2.5)}
  $$
  is a dense subset of $f$. However, since $C$ is was assumed finite, and since the members of $C$ are closed sets, $c$
is also closed, and hence $c$ coincides with $f$. In other words, $C$ is a cover for $f$ in the usual set-theoretic
sense of the word.

The catch here is that abstract semilattices have no operation similar to set-theoretic union, so, in general, it is
impossible to express (2.5). The definition of cover may then be seen as a workaround.

Returning to the general situation, if we are given a cover $C$ for some $f$ in $E$, as well as a character $\varphi$ on
$E$, it is clear that $\varphi(e)\leq\varphi(f)$, for all $e\in C$, so clearly
  $$
  \bigvee_{e\in C}\varphi(e) \leq\varphi(f). \eqno{(2.6)}
  $$
  Here we are using the fact that the semilattice $\{0,1\}$ is actually a Boolean algebra, and hence it also has a join
operation ``$\vee$''.

\state 2.7. Definition. Given a semilattice $E$, and given a character $\varphi$ on $E$, we shall say that $\varphi$ is
\emph{tight} provided (2.6) becomes an equality for every $f$ in $E$, and every cover $C$ for $f$. The subset of $\hatE$
formed by all tight characters is denoted by $\tsp E$.

We should warn the reader that the above definition of tight characters is an equivalent simplification of the general
definition of tightness for representations of semilattices in Boolean algebras \cite[Definition 11.6]{actions}, and, as
such, it does not apply beyond characters. See \cite[11.8]{actions} and \cite{covertojoin} for more on this.

If $\xi$ is a filter, and its characteristic function is a tight character, we shall say that $\xi$ is a \emph{tight
filter}.

On the other hand, an \emph{ultra-filter} is a filter which is not properly contained in any filter. So we shall say
that a character is an \emph{ultra-character} if its support is an ultra-filter.

Thus, under the correspondence between characters and filters, we have that tight characters correspond to tight
filters, while ultra-characters correspond do ultra-filters.

Throughout this work we will often switch between the character model and the filter model of spectra whenever it seems
convenient. In any case we will mostly use the Greek letters $\varphi$, $\chi$ and $\psi$ to denote characters, and
$\xi$, $\eta$ and $\zeta$ to denote filters, so this can be relied upon to determine which point of view is being used.

It is an important result that every ultra-character is tight, and in fact we have:

\state 2.8. Theorem. {\rm(c.f.~\cite[Theorem 12.9]{actions}).} Given a semilattice $E$, one has that $\tsp E$ coincides
with the closure in $\hatE$ of the set of all ultra-characters.

Among other things, the above says that $\tsp E$ is closed in $\hatE$, and hence the former is also locally
compact. Moreover, we have that, for every $e$ in $E$,
  $$
  \D_e := \big\{\varphi\in\tsp E : \varphi(e) = 1\big\} \eqno{(2.9)}
  $$
  is a compact open subset of $\tsp E$, because $\D_e = V_e\cap\tsp E$, where $V_e$ is as in (2.2). Switching from
characters to filters (the first of many times we will do this), we have that
  $$
  \D_e = \big\{\xi\in\tsp E : e\in\xi\big\}. \eqno{(2.10)}
  $$

\medskip A fact that we shall use repeatedly is in order:

\state 2.11. Proposition. Given a nonzero element $e$ in a semilattice $E$, there exists a tight character $\varphi$ on
$E$, such that $\varphi(e) = 1$.

\Proof Viewing the singleton $\{e\}$ as a filter-base, let $\xi$ be the filter generated by it. Employing Zorn's Lemma
we see that there is an ultra-filter $\zeta$ containing $\xi$. If $\varphi$ is the characteristic function of $\zeta$,
we then have that $\varphi$ is an ultra-character, hence tight, and clearly $\varphi(e) = 1$. \endProof

By definition, the topology of $\tsp E$ has a basis of open sets of the form
  $$
  U_{e_1, \ldots, e_n; f_1, \ldots, f_m} = \big\{\varphi\in\tsp E:\varphi(e_i) = 1,\ \varphi(f_j) = 0,\ 1\leq i\leq n, \
1\leq j\leq m\big\},
  $$
  where $n$ and $m$ are non negative integers, and $e_1, \ldots, e_n; f_1, \ldots, f_m$ are elements of $E$. Since the
zero homomorphism is not in $\tsp E$, we may clearly restrict to the situation in which $n>0$. Moreover, setting $e =
\medprod_{i = 1}^n e_i$, it is clear that
  $$
  \varphi(e) = 1\ \reliff\ \varphi(e_i) = 1,\for1\leq i\leq n,
  $$
  so the basic open set above coincides with
  $$
  U_{e; f_1, \ldots, f_m} = \big\{\varphi\in\tsp E:\varphi(e) = 1,\ \varphi(f_j) = 0,\ 1\leq j\leq
m\big\}. \eqno{(2.12)}
  $$

\state 2.13. Proposition. Let $\psi$ be an ultra-character. Then the collection of sets $\D_e$ defined in (2.9) (that
is, the sets of the form (2.12) with $m = 0$), where $e$ ranges in the support of $\psi$, forms a neighborhood basis for
$\psi$.

\Proof It is enough to prove that, given any open neighborhood of $\psi$ of the form $U_{e; f_1, \ldots, f_m}$, there
exists $e'$ in $E$ such that
  $$
  \psi\in U_{e'} \subseteq U_{e; f_1, \ldots, f_m}.
  $$
  Letting $\xi$ be the support of $\psi$, we have that $\xi$ is an ultra-filter. Moreover, $e$ lies in $\xi$, while the
$f_i$ do not. By \cite[12.3]{actions}, for each $i$, there exists some $e_i$ in $\xi$, such that $e_i\perp
f_i$. Observing that, for all characters $\varphi$, this gives
  $$
  \varphi(e_i)\varphi(f_i) = \varphi(e_if_i) = 0,
  $$
  we deduce that
  $$
  \varphi(e_i) = 1 \imply\varphi(f_i) = 0.
  $$
  It is then clear that $e' := e\medprod_{i = 1}^n e_i$ satisfies the desired conditions. \endProof

\bigskip Our interest in semilattices derives from the role they play in the theory of inverse semigroups. So let us now
briefly go over a few preliminaries on inverse semigroups.

If $S$ is a semigroup (that is, a set equipped with a binary associative operation), one says that $S$ is an
\emph{inverse semigroup}, provided for every $s$ in $S$, there exists a unique element $s^*$ in $S$, such that
  $$
  ss^*s = s, \and s^*ss^* = s^*.
  $$
  In this case the set
  $$
  E(S) = \{e\in S: e^2 = e\}
  $$
  is an abelian sub-semigroup, and moreover,
  $$
  e^* = e, \for e\in E(S).
  $$
  It is then easy to see that $E(S)$ is a semilattice, except that it does not necessarily have a zero element. When $S$
itself has a zero element, it is clear that zero is idempotent and hence it serves as a zero element for $E(S)$, as
well.

One usually refers to $E(S)$ as the \emph{idempotent semilattice} of $S$.

\medskip{\bf In this work, the expression \withfont{cmbxsl10}{inverse semigroup} will always implicitly mean an inverse
semigroup with zero.} \medskip

If $S$ and $T$ are inverse semigroups, we say that a map $h:S\to T$ is a homomorphism, provided $h(0) = 0$, and
  $$
  h(s_1s_2) = h(s_1)h(s_2),\for s_1,s_2\in S.
  $$

If $X$ is any set, the \emph{symmetric inverse semigroup} on $X$, denoted ${\cal I}(X)$, is the set of all bijective
maps between subsets of $X$, equipped with the operation of composition of partially defined maps. It is easy to prove
that ${\cal I}(X)$ is indeed an inverse semigroup, with the empty map playing the role of zero.

An \emph{action} of the inverse semigroup $S$ on a given locally compact Hausdorff topological space $X$ is a
homomorphism
  $$
  \theta:S\to{\cal I}(X),
  $$
  such that the domain and range of each $\theta_s$ is open, and
  $$
  \theta_s: \dom(\theta_s) \to\ran(\theta_s)
  $$
  is a homeomorphism. It is also required that $X$ coincide with the union of the domains of all of the
$\theta_s$. Observe that, according to our definition of homomorphism, $\theta_0$ is required to be the empty map,
namely the zero element of ${\cal I}(X)$.

An \emph{inverse semigroup dynamical system} is a triple $(\theta,S,X)$, where $S$ is an inverse semigroup, $X$ is a
locally compact Hausdorff topological space, and $\theta$ is an action of $S$ on $X$. See \cite[Section 4]{actions} for
more details.

Given an inverse semigroup dynamical system $(\theta,S,X)$, the corresponding \emph{groupoid of germs} (see \cite[page
140]{paterson} or \cite[Definition 4.6]{actions} for definitions), sometimes also called the \emph{transformation
groupoid}, will be denoted by $\GermGpd\theta$, or by
  $$
  S\ltimes_\theta X
  $$
  if we want to emphasize all ingredients involved. As it turns out, $\GermGpd\theta$ is always \'etale (but not
necessarily Hausdorff).

If $G$ is a groupoid, we shall denote its unit space by $G\ex0$, and we will say that the maps
  $$
  r:\gamma\in G\mapsto\gamma\gamma^{-1}\in G\ex0,\and d:\gamma\in G\mapsto\gamma^{-1}\gamma\in G\ex0,
  $$
  are the \emph{range} and \emph{domain} (or \emph{source}) maps, respectively. A subset $U\subseteq G$ will be called a
\emph{bisection} provided the restriction of $r$ and $d$ to $U$ are injective. If $G$ is moreover a topological
groupoid, and $U$ is an open bisection, we shall say that $U$ is a \emph{slice}.

En passant we notice that, in a \'etale topological groupoid, the collection of all slices is a basis for its topology.

Regarding the transformation groupoid $\GermGpd\theta$ mentioned above, for every $s\in S$, the set
  $$
  \Delta_s = \big\{[s,x]_\theta:x\in\dom(\theta_s)\big\}, \eqno{(2.14)}
  $$
  where $[s,x]_\theta$ is our notation for the germ of $s$ at $x$, is a slice, called the \emph{fundamental slice}
associated to $s$.

If $G$ is any \'etale groupoid, one usually denotes by $G\op$ the collection of all slices of $G$, and we recall that
$G\op$ has a natural structure of inverse semigroup, with multiplication
  $$
  UV = \big\{\alpha\beta: (\alpha,\beta)\in(U\times V)\cap G^{(2)}\big\}, \eqno{(2.15)}
  $$
  where $G\ex2$ is the set of \emph{composable} elements, namely
  $$
  G\ex2 = \{(\alpha,\beta): d(\alpha) = r(\beta)\}.
  $$

In the context of the groupoid $\GermGpd\theta$ discussed above, we have by \cite[Proposition 7.4]{actions} that
  $$
  \Delta_s \Delta_t = \Delta_{st},\for s,t\in S,
  $$
  so the map
  $$
  \Delta:s\in S\mapsto\Delta_s \in\GermGpd\alpha\op\eqno{(2.16)}
  $$
  is a semigroup homomorphism. Since we are assuming that $\theta_0$ is the empty map, it follows that $\Delta_0$ is the
empty slice, i.e.~the zero element of $\GermGpd\alpha\op$. In other words, $\Delta$ respects the requirement of sending
zero to zero.

\bigskip Fixing an inverse semigroup $S$, and letting $E = E(S)$ be its idempotent semilattice, we shall now describe a
canonical action of $S$ on $\tsp E$, which will in turn lead to the notion of the tight groupoid of $S$.

Keeping (2.9) in mind, given $s\in S$, and given $\varphi\in\D_{s^*s}$, the map
  $$
  \theta_s(\varphi):e\in E\mapsto\varphi(s^*es)\in\{0, 1\} \eqno{(2.17)}
  $$
  is clearly a character. By \cite[Proposition 12.8]{actions} one has that $\theta_s(\varphi)$ is tight, and since
$\theta_s(\varphi)(ss^*) = 1$, it follows that $\theta_s(\varphi)$ lies in $\D_{ss^*}$. This in turn defines a map
  $$
  \theta_s : \D_{s^*s}\to\D_{ss^*},
  $$
  which may be proven to be a homeomorphism. According to \cite[Sections 10 and 12]{actions} we have that $\theta$
defines an action of $S$ on $\tsp E$, so that $(\theta,S,\tsp E)$ is an inverse semigroup dynamical system.

Should we have the filter version of $\tsp E$ in mind, it is easy to prove that, for every tight filter
$\xi\in\D_{s^*s}$, one has that $\theta_s(\xi)$ is the filter generated by the filter-base
  $$
  \eta= \{ses^*: e\in\xi\}. \eqno{(2.18)}
  $$

The corresponding groupoid of germs, denoted by $\Gt{S}$, is then termed the \emph{tight groupoid} of $S$.

\state 2.19. Definition. Given an inverse semigroup $S$, the \emph{tight regular representation} of $S$ is the
homomorphism
  $$
  \Delta:s\in S\mapsto\Delta_s \in\Gt S\op
  $$
  described in (2.16), relative to the dynamical system $(\theta,S,\tsp E)$.

At this point it is perhaps worth remarking that, if a given semilattice $E$ is seen as an inverse semigroup (which it
actually is), then the canonical action of $E$ on $\tsp E$ is trivial, in the sense that all maps
  $$
  \theta_e : \D_{e}\to\D_{e}
  $$
  are the identity. Consequently the tight groupoid $\Gt E$ is naturally isomorphic to the topological space $\tsp E$, a
groupoid in which all elements are units. In turn the semigroup of slices $\Gt E\op$ coincides with $\O(\tsp E)$, namely
the semilattice of all open subsets of $\tsp E$, incidentally the same thing as the the topology of $\tsp E$. For each
$e$ in $E$, the corresponding fundamental slice, as defined in (2.14), turns out to be $\D_e$, while the tight regular
representation of $E$ becomes the homomorphism
  $$
  e\in E\mapsto\D_e\in\O(\tsp E). \eqno{(2.20)}
  $$

Closing this section, let us present a few preliminaries on covariant maps. For this let us fix two inverse semigroup
dynamical systems $(\alpha,S,X)$ and $(\beta,T,Y)$.

\state 2.21. Definition. Given a homomorphism $h:S\to T$, and given a continuous mapping $f:X\to Y$, we shall say that
the pair $(h,f)$ is \emph{covariant} if, for every $s$ in $S$,
  $$
  f\big(\dom(\alpha_s)\big)\subseteq\dom\big(\beta_{h(s)}\big), \eqno{(2.21.1)}
  $$
  and
  $$
  f\big(\alpha_s(x)\big) = \beta_{h(s)}\big(f(x)\big), \for x\in\dom(\alpha_s).
  $$

Covariant pairs induce groupoid homomorphisms on the corresponding groupoids of germs, as we shall now show.

\state 2.22. Proposition. Given a covariant pair $(h,f)$, there exists a continuous homomorphism of groupoids
(a.k.a.~functor)
  $$
  h{\times}f:\GermGpd\alpha\to\GermGpd\beta,
  $$
  such that
  $$
  (h{\times}f)\big([s,x]_\alpha\big) = [h(s),f(x)]_\beta,
  $$
  for every germ $[s,x]_\alpha$ in $\GermGpd\alpha$.

\Proof Suppose that for every $i = 1, 2$, we are given $s_i$ in $S$, and $x_i\in\dom(\alpha_{s_i})$, such that
$[s_1,x_1]_\alpha= [s_2,x_2]_\alpha$. We then have that $x_1 = x_2$, and there exists $e\in E(S_1)$, such that $x_1$
belongs to $\dom(\alpha_e)$, and $s_1e = s_2e$. Therefore \bitem$f := h(e)\in E(T)$, \bitem$f(x_1)\in\dom(\beta_{f})$,
and \bitem$h(s_1)f = h(s_2)f$, \Medskip so we see that
  $$
  \big[h(s_1),f(x_1)\big]_\beta= \big[h(s_2),f(x_1)\big]_\beta,
  $$
  whence the correspondence
  $$
  [s,x]_\alpha\mapsto\big[h(s),f(x)\big]_\beta
  $$
  is well defined, providing a map
  $$
  h{\times}f:\GermGpd\alpha\to\GermGpd\beta,
  $$
  which may be easily shown to satisfy the required conditions. \endProof

Bearing in mind the definition of the inverse semigroup of slices, given near (2.15), a groupoid homomorphism, such as
the one in (2.22), raises the question as to its behavior in face of the corresponding inverse semigroups of slices.

\state 2.23. Lemma. Let $G$ and $H$ be \'etale groupoids, let $\varphi:G\to H$ be a continuous homomorphism, and let
  $$
  \varphi\ex0 :G\ex0 \to H\ex0
  $$
  be obtained by restricting $\varphi$ to the corresponding unit spaces. Assuming that $\varphi\ex0$ is an open map, we
have that $\varphi$ is also an open map. If, in addition, $\varphi\ex0$ is one-to-one, then $\varphi$ sends slices of
$G$ to slices of $H$, and the corresponding map
  $$
  \Phi:U\subseteq G\op\mapsto\varphi(U)\in H\op
  $$
  is an inverse semigroup homomorphism.

\Proof Choosing any open subset $U\subseteq G$, we will prove that $\varphi(U)$ is open. Given $\gamma$ in $U$, let $V$
be any slice containing $\varphi(\gamma)$. Using that $\varphi$ is continuous, we then pick an open subset $W\subseteq
U$, containing $\gamma$, such that $\varphi(W)\subseteq V$. Observing the obvious fact that
  $$
  \d(\varphi(\gamma)) = \varphi\ex0(\d(\gamma)), \for\gamma\in G,
  $$
  where, by abuse of language, we are using $\d$ to refer to both the source map on $G$, and the source map on $H$, we
see that the following diagram is commutative:

\begingroup\noindent\hfill\beginpicture\setcoordinatesystem units <0.025truecm, -0.02truecm> \setplotarea x from -30 to
160, y from -30 to 130 \put{\null} at -30 -30 \put{\null} at -30 130 \put{\null} at 160 -30 \put{\null} at 160 130
\put{$W$} at 0 0 \put{$V$} at 130 0 \arrow<0.11cm> [0.5, 1.8] from 19.5 0 to 110.5 0 \put{$\varphi$} at 65 -15
\put{$\d(W)$} at 0 100 \put{$\d(V)$} at 130 100 \arrow<0.11cm> [0.5, 1.8] from 28.6 100 to 101.4 100 \put{$\varphi\ex0$}
at 65 85 \arrow<0.11cm> [0.5, 1.8] from 0 25 to 0 75 \put{$\d$} at -10 50 \arrow<0.11cm> [0.5, 1.8] from 130 25 to 130
75 \put{$\d$} at 140 50 \endpicture\hfill\null\endgroup

Since $\varphi\ex0\circ\d$ is an open map, we deduce that $\varphi\ex0(\d(W))$, and hence also $\d(\varphi(W))$, is open
in $\d(V)$. As $V$ is a slice, we have that $\d$ is a homeomorphism from $V$ to $\d(V)$, so it follows that $\varphi(W)$
is open in $V$, and hence open in $H$. Finally, observing that
  $$
  \varphi(\gamma)\in\varphi(W)\subseteq\varphi(U),
  $$
  we deduce that $\varphi(U)$ is open.

Assuming that $\varphi\ex0$ is moreover one-to-one, let $U\subseteq G$ be a slice. Since $\d$ is injective on $U$, we
see that $\varphi\ex0\circ\d$ is also injective on $U$, so the same is true for $\d\circ\varphi$. This implies that $\d$
is injective on $\varphi(U)$. By the same reasoning one proves that the range map $r$ is injective on $\varphi(U)$, and
hence that $U$ is a bisection. Having already proved that $\varphi(U)$ is open, we have shown that $\varphi(U)$ is a
slice, as desired.

To conclude we must show that
  $$
  \Phi(UV) = \Phi(U)\Phi(V),
  $$
  whenever $U, V\in G\op$. While the inclusion ``$\subseteq$'' is obvious, the reverse one requires a little
argument. For this, choose $\gamma_1$ and $\gamma_2$ in $G$, and suppose that
  $$
  \big(\varphi(\gamma_1),\varphi(\gamma_2)\big)\in\big(\Phi(U)\times\Phi(V)\big)\cap H\ex2.
  $$
  Therefore $\varphi(\gamma_1)\varphi(\gamma_2)$ is a typical element of $\Phi(U)\Phi(V)$, which we in turn would like
to prove lies in $\Phi(UV)$. Observing that $\d(\varphi(\gamma_1)) = r(\varphi(\gamma_2))$, or, equivalently,
  $$
  \varphi\ex0(\d(\gamma_1)) = \varphi\ex0(r(\gamma_2)),
  $$
  we see that $\d(\gamma_1) = r(\gamma_2)$, meaning that $(\gamma_1,\gamma_2)\in G\ex2$. Therefore $\gamma_1\gamma_2\in
UV$, and
  $$
  \varphi(\gamma_1)\varphi(\gamma_2) = \varphi(\gamma_1\gamma_2) \in\varphi(UV),
  $$
  proving the reverse inclusion. \endProof

We shall now consider covariant pairs $(h, f)$ in which $f$ is a surjective map, but we shall actually need a bit more
in the sense that surjectivity extends to domains.

\state 2.24. Definition. A covariant pair $(h, f)$ will be called a \emph{covariant epimorphism}, if $f$ is surjective
and, in addition,
  $$
  f\big(\dom(\alpha_s)\big) = \dom\big(\beta_{h(s)}\big), \for s\in S.
  $$
  (Compare 2.21.1).

The following is easily proved by simply unraveling the definitions involved.

\state 2.25. Proposition. Given inverse semigroup dynamical systems $(\alpha,S,X)$ and $(\beta,T,Y)$, and given a
covariant pair $(h,f)$, one has that
  $$
  (h{\times}f)(\Delta^\alpha_s) \subseteq\Delta^\beta_{h(s)}, \for s\in S,
  $$
  where we denote the fundamental slices described in (2.14) by $\Delta^\alpha_s$ and $\Delta^\beta_t$ relative to the
actions $\alpha$ and $\beta$, respectively. If $(h, f)$ is moreover a covariant epimorphism, then
  $$
  (h{\times}f)(\Delta^\alpha_s) = \Delta^\beta_{h(s)}, \for s\in S.
  $$

\section 3 Tight order

\firstLine If $e$ and $f$ are sets with $e\not\subseteq f$, then there necessarily exists a nonempty set $g$ such that
  $$
  f\cap g = \varnothing, \and g\subseteq e,
  $$
  namely $g = e\setminus f$. One could then say that $g$ is a \emph{witness} to the fact that $e\not\subseteq
f$. However, if $e$ and $f$ are elements of a semilattice, and $e\not\leq f$, it is not necessarily true that a nonzero
element $g$ exists such that
  $$
  f\perp g,\and g\leq e.
  $$
  In other words, it is possible that, although $e\not\leq f$, there is no witness denying that $e\leq f$.

An example of this situation is $e = 2$, and $f = 1$, in the semilattice $\{0,1,2\}$, with the order induced by the
natural numbers.

\state 3.1. Definition. Given a semilattice $E$, and given $e$ and $f$ in $E$, we say that $e$ is \emph{tightly below}
$f$, in symbols $e \menor f$, if there is no nonzero $g$ in $E$, such that $ f\perp g \leq e. $ This will be called the
\emph{tight order relation} on $E$.

In other words, $e\menor f$ if and only if there is no witness denying that $e\leq f$, regardless of whether or not the
latter is true. If $e\leq f$, then clearly $e\menor f$, but the converse may fail, as seen in the example above.

\fix From now on we fix a semilattice $E$.

\state 3.2. Lemma. Given $e$ and $f$ in $E$, the following are equivalent:
  \item {(i)} $e\menor f$,
  \item {(ii)} for every tight character $\varphi$ on $E$, one has that $\varphi(e)\leq\varphi(f)$,
  \item {(iii)} for every tight filter $\xi$ on $E$, one has that $e\in\xi\relimply f\in\xi$,
  \item {(iv)} for every $g$ in $E$, one has that $g\perp f\imply g\perp e$,
  \item {(v)} for every nonzero $g$ in $E$, such that $g\leq e$, one has that $g\Cap f$.

\Proof \proofImply {i}{ii} Assuming (i), suppose by contradiction that $\varphi$ is a tight character on $E$, such that
$\varphi(e)\not\leq\varphi(f)$, so necessarily
  $$
  \varphi(e) = 1,\and\varphi(f) = 0.
  $$
  Then the singleton $\{ef\}$ cannot possibly be a cover for $e$, or else
  $$
  \varphi(e) = \bigvee_{c\in\{ef\}}\varphi(c) = \varphi(ef) = \varphi(e)\varphi(f) = \varphi(f).
  $$
  Therefore there exists a nonzero $g\leq e$, such that $g\perp ef$. The latter then gives
  $$
  0 = gef = (ge)f = gf,
  $$
  so $g\perp f$, contradicting the fact that $e\menor f$.

\Medskip(ii$\Leftrightarrow$iii)\enspace Obvious.

\Medskip\proofImply{ii}{iv} Again arguing by contradiction, suppose that there exists $g$ in $E$, such that $g\perp f$,
but $ge\neq0$. By (2.11) we may choose a tight character $\varphi$ such that $\varphi(ge) = 1$, so
  $$
  1 = \varphi(ge) = \varphi(g)\varphi(e),
  $$
  whence $\varphi(g) = \varphi(e) = 1$, while
  $$
  0 = \varphi(gf) = \varphi(g) \varphi(f) = \varphi(f) \geq\varphi(e) = 1,
  $$
  a contradiction.

\Medskip\proofImply{iv}{i} Letting $g\in E$ be such that $f\perp g\leq e$, we have that (iv) implies that $g\perp e$, so
  $$
  g = ge = 0.
  $$

We leave the easy proof of the equivalence of (i) and (v) to the reader. \endProof

The tight order relation is one half of a congruence studied by Lawson and Lenz in \cite[Lemma 4.8]{LawsonLenz}. We will
shortly consider the other half.

Notice that the above result implies that ``$\menor$'' is a reflexive and transitive relation, although it might not be
anti-symmetric.

\state 3.3. Definition. Given $e$ and $f$ in $E$, we shall say that $e$ and $f$ are \emph{tightly equivalent}, in
symbols $e\equiva f$, if $e\menor f$, and $f\menor e$.

By (3.2), it is clear that $e\equiva f$ if and only if $\varphi(e) = \varphi(f)$, for every tight character $\varphi$.

\state 3.4. Proposition. Denote the tight regular representation of $E$ by
  $$
  D:e\in E\mapsto\D_e\in\O(\tsp E),
  $$
  as in (2.20). Then, for all $e$ and $f$ in $E$, one has that
  $$
  e\menor f\iff D(e)\leq D(f),
  $$
  and consequently also
  $$
  e\equiva f\iff D(e) = D(f).
  $$
  Therefore $E/{\equiva}$ is a semilattice, naturally isomorphic to the range of $D$. Moreover, for $e$ and $f$ in
$E/{\equiva}$, one has that $e\menor f$ if and only if $e\leq f$, and hence also $e\equiva f$ if and only if $e = f$.

\Proof Given $e$ and $f$ in $E$, notice that the fact that $D(e)\leq D(f)$, that is, $\D_e\subseteq\D_f$, is equivalent
to saying that, for every tight character $\varphi$,
  $$
  \varphi(e) = 1 \relimply\varphi(f) = 1,
  $$
  which in turn is the same as saying that $\varphi(e)\leq\varphi(f)$. That this is equivalent to $e\menor f$ is then a
consequence of (3.2).

In order to prove the last sentence in the statement, we look at the range of $D$ instead of $E/{\equiva}$, so it is
enough to show that, for every $e$ and $f$ in $E$, if $\D_e\not\subseteq\D_f$, then there exists $g$ in $E$ such that
$\D_g$ is not empty, and
  $$
  \D_f\cap\D_g = \varnothing, \and\D_g\subseteq\D_e.
  $$
  To prove this we observe that $\D_e\setminus\D_f$ is a nonempty open set, so it contains an ultra-character $\varphi$
by (2.8). Using (2.13), we see that there exists some $g$ in the support of $\varphi$, such that
  $$
  \varphi\in\D_g\subseteq\D_e\setminus\D_f,
  $$
  concluding the proof. \endProof

\state 3.5. Definition. Suppose that we are given semilattices $E$ and $F$, as well as a homomorphism $h:E\to F$. Then:
  \item {(a)} a subset $R\subseteq F$ is called \emph{large} if, for every $f$ in $F$, there exists a finite set
$C\subseteq R$, such that $e\menor f$, for all $e\in C$, and $\{fe: e\in C\}$ is a cover for $f$.
  \item {(b)} $h$ is said to be \emph{tightly surjective} if its range is a large subset of $F$,
  \item {(c)} $h$ is said to be \emph{tightly injective} if, for all $e_1$ and $e_2$ in $E$, one has that
  $$
  h(e_1)\equiva h(e_2) \relimply e_1\equiva e_2,
  $$
  \item {(d)} the \emph{kernel} of $h$ is defined to be the set
  $$
  \ker(h) = \big\{e\in E: h(e) = 0\big\}.
  $$

A few remarks about these concepts are in order. First of all, ignoring the difference between the tight order and the
usual order in $F$, one gets an incorrect, albeit enlightening interpretation of (3.5.a), saying that every $f$ in $F$
is covered by a finite set contained in $R$.

Regarding (3.5.d), recall that in many areas of Mathematics where there is a meaningful notion of kernel, such as in
Group Theory or Linear Algebra, a morphism is injective if and only if its kernel is trivial. The same is however not
true in the realm of semilattices, as illustrated by the character
  $$
  \varphi:\{0,1,2\}\to\{0,1\},
  $$
  sending both $1$ and $2$ to $1$.

Nevertheless this general principle partly survives in the form of the following result:

\state 3.6. Proposition. Let $E$ and $F$ be semilattices and let $h:E\to F$ be a homomorphism. Then the following are
equivalent:
  \item {(i)} $\ker(h) = \{0\}$,
  \item {(ii)} $h$ is tightly injective,
  \item {(iii)} for every $e$ and $f$ in $E$, one has that $ h(e)\menor h(f) \Rightarrow e\menor f. $

\Proof \proofImply {i}{iii} Suppose that $\ker(h) = \{0\}$, and suppose, by contradiction, that there are elements $e$
and $f$ in $E$ such that $h(e) \menor h(f)$, and yet $e\not\menor f$. Therefore there is a nonzero $g$ in $E$, such that
$ f\perp g \leq e. $ Clearly this implies that
  $$
  h(f)\perp h(g) \leq h(e),
  $$
  and since $h(g)\neq0$, by hypothesis, we see that $h(e)\not\menor h(f)$, a contradiction.

\Medskip\proofImply{iii}{ii} If $h(e)\equiva h(f)$, we have by definition that $h(e)\menor h(f)$, and $h(f)\menor h(e)$,
so (iii) implies that $e\menor f$ and $f\menor e$.

\Medskip\proofImply{ii}{i} If $e\neq0$, we clearly have that $e\not\menor0$ (consider taking $g = e$), so
$e\not\equiva0$, whence (ii) gives
  $$
  h(e)\not\equiva h(0) = 0,
  $$
  so, in particular, $h(e)\neq0$. This shows that $\ker(h) = \{0\}$. \endProof

\bigskip The notion of the tight order relation ``$\menor$'' on semilattices extends to inverse semigroups, but not via
a naive application of its definition, as this would not yield a useful concept. The appropriate generalization is as
follows:

\state 3.7. Definition. Given an inverse semigroup $S$, and given $s,t\in S$, we shall say that:
  \item {(a)} $s\menor t$, when there exists a finite cover\fn{Whenever we speak of a cover in the context of an inverse
semigroup $S$, we shall always view it within the corresponding idempotent semilattice $E(S)$.} $C$ for $s^*s$, such
that $se = te$, for every $e\in C$.
  \item {(b)} $s\equiva t$, when $s\menor t$, and $t\menor s$. \Medskip The relation ``$\menor$'' will be called the
\emph{tight order relation} on $S$.

At first sight this might present a conflict with our previous interpretation of the symbol ``$\menor$'', as defined in
(3.1), but this conflict may be dispelled as follows:

\state 3.8. Proposition. If $e$ and $f$ are idempotent elements in an inverse semigroup $S$, then $e\menor f$ within
$E(S)$, in the sense of Definition (3.1), if and only if $e\menor f$ in the sense of Definition (3.7).

\Proof During this proof, and before we resolve the conflict mentioned above, we will denote the relation introduced in
(3.7) by ``$\menor_1$''.

Assuming that $e\menor_1 f$, suppose by contradiction that there is a nonzero idempotent $g$ such that $f\perp g \leq
e$. Then by hypothesis there exists a finite cover $C$ for $e$, such that $ec = fc$, for every $c\in C$. By the
definition of cover, there is some $c\in C$, such that $g\Cap c$. Therefore
  $$
  0\neq gc = gec = gfc = 0,
  $$
  a contradiction.

Now suppose that $e\menor f$. We then claim that the singleton $C = \{ef\}$ is a cover for $e$. In fact, if $g\leq e$ is
nonzero, then $g$ cannot possibly be orthogonal to $f$, or else $f\perp g \leq e$, contradicting our assumption that
$e\menor f$. Thus
  $$
  0 \neq gf = gef,
  $$
  proving that $g\Cap ef$, and hence proving our claim that $C$ is a cover for $e$. Finally, if $c$ is the single
element of $C$, then
  $$
  ec = eef = fef = fc,
  $$
  proving that $e\menor_1 f$. \endProof

\state 3.9. Proposition. Given an inverse semigroup $S$, and given $s,t\in S$, one has that
  \item {(i)} $s\menor t \relimply s^*s\menor t^*t$,
  \item {(ii)} $s\equiva t \relimply s^*s\equiva t^*t$.

\Proof Assuming that $s\menor t$, take a finite cover $C$ for $s^*s$, such that $se = te$, for every $e\in C$. If, by
contradicting, there is a nonzero idempotent $p$ such that $t^*t\perp p \leq s^*s$, we may find some $e\in C$, such that
$e\Cap p$, and hence that
  $$
  0\neq ep = s^*sep = s^*tep = s^*tt^*tpe = 0,
  $$
  a contradiction. \endProof

Having extended the notion of tight order to inverse semigroups, we may generalize Definition (3.5).

\state 3.10. Definition. Suppose that $S$ and $T$ are inverse semigroups and $h:S\to T$ is a homomorphism. Then
  \item {(a)} a subset $R\subseteq T$ is called \emph{large} if, for every $t$ in $T$, there exists a finite set
$C\subseteq R$, such that $s\menor t$, for all $s\in C$, and $\{s^*st^*t:s\in C\}$ is a cover for $t^*t$.
  \item {(b)} $h$ is said to be \emph{tightly surjective} if its range is a large subset of\/ $T$.
  \item {(c)} $h$ is said to be \emph{tightly injective} if, for all $s_1$ and $s_2$ in $S$, one has that
  $$
  h(s_1)\equiva h(s_2) \relimply s_1\equiva s_2.
  $$

It is easy to see that, if the above definitions are applied to a homomorphism between semilattices, seen as inverse
semigroups, then they are equivalent to the homonimous concepts defined in (3.5). In other words, the definition above
causes no conflicts with previously defined terms.

\medskip Regarding tight injectivity, it would be interesting to decide whether or not (3.6) generalizes to inverse
semigroups, with some appropriate definition of kernel.

\medskip We now present the inverse semigroup generalization of Proposition (3.4), during which it is useful to keep in
mind that, for the case of the tight groupoid $\Gt{S}$, the fundamental slices described in (2.14) become
  $$
  \Delta_s = \big\{[s,\xi]_\theta: \xi\in\D_{s^*s}\big\} \explica{(2.10)}{=} \big\{[s,\xi]_\theta: \xi\in\tsp E, \
s^*s\in\xi\big\}. \eqno{(3.11)}
  $$

\state 3.12. Lemma. Given an inverse semigroup $S$, and given $s$ and $t$ in $S$, one has that:
  \item {(i)} if $s\menor t$, and $\xi$ is a tight filter such that $s^*s\in\xi$, then $t^*t\in\xi$, and
$[s,\xi]_\theta= [t,\xi]_\theta$,
  \item {(ii)} $s\menor t \reliff\Delta_s \subseteq\Delta_t$.

\Proof Regarding the first statement in (i), we have that $s^*s\menor t^*t$, by (3.9), so (3.2) implies that
$t^*t\in\xi$.

Let $C$ be a finite cover for $s^*s$, such that $se = te$, for every $e\in C$. Since $\xi$ is a tight filter, we have
that $e\in\xi$, for some $e$ in $C$, whence
  $$
  [s, \xi]_\theta= [t, \xi]_\theta.
  $$
  This proves (i), and clearly also the forward implication of (ii).

Assuming now that $\Delta_s \subseteq\Delta_t$, and discarding the trivial case in which $s = 0$, pick any $\xi$ in
$\D_{s^*s}$, so that $[s,\xi]_\theta\in\Delta_t$, by hypothesis. Noticing that $[t, \xi]_\theta$ is the only element in
$\Delta_t$ whose source is $\xi$, we deduce that $[s, \xi]_\theta= [t, \xi]_\theta$, so there exists some $e\in\xi$,
such that $se = te$, and we can clearly assume that $e\leq s^*s$. This shows that
  $$
  \xi\in\D_e\subseteq\D_{s^*s},
  $$
  and hence that the collection of the sets $\D_e$, with $e\leq s^*s$, and $se = te$, forms a cover for $\D_{s^*s}$. By
compactness of $\D_{s^*s}$ we may therefore pick a finite set $ \{e_1,e_2,\ldots, e_n\}\subseteq E, $ such that $e_i\leq
s^*s$, and $se_i = te_i$, for all $i$, and
  $$
  \D_{s^*s} = \bigcup_{i = 1}^n \D_{e_i}.
  $$
  We will next prove that $\{e_1,e_2,\ldots, e_n\}$ is a cover for $s^*s$. Indeed, arguing by contradiction, assume that
there exists some nonzero $f\leq s^*s$, such that $ f\perp e_i, $ for all $i$. Choosing any tight filter $\xi$
containing $f$, we deduce that $\xi$ lies in $\D_{s^*s}$, but not in any of the $\D_{e_i}$, a contradiction. \endProof

\state 3.13. Corollary. Given an inverse semigroup $S$, and given $s$ and $t$ in S, one has that
  $$
  s\equiva t \reliff\Delta_s = \Delta_t.
  $$
  Consequently the map
  $$
  \Pi: S/{\equiva} \to\Gt S\op,
  $$
  obtained by factoring $\Delta$ through ``$\kern2pt \equiva$'', is an isomorphism onto its range.

The following are simple technical consequences of the above Corollary for everyday use.

\state 3.14. Lemma. Given an inverse semigroup $S$, suppose that we are given $e, f, r, s,t\in S$, with $e$ and $f$
idempotent. Then
  \item {(i)} if $sf\equiva tf$, and $e\menor f$, then $se\equiva te$,
  \item {(ii)} if $s\menor t$, then $sr\menor tr$.

\Proof (i)\enspace Since $e\menor f$, we have that $\Delta_e\subseteq\Delta_f$, so
  $$
  \Delta_{se} = \Delta_s \Delta_e = \Delta_s \Delta_f\Delta_e = \Delta_{sf}\Delta_e = \Delta_{tf}\Delta_e =
\Delta_t\Delta_f\Delta_e = \Delta_t\Delta_e = \Delta_{te},
  $$
  so $se\equiva te$.

\Medskip(ii)\enspace If $s\menor t$, then $\Delta_s \subseteq\Delta_t$, so
  $$
  \Delta_{sr} = \Delta_s \Delta_r \subseteq\Delta_t\Delta_r = \Delta_{tr},
  $$
  so $sr\menor tr$. \endProof

The above are just a few basic facts that we will actually use below, but, of course, there are lots of other similar
elementary properties one may prove with the same technique.

\section 4 Tight semilattice homomorphisms and dual maps

\firstLine Our goal here is to introduce a dual, or transposed, map $\hath$ from $\tsp F$ to $\tsp E$, for a given
semilattice homomorphism $h:E\to F$. However, the most naive attempt at this meets the obstruction that, given a
character $\psi\in\tsp F$, the composition $\psi\circ h$ may fail to be tight. Even worse, $\psi\circ h$ may vanish, and
hence it would not even be character.

\fix Throughout this section, we shall fix semilattices $E$ and $F$, as well as a homomorphism
  $$
  h:E\to F.
  $$

\state 4.1. Proposition. The following are equivalent:
  \item {(i)} for every tight character $\psi$ on $F$, one has that $\psi\circ h$ is nonzero (hence a character),
  \item {(ii)} for every $f$ in $F$, there exists a finite set $C\subseteq E$, such that $\{fh(e): e\in C\}$ is a cover
for $f$.

\Proof \proofImply {i}{ii} Suppose by contradiction that (i) holds but (ii) fails for a given $f$ in $F$. Given any
finite subset $C\subseteq E$, we then put
  $$
  U_C = \big\{\psi\in\tsp F: \psi(f) = 1, \text{ and } \psi(h(e)) = 0, \ \forall e\in C\big\},
  $$
  and we claim that $U_C$ is nonempty. To see this, observe that our assumptions imply that $\{fh(e): e\in C\}$ is not a
cover for $f$, so there exists some nonzero $g\leq f$, such that $g\perp fh(e)$, for all $e$ in $C$. Choosing any tight
character $\psi$ on $F$ such that $\psi(g) = 1$, one can easily prove that $\psi$ lies in $U_C$, proving the claim.

If $C_1$ and $C_2$ are finite subsets of $E$, observe that
  $$
  U_{C_1}\cap U_{C_2} = U_{C_1\cup C_2},
  $$
  which in turn implies that the $U_C$ satisfy the finite intersection property. As each $U_C$ is compact, we conclude
that
  $$
  \bigcap U_C\neq\varnothing,
  $$
  where the intersection ranges over all finite subsets $C\subseteq E$. Choosing some $\psi$ in the above intersection,
we have that $\psi(f) = 1$, and $\psi(h(e)) = 0$, for all $e$ in $E$. So, in particular $\psi\circ h = 0$, contradicting
condition (i).

\Medskip\proofImply{ii}{i} Let $\psi$ be a tight character on $F$, and pick $f$ in $F$ such that $\psi(f) = 1$. Assuming
(ii), choose a finite set $C\subseteq E$, such that $\{fh(e): e\in C\}$ is a cover for $f$. By tightness we then have
that
  $$
  1 = \psi(f) = \bigvee_{e\in C} \psi(fh(e)),
  $$
  so there is some $e\in C$ such that
  $$
  1 = \psi(fh(e)) = \psi(f)\psi(h(e)) = \psi(h(e)),
  $$
  whence $\psi\circ h$ is not identically zero. \endProof

We should point out that condition (3.5.a) looks very much like condition (4.1.ii), with the significant difference that
the latter does not require that $h(e)\menor f$. For example, if $E$ is the semilattice $\{0,e,f, 1\}$, with $0\leq e
\leq1$, $0\leq f \leq1$, and $e\perp f$, then the inclusion of the sub-semilattice $\{0,1\}$ into $E$ satisfies
(4.1.ii), but it is not tightly surjective.

The previous result pinpoints the precise conditions under which the map
  $$
  \hath: \psi\in\tsp F \mapsto\psi\circ h \in\hatE\eqno{(4.2)}
  $$
  is well defined, but, due to our interest in tight spectra, it would be desirable for $\hath$ to take values in $\tsp
E$.

\state 4.3. Proposition. Assume that $h$ satisfies the equivalent conditions of (4.1). Then the following are
equivalent:
  \item {(i)} for every tight character $\psi$ on $F$, one has that $\psi\circ h$ is a tight character on E,
  \item {(ii)} if $\{e_1,\ldots,e_n\}$ is a cover for a given $e\in E$, then $\{h(e_1),\ldots,h(e_n)\}$ is a cover for
$h(e)$.

\Proof \proofImply {i}{ii} Assuming (i), let $\{e_1,\ldots,e_n\}$ be a cover for $e\in E$. If, by contradiction,
$\{h(e_1),\ldots,h(e_n)\}$ is not a cover for $h(e)$, then there exists a nonzero $g\leq h(e)$, such that $g\perp
h(e_i)$, for all $i$. Choosing a tight character $\psi$ on $F$, such that $\psi(g) = 1$, we have that $\psi\circ h$ is
tight by hypothesis, so the fact that $\{e_1,\ldots,e_n\}$ covers $e$ yields
  $$
  \bigvee_{i = 1}^n \psi(h(e_i)) = \psi(h(e)) \geq\psi(g) = 1.
  $$
  However, for every $i$ we have that
  $$
  \psi(h(e_i)) = \psi(g) \psi(h(e_i)) = \psi(gh(e_i)) = 0,
  $$
  because $g\perp h(e_i)$. This is a contradiction, so (ii) is proved.

\Medskip\proofImply{ii}{i} Given a tight character $\psi$ on $F$, one has that $\psi\circ h$ is nonzero, hence a
character, by hypothesis. The fact that $\psi\circ h$ is tight is then proved by a routine verification of the
definition. \endProof

Putting together (4.1) and (4.3), we obtain the following immediate consequence:

\state 4.4. Corollary. Given semilattices $E$ and $F$, as well as a homomorphism $h:E\to F$, the following are
equivalent:
  \item {(i)} for every tight character $\psi$ on $F$, one has that $\psi\circ h$ is a (nonzero) tight character on E,
  \item {(ii)} $h$ satisfies (4.1.ii) and (4.3.ii).

It is therefore natural to restrict one's attention to semilattice homomorphisms satisfying the above conditions. To
make this point a bit more explicit we give the following:

\state 4.5. Definition. A semilattice homomorphism $h:E\to F$ will be called \emph{tight} if it satisfies the equivalent
conditions of (4.4).

Notice that any character $\varphi:E\to\{0,1\}$ is a homomorphism, so one may ask whether or not it satisfies the above
definition, and it is easy to see that the answer is affirmative if and only $\varphi$ is a tight character in the usual
sense of the word. Put in another way, the above definition causes no conflict with the usual definition of tightness
for characters.

If $h$ is a tight homomorphism then (4.2) clearly provides a map
  $$
  \hath: \psi\in\tsp F \mapsto\psi\circ h \in\tsp E, \eqno{(4.6)}
  $$
  which will now be the main focus of this section.

Recall from section (2) that there is a natural bijective correspondence between the set of characters and the set of
filters. If we therefore view tight spectra as sets of tight filters, the appropriate interpretation of the map $\hath$,
above, is as follows. For any tight filter $\eta$ on $F$, one has
  $$
  \hath(\eta) = \{e\in E: h(e)\in\eta\}. \eqno{(4.7)}
  $$

\state 4.8. Proposition. Let $h:E\to F$ be a tight semilattice homomorphism. Then
  \item {(i)} $\hath$ is continuous,
  \item {(ii)} if $\hath$ is bijective, then $\hath^{-1}$ is continuous and hence $\hath$ is a homeomorphism.

\Proof We leave it for the reader to prove (i). Assuming that $\hath$ is bijective, pick $\varphi_0$ in $\tsp E$, and
choose any $e$ in $E$, such that $\varphi_0(e) = 1$. The obvious fact that, for every $\psi$ in $\tsp F$, one has that
  $$
  \psi(h(e)) = 1 \iff\hath(\psi)(e) = 1,
  $$
  implies that $\hath$ maps the set
  $$
  U = \big\{\psi\in\tsp F: \psi(h(e)) = 1\big\}
  $$
  onto the set
  $$
  V = \big\{\varphi\in\tsp E: \varphi(e) = 1\big\}.
  $$
  Since $U$ and $V$ are compact, and since $\hath$ is continuous, it follows that $\hath^{-1}$ is continuous as a map
from $V$ to $U$. As $V$ is also an open neighborhood of $\varphi_0$, we deduce that $\hath^{-1}$, viewed as a map from
$\tsp E$ to $\tsp F$, is continuous at $\varphi_0$. \endProof

With the next result we start investigating conditions on $h$ which will later characterize the fact that $\hath$ is a
homeomorphism.

\state 4.9. Theorem. If $h:E\to F$ is a tight semilattice homomorphism, then $\hath$ is surjective if and only if $h$ is
tightly injective.

\Proof Suppose that $h$ is tightly injective. Given $\varphi\in\tsp E$, we need to prove that $\varphi$ lies in the
range of $\hath$.

Suppose first that $\varphi$ is an ultra-character, and let $\xi$ be the support of $\varphi$, so that $\xi$ is an
ultra-filter. It is easy to see that
  $$
  \eta_0 := \{h(e): e\in\xi\}
  $$
  is a filter-base, the fact that $\eta_0$ does not contain zero being a consequence of $\ker(h)$ being trivial. Letting
$\eta$ be the filter generated by $\eta_0$, we may use Zorn's lemma to produce an ultra-filter $\eta'$ on $F$ containing
$\eta$. We next consider the filter
  $$
  \xi' = \hath(\eta') \explica{(4.7)}{=} \{e\in E: h(e)\in\eta'\},
  $$
  It is elementary to check that $\xi\subseteq\xi'$, so $\xi= \xi'$ by maximality. Letting $\psi$ be the characteristic
function of $\eta'$, we have that $\psi$ is an ultra-character, hence tight. In addition, for every $e$ in $E$ one has
that
  $$
  \psi(h(e)) = [h(e)\in\eta'] = [e\in\xi'] = [e\in\xi] = \varphi(e),
  $$
  proving that $\varphi= \psi\circ h = \hath(\psi)$, and hence that $\varphi$ lies in the range of $\hath$.

Now let $\varphi$ be any tight character. Recalling that the ultra-characters are dense in $\tsp E$, we may find a net
$\{\varphi_i\}_i$ of ultra-characters such that $\varphi= \lim_i\varphi_i$. Using the first part of this proof we may
pick ultra-characters $\psi_i$ on $F$ such that $\varphi_i = \psi_i\circ h$, for every $i$. Fixing any $e_0$ in $E$ such
that $\varphi(e_0) = 1$, observe that
  $$
  1 = \varphi(e_0) = \lim_i \varphi_i(e_0) = \lim_i \psi_i(h(e_0)),
  $$
  so, upon passing to a subnet, we may assume that $\psi_i(h(e_0)) = 1$, for all $i$. The subset of $\tsp F$ given by
  $$
  \big\{\psi\in\tsp F : \psi(h(e_0)) = 1\big\},
  $$
  (see also 2.9), is therefore a compact set to which all of the $\psi_i$ belong. So there is a converging subnet
$\{\psi_{i_j}\}_j$, whose limit we denote by $\psi$. It then follows that, for all $e$ in $E$, we have that
  $$
  \psi(h(e)) = \lim_j \psi_{i_j}(h(e)) = \lim_j \varphi_{i_j}(e) = \varphi(e),
  $$
  proving that $\varphi= \psi\circ h$, and hence that $\hath$ is surjective.

In order to prove the other implication, assume that $\hath$ is surjective. Given any nonzero $e$ in $E$, choose any
tight character $\varphi$ on $E$, such that $\varphi(e) = 1$. Taking $\psi$ in $\tsp F$ with $\hath(\psi) = \varphi$, we
have that
  $$
  1 = \varphi(e) = \hath(\psi)(e) = \psi(h(e)),
  $$
  so $h(e) \neq0$, proving that $\ker(h) = 0$, and hence that $h$ is tightly injective. \endProof

Of course we also want to characterize when is $\hath$ injective, but this requires a bit more work.

\section 5 Spectral order

\firstLine Given a semilattice $E$, we have already mentioned that every ultra-filter on $E$ is tight, although the
converse doesn't always hold. If $\xi$ is a non-maximal tight filter, then Zorn's lemma yields an ultra-filter $\zeta$
such that $\xi\varsubsetneq\zeta$. This said we see that $\tsp E$ may be equipped with a nontrivial partial order
relation, namely the order of inclusion.

On the other hand, if $\varphi$ and $\psi$ are tight characters, and $\xi$ and $\eta$ are their respective supports, its
is clear that
  $$
  \xi\subseteq\eta\iff\varphi\leq\eta,
  $$
  where the inequality in the right-hand-side refers to pointwise order. If one prefers to see $\tsp E$ as the set of
all tight characters, as opposed to tight filters, then the relevant order is pointwise order.

To the best of our knowledge, this order relation has not played any relevant role in the literature so far, but, since
it will be crucial for our goals, it is worth emphasizing it in the following:

\state 5.1. Definition. If $\varphi$ and $\psi$ are tight characters on the semilattice $E$, we shall say that
$\varphi\leq\psi$ if $\varphi(e)\leq\psi(e)$, for all $e$ in $E$. This order relation will be called the \emph{spectral
order}, and henceforth we shall always view $\tsp E$ as a partially ordered set equipped with this order relation.

Order relations obviously play a fundamental role in the study of semilattices. However the order relation just
introduced will play a very specific role and, in particular, it does not make $\tsp E$ into a semilattice, as
elementary examples show.

There are many important examples of semilattices in which every tight filter is maximal (see e.g. \cite[Theorem
3.4]{reconstru}). If $E$ is such a semilattice, the above order relation on $\tsp E$ is clearly not very enlightening as
$\varphi\leq\psi$ if and only if $\varphi= \psi$. This might perhaps explain why this order relation has not received
much attention so far. Nevertheless what follows will bring it to the spotlight.

We now make a short digression to settle some terminology regarding order relations to be used in the sequel.

\state 5.2. Definition. If $A$ and $B$ are partially ordered sets, and if $f:A\to B$ is a map, we shall say that:
  \item {(a)} a subset $U\subseteq A$ is an \emph{up-set} if, for every $x$ and $y$ in $A$, one has that $y\geq x\in U$
implies that $y\in U$,
  \item {(b)} $f$ is \emph{order-preserving} if, for all $a_1,a_2\in A$, one has that
  $$
  a_1\leq a_2 \imply f(a_1)\leq f(a_2),
  $$
  \item {(c)} $f$ is \emph{order-injective} if, for all $a_1,a_2\in A$, one has that
  $$
  f(a_1)\leq f(a_2) \imply a_1\leq a_2,
  $$
  \item {(d)} $f$ is an \emph{order-isomorphism} if it is bijective, order-preserving and order-injective. \Medskip If,
besides being a partially ordered set, $A$ is also a topological space, we will say that $A$ is an \emph{ordered
topological space} (we shall require no compatibility between the order and the topological structures). Assuming that
$A$ and $B$ are ordered topological spaces, we shall moreover say that: \smallskip
  \item {(e)} $f$ is an \emph{order-homeomorphism} if $f$ is both an order-isomorphism and a homeomorphism,
  \item {(f)} $A$ and $B$ are \emph{order-homeomorphic} if there exists an order-homeomorphism between $A$ and $B$.

\medskip

Filters are clearly examples of up-sets, but we will shortly deal with up-sets that are not subsets of any semilattice.

Using the fact that partial orders are reflexive and anti-symmetric, it is easy to see that every order-injective map is
injective. It should however be noted that (3.6.iii) does not qualify as an example of order-injectivity due to the fact
that ``$\menor$'' is not an anti-symmetric relation, hence not a partial order. In particular (3.6.iii) does not imply
that $h$ is injective.

\state 5.3. Proposition. If $h:E\to F$ is a tight semilattice homomorphism, then $\hath$ is order-preserving.

\Proof Given $\psi_1, \psi_2\in\tsp F$, such that $\psi_1\leq\psi_2$, then, for all $e$ in $E$, we have
  $$
  \hath(\psi_1)(e) = \psi_1(h(e))\leq\psi_2(h(e)) = \hath(\psi_1)(e),
  $$
  so $\hath(\psi_1)\leq\hath(\psi_1)$. \endProof

What follows is our first important use of the spectral order.

\state 5.4. Theorem. Let $h:E\to F$ be a tight semilattice homomorphism. Then $\hath$ is order-injective if and only if
$h$ is tightly surjective.

\Proof Assume that $\hath$ is order-injective. Given $f$ in $F$, let
  $$
  S = \{e\in E: h(e)\menor f\}.
  $$
  Arguing by contradiction, assume that, for every finite subset $C\subseteq S$, one has that $\{fh(e): e\in C\}$ does
not cover $f$. Therefore, for each such $C$, there exist a nonzero $g\leq f$, such that $g\perp fh(e)$, for all $e\in
C$. Choosing any tight character $\psi$ such that $\psi(g) = 1$, we have that
  $$
  \psi(f) = 1, \and\psi(h(e)) = 0, \for e\in C.
  $$
  In particular this says that the set
  $$
  U_C = \big\{\psi\in\tsp F: \psi(f) = 1,\text{ and } \psi(h(e)) = 0,\ \forall e\in C\big\}
  $$
  is nonempty. If $C_1$ and $C_2$ are finite subsets of $S$, it is easy to see that
  $$
  U_{C_1}\cap U_{C_1} = U_{C_1\cup C_2},
  $$
  so we see that the collection of all sets of the form $U_C$, as above, satisfies the finite intersection
property. Since they are all compact, their intersection is nonempty, so we may choose some
  $$
  \psi\in\bigcap U_C,
  $$
  where the intersection ranges over all finite subsets $C\subseteq S$. We then have that
  $$
  \psi(f) = 1, \and\psi(h(e)) = 0, \for e\in S. \eqno{(5.4.1)}
  $$
  Letting $\eta$ be the support of $\hath(\psi)$, we then have that,
  $$
  e\in\eta\imply\hath(\psi)(e) = 1 \imply\psi(h(e)) = 1 \explica{(5.4.1)}{\imply} e\notin S \imply h(e)\not\menor f.
  $$

For each $e$ in $\eta$ we then have that there exists a nonzero $g$ in $F$, such that $f\perp g \leq h(e)$. Choosing any
tight character $\chi$ on $F$ satisfying $\chi(g) = 1$, we have that
  $$
  \chi(f) = 0, \and\chi(h(e)) = 1.
  $$
  In particular this says that, for every $e\in\eta$, the set
  $$
  V_e = \big\{\chi\in\tsp F: \chi(f) = 0,\text{ and } \chi(h(e)) = 1\big\}
  $$
  is nonempty. Given $e_1$ and $e_2$ in $\eta$, it is easy to see that
  $$
  V_{e_1}\cap V_{e_2} = V_{e_1e_2},
  $$
  from where we deduce that the collection of all sets of the form $V_e$, with $e\in\eta$, satisfies the finite
intersection property. Since they are all compact, their intersection is nonempty, so we may choose some
  $$
  \psi'\in\bigcap_{e\in\eta} V_e, \eqno{(5.4.2)}
  $$
  and we let $\eta'$ be the support of $\hath(\psi')$. For every $e$ in $\eta$ we then have that
  $$
  1\explica{(5.4.2)}{=} \psi'(h(e)) = \hath(\psi')(e) = [e\in\eta'],
  $$
  so that $e\in\eta'$. The conclusion is then that $\eta\subseteq\eta'$, whence $\hath(\psi) \leq\hath(\psi')$. The fact
that $\hath$ is order-injective then gives $\psi\leq\psi'$, whence
  $$
  1 \explica{(5.4.1)}{=} \psi(f) \leq\psi'(f) \explica{(5.4.2)}{=}0,
  $$
  a contradiction.

In order to prove the converse, let $\varphi, \psi\in\tsp F$, and assume that $\hath(\varphi)\leq\hath(\psi)$. Given any
$f$ in $F$, we are then required to prove that $\varphi(f)\leq\psi(f)$, which incidentally is equivalent to proving that
  $$
  \varphi(f) = 1\imply\psi(f) = 1,
  $$
  so we duly assume that $\varphi(f) = 1$. By hypothesis we have that there exists a finite set $C\subseteq E$, such
that $h(e)\menor f$, for all $e\in C$, and $\{fh(e): e\in C\}$ is a cover for $f$. As $\varphi$ is tight, we conclude
that
  $$
  1 = \varphi(f) = \bigvee_{e\in C} \varphi(fh(e)) = \bigvee_{e\in C} \varphi(h(e)),
  $$
  so there exists some $e$ in $C$ satisfying $\varphi(h(e)) = 1$. Consequently
  $$
  1 = \varphi(h(e)) = \hath(\varphi)(e) \leq\hath(\psi)(e) = \psi(h(e)) \explica{(3.2)}{\leq} \psi(f),
  $$
  whence $\psi(f) = 1$, concluding the proof. \endProof

Yes, at the end of Section (4) we did express a desire to characterize when is $\hath$ injective, while the above result
clearly falls a bit short, only characterizing order-injectivity. However, the fundamental importance of the spectral
order will instead reorient our interest to order-injectivity. The question of when is $\hath$ injective can
nevertheless be left as an open problem for the interested reader.

Based on our previous results it is therefore interesting to focus on the class of homomorphisms $h$ for which one may
apply: \bitem Corollary (4.4), requiring \underbar{(4.1.ii)} and \underbar{(4.3.ii)}, and guaranteeing that $\hath$ is a
well defined map from $\tsp F$ to $\tsp E$, \bitem Theorem (4.9), requiring that $h$ be \underbar{tightly injective},
and guaranteeing that $\hath$ is surjective, \bitem Theorem (5.4), requiring that $h$ be \underbar{tightly surjective},
and guaranteeing that $\hath$ is order-injective.

\Medskip If these are taken independently, it might look like one would need all four underlined conditions to
characterize the above class of homomorphisms. However, if taken together, we get a 50{\percent} discount.

\state 5.5. Proposition. Let $h:E\to F$ be a semilattice homomorphism which is both tightly injective and tightly
surjective. Then $h$ satisfies (4.1.ii) and (4.3.ii).

\Proof Tight surjectivity is obviously stronger than (4.1.ii), as already mentioned. Regarding (4.3.ii), let
$\{e_1,\ldots,e_n\}$ be a cover for a given $e\in E$, and suppose by contradiction that $\{h(e_1),\ldots,h(e_n)\}$ is
not a cover for $h(e)$. Then there exists a nonzero $f\leq h(e)$, such that $f\perp h(e_i)$, for all $i$. Since $h$ is
tightly surjective, there exists a finite set $C\subseteq E$, such that $h(c)\menor f$, for all $c\in C$, and
  $$
  \{fh(c): c\in C\} \eqno{(5.5.1)}
  $$
  is a cover for $f$. For $c\in C$, observe that
  $$
  h(c)\menor f \leq h(e),
  $$
  so $h(c)\menor h(e)$, and hence $c\menor e$, due to the fact that $h$ is tightly injective.

We then claim that $ce = 0$, for all $c\in C$. Embedding a proof by contradiction within another, suppose that
$ce\neq0$, for some $c$ in $C$. Then, as $ce\leq e$, there would be some $i$, such that $ce\Cap e_i$, whence
$h(cee_i)\neq0$, again because $h$ is tightly injective. This would in turn give
  $$
  0\neq h(cee_i) \leq h(c)h(e_i) \explica{(3.14.ii)}{\menor} fh(e_i) = 0,
  $$
  a contradiction, proving that $ce = 0$. The combined facts that $e\perp c \leq c$, and that $c\menor e$, then imply
that $c = 0$. This says that the cover for $f$ mentioned in (5.5.1) is the singleton $\{0\}$, which can only mean that
$f = 0$, a contradiction. \endProof

Having performed the above cleanup, we may now introduce one of the central concepts of this work.

\state 5.6. Definition. Let $E$ and $F$ be semilattices and let $h:E\to F$ be a homomorphism. Then $h$ will be called a
\emph{consonance} if $h$ is both tightly injective and tightly surjective.

As a consequence we have the following key result.

\state 5.7. Corollary. Let $h:E\to F$ be a semilattice homomorphism. Then the following are equivalent:
  \item {(i)} $h$ is tight and $\hath$ is an order-isomorphism,
  \item {(ii)} $h$ is a consonance. \medskip\noindent In this case $\hath$ is also an order-homeomorphism.

\Proof The equivalence (i$\Leftrightarrow$ii) clearly follows from (4.4), (4.9), (5.3), (5.4) and (5.5), while the last
sentence in the statement follows from (4.8). \endProof

If $h$ is a consonance, we therefore see that $\hath$ is invertible, so it would be desirable to obtain an explicit
formula for the inverse of $\hath$, and our next result does exactly that.

\state 5.8. Proposition. Let $h:E\to F$ be a consonance. Then, for every tight filter $\xi$ on $E$, the set
  $$
  \ch(\xi) = \big\{f\in F: h(e)\menor f, \text{ for some } e\in\xi\big\}
  $$
  is a tight filter on $F$, and $ \ch(\xi) = \hath^{-1}(\xi). $ Consequently $\ch$ is an order-homeomorphism from $\tsp
E$ to $\tsp F$.

\Proof Given $\xi$ as above, let us shorten our notation by writing $\eta= \hath^{-1}(\xi)$. We must then prove that $
\ch(\xi) = \eta. $ Given $f$ in $\eta$, and using that $h$ is tightly surjective, pick a finite set $C\subseteq E$, such
that $h(e)\menor f$, for all $e\in C$, and $\{fh(e): e\in C\}$ is a cover for $f$. Since $f$ lies in $\eta$, the
tightness condition implies that $fh(e)\in\eta$, for some $e$ in $C$. This clearly implies that $h(e)$ lies in $\eta$,
so (4.7) implies that $e$ lies in $\hath(\eta) = \xi$, hence proving that $f$ belongs to $\ch(\xi)$, and showing that
$\eta\subseteq\ch(\xi)$.

In order to prove the reverse inclusion, pick $f$ in $\ch(\xi)$, and take $e$ in $\xi$ such that $h(e)\menor f$. So,
  $$
  e\in\xi= \hath(\eta) \relimply h(e)\in\eta\explica{(3.2)}{\relimply} f\in\eta,
  $$
  as required. \endProof

\section 6 Inverse semigroup homomorphisms

\firstLine The plan for the rest of this work is to extend our study of semilattice homomorphisms started in section (4)
to the case of inverse semigroups. Thus, in order to be able to quickly refer to the general setup we will adopt from
now on, we make the following:

\state 6.1. Hypothesis. \rm From now on we will fix inverse semigroups $S$ and $T$ (as always, assumed to have zero),
and a homomorphism
  $$
  h:S\to T,
  $$
  (as always, assumed to send zero to zero). We will moreover denote the idempotent semilattices of $S$ and $T$ by $E$
and $F$, respectively, and we will let
  $$
  \hz:E\to F
  $$
  be the map obtained by restricting $h$. The standard action of $S$ on $\tsp E$ will be denoted by $\alpha$, and the
standard action of\/ $T$ on $\tsp F$, by $\beta$ (see 2.17). Accordingly we shall denote the corresponding fundamental
slices introduced in (2.14) by
  $$
  \Delta^\alpha_s = \big\{[s,\xi]_\alpha:\xi\in\dom(\alpha_s)\big\}, \and\Delta^\beta_t =
\big\{[t,\eta]_\beta:\eta\in\dom(\beta_t)\big\}.
  $$

In the present section we want to focus on aspects of inverse semigroups that transcend their idempotent
semilattices. So, in order not to worry about aspects that are specific to semilattices, we will often put ourselves in
the very comfortable situation where the restricted map $ \hz$ is assumed to satisfy all of the nice properties so far
mentioned regarding semilattice homomorphisms, namely that $\hz$ is a consonance.

Incidentally, we would also like to point out that, should one decide to formalize the passage from a semilattice to its
tight spectrum through a functor, our emphasis on the correspondence $h\to\hath$, above, definitely suggests a
contravariant functor. However, if we want to extend this to the correspondence sending an inverse semigroup to its
tight groupoid, which is roughly what we want to do next, it would not be a good idea to insist on contravariance. The
reason being that, if we take a group $G$, seen as an inverse semigroup (after the addition of a zero), then its tight
groupoid is $G$, itself, and this doesn't seem to fit well in a contravariant functor. This dilemma in fact shows up in
other places, most notably in the construction of C*-algebras from groupoids.

In order to sidestep this issue altogether we will henceforth focus on $\ch$, rather than on $\hath$, although we will
have to pay the price of dealing only with homomorphisms $h$ for which $\hz$ is a consonance.

Assuming that $\hz$ is a consonance, (5.7) implies that the dual map $\hat\hz$ is a homeomorphism from $\tsp F$ to $\tsp
F$, and hence so is the map
  $$
  \chz:\tsp E\to\tsp F
  $$
  introduced in (5.8).

\state 6.2. Proposition. Under (6.1), and assuming that $\hz$ is a consonance, one has that $(h, \chz)$ is a covariant
epimorphism in the sense of (2.24).

\Proof Fixing $s$ in $S$, and given $\xi\in\tsp E$, we claim that
  $$
  s^*s\in\xi\reliff\hz(s^*s)\in\chz(\xi).
  $$
  Letting $\eta= \chz(\xi)$, we have that $\hat\hz(\eta) = \xi$, and the above is equivalent to
  $$
  s^*s\in\hat\hz(\eta) \reliff\hz(s^*s)\in\eta,
  $$
  which is evident, and hence the claim is proved, which in turn translates into saying that
  $$
  \xi\in\dom(\alpha_s) \reliff\chz(\xi)\in\dom(\beta_{h(s)}).
  $$
  So, as $\chz$ is bijective, we deduce that
  $$
  \chz(\dom(\alpha_s)) = \dom(\beta_{h(s)}), \eqno{(6.2.1)}
  $$
  and (2.21.1) follows (we shouldn't yet brag that $(h, \chz)$ is a covariant epimorphism because we do not even know
that it is covariant). We will now check that
  $$
  \chz\big(\alpha_s(\xi)\big) = \beta_{h(s)}\big(\chz(\xi)\big),
  $$
  for all $\xi$ in $\dom(\alpha_s)$. Recalling that $\chz$ is the inverse of $\hathz$, and by substituting
$\hathz(\eta)$ for $\xi$, with $\eta$ in $\dom(\beta_{h(s)})$, we see that this is equivalent to saying that
  $$
  \alpha_s\big(\hathz(\eta)\big) = \hathz\big(\beta_{h(s)}(\eta)\big).
  $$
  The character point of view makes this task easier, so we view $\eta$ as a character on $F$, according to (2.1). For
every $e$ in $E$, we then have that
  $$
  \hathz\big(\beta_{h(s)}(\eta) \big)(e) = \beta_{h(s)}(\eta)\big(h(e)\big) = \eta\big(h(s)^*h(e)h(s)\big) =
  $$
  $$
 = \eta\big(h(s^*es)\big) = \hathz(\eta)(s^*es) = \alpha_s \big(\hathz(\eta) \big)(e),
  $$
  and we are done proving that $(h,\chz)$ is a covariant pair. That this is actu\-ally a covariant epimorphism then
follows immediately from (6.2.1). \endProof

From (2.22) we then obtain a continuous homomorphism between tight groupoids
  $$
  h{\times}\chz: \Gt{S}\to\Gt{T},
  $$
  such that
  $$
  (h{\times}\chz)\big([s,\xi]_\alpha\big) = [h(s),\chz(\xi)]_\beta,
  $$
  for every germ $[s,\xi]_\alpha$ in $\Gt{S}$.

\state 6.3. Notation. Under (6.1), and assuming that $\hz$ is a consonance, the map $h{\times}\chz$ defined above will
be denoted by $\ch$.

Of course this notation was already used in (5.8), but, if we treat semilattices as inverse semigroups in (6.3), it is
easy to see that the two notations refer to the same mathematical object. In other words, the above abuse of language
causes no conflicts.

In any case, the restriction of $\ch$ to the corresponding unit spaces, namely $\chz$, is a homeomorphism by our
assumptions, and hence evidently also an injective open map. So (2.23) applies and we get an inverse semigroup
homomorphism
  $$
  \hzao:U\subseteq\Gt{S}\op\mapsto\ch(U)\in\Gt{T}\op. \eqno{(6.4)}
  $$

Since $(h, \chz)$ is also a covariant epimorphism, we may use the last conclusion of (2.25) to produce the following
commutative diagram:

\begingroup\noindent\hfill\beginpicture\setcoordinatesystem units <0.025truecm, -0.02truecm> \setplotarea x from -40 to
150, y from -40 to 180 \put{\null} at -40 -40 \put{\null} at -40 180 \put{\null} at 150 -40 \put{\null} at 150 180
\put{$S$} at 0 0 \put{$T$} at 110 0 \arrow<0.11cm> [0.5, 1.8] from 16.5 0 to 93.5 0 \put{$h$} at 55 -18 \put{$\Gt S\op$}
at 0 100 \put{$\Gt T\op$} at 110 100 \arrow<0.11cm> [0.5, 1.8] from 33 100 to 77 100 \put{$\hzao$} at 55 82
\arrow<0.11cm> [0.5, 1.8] from 0 25 to 0 75 \put{$\Delta^\alpha$} at -15 50 \arrow<0.11cm> [0.5, 1.8] from 110 25 to 110
75 \put{$\Delta^\beta$} at 125 50 \put{Diagram 6.5} at 55 140 \endpicture\hfill\null\endgroup

We would now like to generalize (4.9) to the present setting, but before that we'd like to point out that the play of
words between injectivity/surjectivity appearing in its statement is partly due to the fact that the correspondence
$h\to\hath$ is contravariant. Since the statement of our generalization below will emphasize $\ch$, rather than $\hath$,
this play of words will not appear.

\state 6.6. Proposition. Under (6.1), and assuming that $\hz$ is a consonance, one has that $\ch$ is injective if and
only if $h$ is tightly injective.

\Proof Assume that $h$ is tightly injective and that $[s,\xi]_\alpha$ and $[t,\eta]_\alpha$ are germs in $\Gt{S}$, such
that
  $$
  \ch([s,\xi]_\alpha) = \ch([t,\eta]_\alpha),
  $$
  or, equivalently, that
  $$
  [h(s),\chz(\xi)]_\beta= [h(t),\chz(\eta)]_\beta.
  $$
  Then $\chz(\xi) = \chz(\eta)$, so $\xi= \eta$, because $\chz$ is a bijective map. Moreover, there exists $f$ in $F$
(the idempotent semilattice of $T$), such that $\chz(\xi)$ lies in the domain of $\beta_f$, and
  $$
  h(s)f = h(t)f.
  $$
  To say that $\chz(\xi)$ is in $\dom(\beta_f)$ is to say that $f\in\chz(\xi)$, which in turn implies (see 5.8) that
there exists $e$ in $\xi$, such that $h(e)\menor f$. Employing (3.14.i) we get
  $$
  h(s)h(e)\equiva h(t)h(e),
  $$
  which is when the tight injectivity hypothesis intervenes to give $se\equiva te$. In particular, there exists a finite
cover $C$ of $(se)^*se$, such that
  $$
  sec = tec, \for c\in C. \eqno{(6.6.1)}
  $$
  Observe that $e$ lies in $\xi$, and so does $s^*s$, by reason of $[s,\xi]_\alpha$ being a legitimate germ. Therefore
  $$
  (se)^*se = s^*se \in\xi,
  $$
  so the tightness of $\xi$ implies that there exists some $c$ in $C$, lying in $\xi$. It follows that
$\xi\in\dom(\alpha_{ec})$, and (6.6.1) then implies that $[s,\xi]_\alpha= [t,\eta]_\alpha$, concluding the proof that
$\ch$ is injective.

Conversely, assume that $\ch$ is injective, and that $r,s\in S$ are such that $h(r)\equiva h(s)$. Then, referring to
Diagram (6.5), we have that
  $$
  \hzao(\Delta^\alpha_r) = \Delta^\beta_{h(r)} \explica{(3.13)}{=} \Delta^\beta_{h(s)} = \hzao(\Delta^\alpha_s).
  $$
  As $\ch$ is assumed to be injective, it is obvious that $\hzao$ is also injective, whence $\Delta^\alpha_r =
\Delta^\alpha_s$, and again by (3.13), we get $r\equiva s$. This concludes the proof. \endProof

We next generalize (5.4) to inverse semigroup homomorphisms.

\state 6.7. Proposition. Under (6.1), and assuming that $\hz$ is a consonance, one has that $\ch$ is surjective if and
only if $h$ is tightly surjective.

\Proof Let $[t, \eta]_\beta$ be any germ in $\Gt{T}$, so that $t^*t\in\eta$. Assuming that $h$ is tightly surjective,
choose a finite set $C\subseteq S$, such that $h(s)\menor t$, for every $s$ in $C$, and such that $ \{h(s^*s)t^*t:s\in
C\} $ is a cover for $t^*t$. Thanks to the tightness of $\eta$, we may pick some $s$ in $C$, such that $h(s^*s)t^*t$
lies in $\eta$, and hence so does $h(s^*s)$. This means that $\eta$ belongs to $\D_{h(s)^*h(s)}$, and then (3.12.i)
gives
  $$
  [h(s),\eta]_\beta= [t,\eta]_\beta.
  $$

Letting $\xi= \hathz(\eta)$, and temporarily interpreting $\xi$ and $\eta$ as characters, according to (2.1), we have
  $$
  \xi(s^*s) = \hathz(\eta)(s^*s) = \eta(h(s^*s)) = 1,
  $$
  which, back to filters, means that $s^*s\in\xi$. Therefore $[s,\xi]_\alpha$ is a well formed germ in $\Gt{S}$, and
  $$
  \ch([s,\xi]_\alpha) = [h(s),\chz(\xi)]_\beta= [h(s),\eta]_\beta= [t,\eta]_\beta.
  $$
  This proves that $\ch$ is surjective.

In order to prove the other implication, suppose that $\ch$ is surjective. Given $t\in T$, consider the subset $M$ of
$S$ defined by
  $$
  M = \big\{s\in S: h(s)\menor t\big\},
  $$
  and we claim that
  $$
  \Delta^\beta_t = \bigcup_{s\in M} \ch(\Delta^\alpha_s), \eqno{(6.7.1)}
  $$
  where $\Delta^\alpha_s$ and $\Delta^\beta_t$ are as in (6.1).

In order to prove that each $\ch(\Delta^\alpha_s)$ is contained in $\Delta^\beta_t$, pick $s$ in $M$, and observe that,
for each $\xi\in\dom(\alpha_s)$, we have that
  $$
  \ch([s,\xi]_\alpha) = [h(s),\chz(\xi)]_\beta= [t,\chz(\xi)]_\beta\in\Delta^\beta_t,
  $$
  where the last equality follows from (3.12.i).

In order to prove the inclusion ``$\subseteq$'' in (6.7.1), pick $\eta$ in $\dom(\beta_t)$, so that $t^*t\in\eta$. As we
are assuming that $\ch$ is surjective, we may choose $[s,\xi]_\alpha$ in $\Gt{S}$, such that
  $$
  \ch\big([s,\xi]_\alpha\big) = [h(s),\chz(\xi)]_\alpha= [t,\eta]_\beta.
  $$
  If follows that $\chz(\xi) = \eta$, and there exists some idempotent $f\in\eta$, such that $h(s)f = tf$.

Given that $f\in\chz(\xi)$, we may pick $e$ in $\xi$, with $h(e)\menor f$. Setting $s' = se$, we then have that
  $$
  h(s') = h(s)h(e) \menor h(s)f = tf \leq t,
  $$
  so we see that $s'$ lies in $M$. Noticing that $[s,\xi]_\alpha$ is clearly equal to $[s',\xi]_\alpha$, we evidently
have that
  $$
  \ch\big([s',\xi]_\alpha\big) = [t,\eta]_\beta.
  $$
  This concludes the proof of claim (6.7.1).

Observe that the restriction of $\ch$ to the unit spaces of the corresponding groupoids is $\chz$, which is a
homeomorphism, and hence also an open map. So (2.23) implies that $\ch$ is an open map, hence each
$\ch(\Delta^\alpha_s)$ appearing in (6.7.1) is open (besides being compact). Since $\Delta^\beta_t$ is compact (besides
being open), we may find a finite subset $C\subseteq M$, such that
  $$
  \Delta^\beta_t = \bigcup_{s\in C} \ch(\Delta^\alpha_s).
  $$
  Since $(h, \chz)$ is a covariant epimorphism, the last conclusion of (2.25) gives $\ch(\Delta^\alpha_s) =
\Delta^\beta_{h(s)}$, whence
  $$
  \Delta^\beta_t = \bigcup_{s\in C} \Delta^\beta_{h(s)}. \eqno{(6.7.2)}
  $$

In order to prove that $h$ is tightly surjective, meaning that its range is large, it then suffices to prove the last
part of (3.10.a), namely that
  $$
  \{h(s^*s)t^*t:s\in C\}
  $$
  is a cover for $t^*t$. So we suppose by contradiction that $f$ is a nonzero element of $F$, such that $f\leq t^*t$,
while
  $$
  f\perp h(s^*s)t^*t, \for s \in C.
  $$
  Choosing $\eta\in\tsp F$, such that $f\in\eta$, we have that necessarily $t^*t\in\eta$. So, denoting the source map of
the groupoid $\Gt{T}$ by $\d$, we have that
  $$
  \eta\in\dom(\beta_t) = \d(\Delta^\beta_t) \explica{(6.7.2)}{=} \bigcup_{s\in C} \d(\Delta^\beta_{h(s)}),
  $$
  and it follows that there exists some $s$ in $C$, such that $\eta\in\d(\Delta^\beta_{h(s)})$, meaning that
$h(s^*s)\in\eta$. This is however impossible, since it leads to
  $$
  0 = f h(s^*s)t^*t \in\eta. \closeProof
  $$
  \endProof

\section 7 Consonance for inverse semigroups

\firstLine Summarizing our work so far, we have developed a theory of consonances for semilattices in sections (3)-(5)
and, in section (6) we have studied homomorphisms of inverse semigroups whose restrictions to the corresponding
idempotent semilattices are consonances. Now we would like to put all of this together in order to arrive at a theory of
consonances for inverse semigroups. Our setup will henceforth be (6.1), but we will no longer insist in assuming that
$\hz$ is a consonance. Instead this will come as a consequence of the hypotheses made on $h$.

We are therefore ready to introduce the main concept of this work.

\state 7.1. Definition. An inverse semigroup homomorphism $h:S\to T$ will be called a \emph{consonance} if it is both
tightly injective and tightly surjective.

There is a little issue with respect to the behavior of the properties listed above when homomorphisms are restricted to
their idempotent semilattices, especially when the restriction process affects both domains and ranges. This is however
of easy solution.

\state 7.2. Proposition. Under (6.1), one has that:
  \item {(i)} if $h$ is tightly injective, then $\hz$ is tightly injective,
  \item {(ii)} if $h$ is tightly surjective, then $\hz$ is tightly surjective. \Medskip Therefore, if $h$ is a
consonance, then so is $\hz$.

\Proof Point (i) follows immediately from (3.8). On the other hand, proving (ii) amounts to verifying that $\hz(E)$ is a
large subset relative to $F$, given that $h(S)$ is a large subset relative to $T$.

Given $f$ in $F$, and seeing $f$ as an element of $T$, take a finite subset $C\subseteq h(S)$, such that $t\menor f$,
for all $t\in C$, and $\{ft^*t:t\in C\}$ is a cover for $f$. Observing that (3.12.ii) implies that $t^*t\menor f$, for
all $t\in C$, we see that the set
  $$
  C' = \{t^*t:t\in C\}
  $$
  is a subset of $\hz(E)$ satisfying condition (3.5.a) with respect to $f$. So $\hz$ is tightly surjective. \endProof

The relevance of the previous Proposition is that, when $h$ is a consonance, we obtain the privilege of using the
results of Section (6), where $\hz$ was always assumed to be a consonance. Therefore we immediately get the following
main result from (6.6), (6.7) and (2.23).

\state 7.3. Theorem. If $h:S\to T$ is a consonance, then the map
  $$
  \ch: \Gt{S}\to\Gt{T}
  $$
  introduced in (6.3) is well defined, and it is an isomorphism of topological groupoids. Conversely, if $\hz$ is a
consonance (so that $\ch$ is well defined), and if $\ch$ is bijective, then $h$ itself is a consonance.

The reader might have noticed that our main result on semilattice homomorphisms, namely Corollary (5.7), requires that
$\ch$ be an order-isomorphism in exchange for $h$ to be a consonance. However the last sentence of Theorem (7.3)
requires no order-related condition\fn{Of course, the hypothesis that $\hz$ be a consonance does involve a lot of
order-related information.} of $\ch$ to get that $h$ is a consonance. After all, the domain and range of $\ch$, namely
$\Gt{S}$ and $\Gt{T}$, are so far devoid of any order relation, so there is really nothing one can say regarding order
preservation for $\ch$.

In fact we shall see that tight groupoids are naturally equipped with a partial order relation and our next section will
be devoted to its investigation.

\section 8 Order on tight groupoids

\firstLine In this section we will equip tight groupoids with a very important order relation, which will complement
Theorem (7.3), and play a key role from now on.

\state 8.1. Definition. Let $S$ be an inverse semigroup, let $E$ be its idempotent semilattice, and let $\Gt{S}$ be its
tight groupoid. Given $\gamma,\delta\in\Gt{S}$, we will say that $\gamma\leq\delta$ if the following two conditions
hold:
  \item {(a)} $d(\gamma)\leq d(\delta)$, relative to the spectral order relation on $\tsp E$, and
  \item {(b)} there exists $s$ in $S$, such that both $\gamma$ and $\delta$ lie in $\Delta_s $. \Medskip This relation
will be called the \emph{spectral order} on $S$.

It then clearly follows that $\gamma\leq\delta$ if and only if there exists $s\in S$, as well as $\xi, \eta\in\tsp E$,
such that $s^*s\in\xi\subseteq\eta$, while
  $$
  \gamma= \gt s\xi, \and\delta= \gt s\eta, \eqno{(8.2)}
  $$

The reader's first impression regarding this relation might be the same as the author's, in the sense that nothing seems
to indicate that it is a transitive relation. Fortunately, after a while, one realizes that:

\state 8.3. Proposition. The relation ``$\leq$'' defined above is a partial order.

\Proof Reflexivity and anti-symmetry are trivial, so we focus on transitivity. Suppose that $\gamma$, $\delta$ and
$\varepsilon$ are elements in $\Gt{S}$, such that $\gamma\leq\delta$, and $\delta\leq\varepsilon$. Based on (8.2) we may
then write
  $$
  \gamma= \gt s\xi, \quad\delta= \gt s\eta= \gt t\eta, \and\varepsilon= \gt t\zeta,
  $$
  with
  $$
  s,t\in S, \quad\xi,\eta\in\D_{s^*s}, \quad\eta, \zeta\in\D_{t^*t}, \and\xi\subseteq\eta\subseteq\zeta.
  $$
  Observing that $s^*s\in\xi\subseteq\zeta$, it then follow that $\zeta\in\D_{s^*s}$, so the germ $\gt s\zeta$ is well
formed, and we claim that $\varepsilon= \gt s\zeta$. To see this, notice that, due to $\gt s\eta= \gt t\eta$, there is
an idempotent $e$ such that $\eta\in\D_e$, that is, $e\in\eta$, and $se = te$.

As $\eta\subseteq\zeta$, it follows that also $e\in\zeta$, so $\zeta\in\D_e$, and, recalling that $se = te$, we have
shown that $\gt s\zeta= \gt t\zeta$, as required.

Finally, since we may write $\gamma= \gt s\xi$, and $\varepsilon= \gt s\zeta$, and since $\xi\subseteq\zeta$, we have
proved that $\gamma\leq\varepsilon$. \endProof

The presence of the spectral order on tight groupoids will require us to consider this new category of objects, so it is
perhaps worth naming things properly.

\state 8.4. Definition. By a \emph{\otg/} we shall mean a topological groupoid equipped with a partial order relation
(no compatibility is required\fn{The purpose of Definition (8.4) is simply to introduce some new terminology and not
(yet) to initiate a general theory of ordered groupoids since we will only deal with the spectral order on tight
groupoids in the present section. For this reason we have decided to use the perhaps slightly awkward expression
``\otg/'' to reserve the more elegant ``ordered topological groupoid'' for later. The spectral order relation has many
interesting properties that will inspire an updated Definition to be given in section (12) with the purpose of
characterizing tight groupoids as those ordered groupoids satisfying the new set of axioms. We will also compare our
updated definition with Ehresmann's classical definition of ordered groupoids.} between the order relation and the
topological groupoid structure). If $G$ and $H$ are \otgs/, then a function
  $$
  \varphi:G\to H
  $$
  will be called an \emph{isomorphism} if it is, at the same time, an isomorphism of groupoids, a homeomorphism of
topological spaces, and an order-isomorphism. If such a function exists, then $G$ and $H$ will be said to be
\emph{isomorphic}.

Now that tight groupoids are equipped with their spectral order, we may complement Theorem (7.3), with the added bonus
that no extra hypothesis is needed for $\ch$ to respect the corresponding order relations.

\state 8.5. Theorem. If $h:S\to T$ is a consonance then $\ch$ is an isomorphism of \otgs/ from $\Gt{S}$ to $\Gt{T}$.

\Proof As we already know that $\ch$ is an isomorphism of topological groupoids, it suffices to prove it to be an
order-isomorphism. Given $\gamma, \delta\in\Gt{S}$. with $\gamma\leq\delta$, write $\gamma= \ga s\xi$, and $\delta= \ga
s\eta$, with $s^*s\in\xi\subseteq\eta$, as in (8.2). Then, by definition,
  $$
  \ch(\gamma) = \gb{h(s)}{\chz(\xi)}, \and\ch(\delta) = \gb{h(s)}{\chz(\eta)}.
  $$
  Since $\hz$ is also a consonance by (7.2), we have that $\chz$ is order-preserving whence
$\chz(\xi)\subseteq\chz(\eta)$, so it follows that $\ch(\gamma)\leq\ch(\delta)$. This proves that $\ch$ is
order-preserving.

Considering that $\ch$ is bijective by (7.3), in order to prove it to be order-injective, it is clearly enough to prove
that $\ch^{-1}$ is order-preserving. For this, let $\gamma'$ and $\delta'$ be elements of $\Gt{T}$, such that
$\gamma'\leq\delta'$, and write
  $$
  \gamma' = \gb t{\xi'}, \and\delta' = \gb t{\eta'},
  $$
  with $t^*t\in\xi'\subseteq\eta'$. Using that $h$ is tightly surjective, let $C\subseteq S$ be a finite set such that
$h(s)\menor t$, for all $s\in C$, and $\{h(s)^*h(s)t^*t:s\in C\}$ is a cover for $t^*t$. Since $\xi'$ is tight, we
deduce that $h(s)^*h(s)t^*t\in\xi'$, for some $s$ in $C$, and hence also that $h(s)^*h(s)\in\xi'$. As already noticed,
$\hz$ is a consonance, so $\chz$ is an order-homeomorphism by (5.8), and then we may write $\xi' = \chz(\xi)$, and
$\eta' = \chz(\eta)$, where $\xi$ and $\eta$ are elements of $\tsp E$, necessarily satisfying $\xi\subseteq\eta$. Given
that
  $$
  \hz(s^*s) = h(s)^*h(s) \in\xi' = \chz(\xi),
  $$
  we may find $e\in\xi$, such that $\hz(e) \menor\hz(s^*s)$, and hence also that $e \menor s^*s$, due to $\hz$ being
tightly injective. It then follows from (3.2) that
  $$
  s^*s\in\xi\subseteq\eta,
  $$
  so the germs
  $$
  \gamma:= \ga s\xi, \and\delta:= \ga s\eta
  $$
  are well formed and clearly $\gamma\leq\delta$. On the other hand,
  $$
  \ch(\gamma) = \ch(\ga s\xi) = \gb{h(s)}{\chz(\xi)} = \gb{h(s)}{\xi'} \explica{(3.12.i)}{=} \gb t{\xi'} = \gamma',
  $$
  and similarly $\ch(\delta) = \delta'$. So
  $$
  \ch^{-1}(\gamma') = \gamma\leq\delta= \ch^{-1}(\delta'),
  $$
  proving that $\ch^{-1}$ is order-preserving. \endProof

The next two results form the foundations for the proof of Theorem (9.5), one of the main results of this work.

\state 8.6. Proposition. If $S$ is an inverse semigroup, and $U\subseteq\Gt{S}$ is an open subset, then $U$ is an up-set
(Definition 5.2.a) relative to the spectral order if and only if there is family $\{s_i\}_{i\in I}$ of elements of $S$
such that
  $$
  U = \bigcup_{i\in I}\Delta_{s_i}.
  $$

\Proof We begin by proving the easier ``if'' part. Noticing that the union of up-sets is automatically an up-set, it
suffices to prove that $\Delta_s $ is an up-set for every $s$ in $S$. Fixing such an $s$, and supposing that $\gamma$
and $\delta$ are elements in $\Gt{S}$ such that
  $$
  \delta\geq\gamma\in\Delta_s,
  $$
  we need to prove that $\delta$ lies in $\Delta_s $. By (8.2) there is some $t$ in $S$, as well as $\xi$ and $\eta$ in
$\tsp E$, such that $t^*t\in\xi\subseteq\eta$, while
  $$
  \gamma= \gt t\xi, \and\delta= \gt t\eta,
  $$
  However, since $\gamma$ lies in $\Delta_s $, we may also represent it as
  $$
  \gamma= \gt s{\d(\gamma)} = \gt s\xi.
  $$
  From the equality $\gt t\xi= \gt s\xi$, it then follows that there exists some idempotent element $e$ such that
$\xi\in\D_e$, and $te = se$, so that
  $$
  e\in\xi\subseteq\eta.
  $$
  This said it is now easy to see that
  $$
  \delta= \gt t\eta= \gt s\eta\in\Delta_s,
  $$
  as required. Turning now to the ``only if'' part, and given $\gamma\in U$, it suffices to prove that there exists $s$
in $S$, such that
  $$
  \gamma\in\Delta_s \subseteq U.
  $$

Writing $\gamma= [s,\xi]_\theta$, with $\xi\in\D_{s^*s}$, we claim that there exists some $e\in\xi$, such that
$\Delta_{se}\subseteq U$. In order to verify this, assume by contradiction that $\Delta_{se}\sm U$ is nonempty for every
$e$ in $\xi$.

Recall that the source of the slice $\Delta_s $ is $\D_{s^*s}$, and hence the map
  $$
  d:\Delta_s \to\D_{s^*s}
  $$
  is a bijection. Therefore, for every $e_1$ and $e_2$ in $\xi$, we have that
  $$
  d\big(\Delta_{se_1}\cap\Delta_{se_2}\big) = d\big(\Delta_{se_1}\big) \cap d\big(\Delta_{se_2}\big) = \D_{s^*se_1}
\cap\D_{s^*se_2} = \D_{s^*se_1e_2} = d\big(\Delta_{se_1e_2}\big).
  $$
  Consequently $\Delta_{se_1}\cap\Delta_{se_2} = \Delta_{se_1e_2}$, so
  $$
  \big(\Delta_{se_1}\sm U\big)\cap\big(\Delta_{se_1}\sm U\big) = \Delta_{se_1e_2}\sm U,
  $$
  from where we see that the family of sets of the form $\Delta_{se}\sm U$, with $e$ in $\xi$, satisfies the finite
intersection property. Noticing that these are all closed subsets of the compact Hausdorff space $\Delta_s $, we may
choose an element $\delta$ belonging to the intersection of the whole family. Clearly $\delta$ lies in $\Delta_s $, so
we may write $\delta= [s,\eta]_\theta$, with $\eta\in\D_{s^*s}$. Since $\delta\in\Delta_{se}$, for every $e$ in $\xi$,
we deduce that
  $$
  \eta= d(\delta)\in\D_{s^*se},
  $$
  meaning that $s^*se\in\eta$, whence $e\in\eta$. The conclusion is then that $\xi\subseteq\eta$, so
  $$
  \gamma= \gt s\xi\leq\gt s\eta= \delta,
  $$
  Having assumed that $U$ is an up-set, we get that $\delta\in U$, a contradiction. Summarizing, we have proved our
claim that there exists some $e\in\xi$, such that $\Delta_{se}\subseteq U$, so
  $$
  \gamma= \gt s\xi= \gt{se}\xi\in\Delta_{se}\subseteq U,
  $$
  as desired. \endProof

If $A$ is any partially ordered set, observe that the intersection of any two up-sets is obviously an up-set. The next
crucial result is a not so obvious generalization of this fact.

\state 8.7. Proposition. If $S$ is an inverse semigroup, and $U_1$ and $U_2$ are up-sets in $\Gt{S}$, then so are
$U_1U_2$ and $U_1^{-1}$.

\Proof Take a pair of elements
  $$
  (\gamma_1, \gamma_2)\in(U_1\times U_2)\cap\Gt{S}\ex2,
  $$
  so that $\gamma_1\gamma_2$ is a typical element of $U_1U_2$. Supposing that $\delta$ is another element of $\Gt{S}$
such that $\gamma_1\gamma_2\leq\delta$, we must prove that $\delta\in U_1U_2$. According to (8.2), there exists $s\in
S$, as well as $\xi, \eta\in\tsp E$, such that $s^*s\in\xi\subseteq\eta$, and
  $$
  \gamma_1\gamma_2 = \gt s\xi, \and\delta= \gt s\eta.
  $$
  Let us moreover write $\gamma_i = \gt{s_i}{\xi_i}$, with $s_i^*s_i\in\xi_i$, for $i = 1,2$. Then necessarily $\xi_1 =
\theta_{s_2}(\xi_2)$, and
  $$
  \gt s\xi= \gamma_1\gamma_2 = \gt{s_1}{\xi_1}\gt{s_2}{\xi_2} = \gt{s_1s_2}{\xi_2}.
  $$
  As a consequence we have that $\xi= \xi_2$, and there exists $e\in\xi$, such that $se = s_1s_2e$. As
$\xi\subseteq\eta$, we deduce that $e\in\eta$, as well, whence
  $$
  \delta= \gt s\eta= \gt{se}\eta= \gt{s_1s_2e}{\eta} = \gt{s_1s_2}{\eta} = \gt{s_1}{\theta_{s_2}(\eta)}\gt{s_2}{\eta},
\eqno{(8.7.1)}
  $$
  where we are also using that
  $$
  s_2^*s_2\in\xi_2 = \xi\subseteq\eta,
  $$
  to make sure $\gt{s_2}{\eta}$ is a legitimate germ in $\Gt{S}$. The legitimacy of the germ
$\gt{s_1}{\theta_{s_2}(\eta)}$, above, is perhaps also under question, but we may remedy this as follows: noticing that
both $s^*s$ and $e$ lie in $\xi$, we see that
  $$
  \eta\supseteq\xi\ni s^*se = (se)^*se = (s_1s_2e)^*s_1s_2e = es_2^*s_1^*s_1s_2e \leq s_2^*s_1^*s_1s_2,
  $$
  so $s_2^*s_1^*s_1s_2$ belongs to $\eta$. Temporarily reverting to the character model for $\tsp E$, we then have that
  $$
  \theta_{s_2}(\eta)(s_1^*s_1) \explica{(2.17)}{=} \eta(s_2^*s_1^*s_1s_2) = 1,
  $$
  which, back to filters, means that $s_1^*s_1 \in\theta_{s_2}(\eta)$, validating the so far suspicious germ
above. Setting
  $$
  \zeta_1 = \gt{s_1}{\theta_{s_2}(\eta)}, \and\zeta_2 = \gt{s_2}{\eta},
  $$
  we then have that
  $$
  U_2\ni\gamma_2 = \gt{s_2}{\xi_2} = \gt{s_2}{\xi} \leq\gt{s_2}{\eta} = \zeta_2,
  $$
  so $\zeta_2$ lies in $U_2$, by virtue of $U_2$ being an up-set. Also
  $$
  U_1\ni\gamma_1 = \gt{s_1}{\xi_1} = \gt{s_1}{\theta_{s_2}(\xi_2)} = \gt{s_1}{\theta_{s_2}(\xi)}
\leq\gt{s_1}{\theta_{s_2}(\eta)} = \zeta_1,
  $$
  where we have used the easy to prove fact that $\theta_{s_2}$ is an order-preserving map. Therefore this shows that
$\zeta_1$ lies in $U_1$, whence
  $$
  \delta\explica{(8.7.1)}{=} \zeta_1\zeta_2\in U_1U_2,
  $$
  as needed. We leave the proof that $U_1^{-1}$ is an up-set for the reader. \endProof

We have already exploited the fact that products and inverses of slices are slices (in viewing $\Gt S\op$ as an inverse
semigroup), and now we know that products and inverses of up-sets are up-sets. Compactness is another property preserved
under products and inverses, so we see that the collection of subsets of $\Gt S$ which are, at the same time, compact,
up-sets, and slices, is closed under multiplication and inversion. In other words, this collection is an inverse
sub-semigroup of $\Gt S\op$.

\state 8.8. Definition. Let $S$ be an inverse semigroup.
  \item {(a)} A subset $U\subseteq\Gt S$ that is both a slice and an up-set relative to the spectral order will be
called an \emph{up-slice},
  \item {(b)} The inverse sub-semigroup of $\Gt S\op$ formed by all compact up-slices will be called the \emph{tight
envelope} of $S$, and it will be denoted by $\Cpl S$.

Since (8.6) implies that each $\Delta_s$ is a compact up-slice, it follows that the range of the tight regular
representation defined in (2.19) is contained in $\Cpl S$.

This will lead us to consider several inverse sub-semigroups of $\Gt S\op$ containing the range of the tight regular
representation of $S$, all of which share some special features which we would now like to discuss.

As we deal with these semigroups, especially when one is a sub-semigroup of the other, we need to take into account that
many of our concepts introduced above strongly depend on the environment in which they occur. For example, if $E$ is a
semilattice, if $F\subseteq E$ is a sub-semilattice, and if $C\subseteq F$ is a cover for some $f$ in $F$, as seen from
the perspective of $F$, then $C$ is not necessarily a cover for $f$ relative to $E$. Since the tight order relation is
based on the notion of covers, it also depends on the environment.

\state 8.9. Definition. An inverse semigroup $S$ shall be called \emph{flat} if its tight order ``$\menor$'' coincides
with the usual order, and consequently ``$\equiva$'' is tantamount to equality.

\state 8.10. Proposition. Given an inverse semigroup $S$, with idempotent semilattice denoted by $E$, let $\Tcal$ be any
inverse sub-semigroup of $\Gt S\op$, such that
  $$
  \Delta(S)\subseteq\Tcal\subseteq\Cpl S,
  $$
  and denote the idempotent semilattice\fn{Observe that the elements of $\Ecal$ are compact open up-sets contained in
the unit space of $\Gt S$, a.k.a. the tight spectrum $\tsp E$.} of\/ $\Tcal$\/ by $\Ecal$. Then
  \item {(i)} If $\{U_1,\cdots,U_n\}$ is a cover for some $U\in\Ecal$, then $U = \bigcup_{i = 1}^nU_i$,
  \item {(ii)} $\Tcal$ is flat.

\Proof Referring to our discussion above, throughout this proof, whenever we speak of covers, or of the tight order
relation, these notions should be interpreted, not within $\Cpl S$, but of course within $\Tcal$.

Regarding (i), suppose by contradiction that
  $$
  V := U \setminus\bigcup_{i = 1}^nU_i
  $$
  is nonempty. Since $V$ is an open subset of $\tsp E$, there exists some ultra-filter $\xi$ in $V$, and, by (2.13), we
may find some $e$ in $\xi$ such that $\xi\in\D_e\subseteq V$, and clearly
  $$
  \D_e \perp U_i,
  $$
  for every $i$. Using the crucial hypothesis that $\Tcal$ contains the range of the regular representation, we have
that
  $$
  \D_e = \Delta_e \in\Ecal,
  $$
  so the cover property implies that $\D_e \Cap U_i$, for some $i$, a contradiction, hence concluding the proof of (i).

As for (ii), our task is then to prove that, if $U_1$ and $U_2$ are members of $\Tcal$, then
  $$
  U_1\menor U_2 \reliff U_1\subseteq U_2.
  $$
  We will only prove the forward implication, because the converse one is trivial. Assuming first that $U_1$ and $U_2$
are idempotent elements in $\Tcal$, and arguing by contradiction, suppose that
  $$
  V := U_1 \setminus U_2
  $$
  is nonempty. Since $V$ is then an open subset of $\tsp E$, there exists some ultra-filter $\xi$ in $V$, and, by
(2.13), we may find some $e$ in $\xi$ such that $\xi\in\D_e\subseteq V$. This may then be reinterpreted as saying that
  $$
  U_2 \perp\D_e \leq U_1.
  $$
  Observing that $\D_e \in\Ecal$, as seen above, this contradicts the fact $U_1\menor U_2$, and hence proves our claim
that $U_1\subseteq U_2$. In the general case, no longer assuming that $U_1$ and $U_2$ are idempotents, we have that $
U_1^*U_1 \menor U_2^*U_2, $ by (3.9), so we deduce from the above that $U_1^*U_1 \subseteq U_2^*U_2$.

From the assumption that $U_1\menor U_2$, we have that there exists a finite cover $C$ for $U_1^*U_1$, such that $U_1U =
U_2U$, for all $U$ in $C$. This said, choose any $\gamma$ in $U_1$, so that
  $$
  d(\gamma) = \gamma^{-1}\gamma\in U_1^*U_1.
  $$
  Using (i) we have that $U_1^*U_1$ coincides with the union of all members of $C$, and hence that $d(\gamma)$ lies in
$U$, for some $U$ in $C$. From this we deduce that
  $$
  \gamma= \gamma d(\gamma) \in U_1U = U_2U\subseteq U_2,
  $$
  thus proving that $U_1\subseteq U_2$. \endProof

Most of the time, restricting a function consists in reducing its domain but, of course, one can also reduce the
codomain, as long as the new codomain is known to contain the range. Obviously nobody in their right mind would ever try
to do something like that to a surjective function, but, when speaking of tightly surjective ones, this doesn't only
make sense -- it raises a lot of questions! On the other hand, when restricting the codomain of an injective function,
nobody will ever give a second though about the injectivity of the new function, but, when speaking of tightly injective
ones, this is not so straightforward.

\state 8.11. Proposition. Given an inverse semigroup $S$, let $\Tcal\subseteq\Cpl S$ be any inverse sub-semigroup
containing the range of the regular representation. Also let
  $$
  h:S \to\Tcal
  $$
  be the tight regular representation of $S$, seen as a map into $\Tcal$. In other words, $h$ is the
codomain-restriction of the tight regular representation to $\Tcal$. Then $h$ is a consonance. In case $\Tcal= \Cpl S$,
we will denote $h$ by $\rho$, and call it the \emph{fundamental consonance} of $S$. Consequently $\Gt{S}$, $\Gt{\Tcal}$,
and $\Gt{\Cpl S}$ are all isomorphic as \otgs/.

\begingroup\noindent\hfill\beginpicture\setcoordinatesystem units <0.025truecm, 0.025truecm> \setplotarea x from -10 to
130, y from 150 to -20 \put{\null} at -10 150 \put{\null} at -10 -20 \put{\null} at 130 150 \put{\null} at 130 -20
\put{$\Delta(S)$} at 120 0 \put{$\Tcal$} at 120 46.667 \put{\rotatebox{90}{$\subseteq$}} at 120 23.333 \put{$\Cpl{S}$}
at 120 93.333 \put{\rotatebox{90}{$\subseteq$}} at 120 70 \put{$\Gt{S}\op$} at 120 140 \put{\rotatebox{90}{$\subseteq$}}
at 120 116.667 \put{$S$} at 0 46.667 \arrow<0.11cm> [0.5, 1.8] from 14.208 57.718 to 100.266 124.651 \put{$\Delta$} at
65.854 115.622 \arrow<0.11cm> [0.5, 1.8] from 16.776 53.191 to 103.224 86.809 \put{$\rho$} at 73.665 86.044
\arrow<0.11cm> [0.5, 1.8] from 18 46.667 to 102 46.667 \put{$h$} at 76.8 58.667 \endpicture\hfill\null\endgroup

\Proof Focusing on the proof of the tight injectivity of $h$, suppose that $s_1$ and $s_2$ are elements in $S$ such that
$h(s_1)\equiva h(s_2)$, that is, $\Delta_{s_1}\equiva\Delta_{s_2}$. Using (8.10) we then have that $\Delta_{s_1} =
\Delta_{s_2}$, so (3.13) implies that $s_1\equiva s_2$, as required. Regarding tight surjectivity, given any $U$ in
$\Tcal$, we have that $U$ is a compact up-slice, so (8.6) allows us to write
  $$
  U = \bigcup_{i\in I}\Delta_{s_i},
  $$
  for some collection $\{s_i\}_{i\in I}$ of elements in $S$, which may be assumed finite due to the fact that $U$ is
compact. With an eye on (3.10.a), we deduce that $\Delta_{s_i}\menor U$, because actually $\Delta_{s_i}\subseteq U$, and
clearly
  $$
  U^*U = d(U) = \bigcup_{i\in I}d(\Delta_{s_i}) = \bigcup_{i\in I}\Delta_{s_i^*}\Delta_{s_i},
  $$
  which we may reinterpret as saying that
  $$
  \big\{\Delta_{s_i^*}\Delta_{s_i}U^*U: i\in I \big\}
  $$
  is a cover\fn{Actually $\Delta_{s_i^*}\Delta_{s_i}U^*U = \Delta_{s_i^*}\Delta_{s_i}$ but we have written this as such
in order to fit the pattern in (3.10.a).} for $U^*U$. This proves that $h$ is tightly surjectivity while the last
sentence of the statement follows from (8.5). \endProof

\section 9 Consonance as an equivalence relation

\firstLine Having defined the notion of consonance for homomorphisms between inverse semigroups in (7.1), one could
attempt to say that two inverse semigroups $S$ and $T$ are \emph{consonant} if there exists a consonance from $S$ to
$T$. However, even though the term ``consonance'' might evoke an impression of equivalence, it is obvious that this is
not an equivalence relation.

On the other hand, one of our main results, namely Theorem (8.5), implies that consonant inverse semigroups have
isomorphic tight groupoids, and, of course, isomorphism is a bona fide equivalence relation.

In order to get out of this conundrum, we shall abandon the above attempt at defining consonant inverse semigroups and
instead adopt the following:

\state 9.1. Definition. Given inverse semigroups $S$ and $T$, we shall say that $S$ and $T$ are \emph{consonant} if
there exists a sequence of inverse semigroups
  $$
  S = S_0,\ S_1, \ S_2,\ \cdots,\ S_n = T,
  $$
  such that, for every $i\geq1$, either there exists a consonance from $S_{i-1}$ to $S_i$, or one from $S_i$ to
$S_{i-1}$.

In other words, this is the equivalence relation generated by the abandoned notion of consonant inverse semigroups
above.

Since isomorphism of \otgs/ is an equivalence relation, Theorem (8.5) implies that the tight groupoids of consonant
inverse semigroups are isomorphic.

The trouble, however, is that equivalence relations generated by non-symmetric relations, as the one above, often
present an uncomfortable situation as the length of the ``zig-zag'' in (9.1) may be unpredictable. Shortening this
zig-zag is therefore the main goal of this section.

\medskip We shall start by considering the following specific situation: suppose that $S_1$ and $S_2$ are inverse
semigroups such that $\Gt{S_1}$ is isomorphic to $\Gt{S_2}$, in the sense that there exists an isomorphism of \otgs/
  $$
  \varphi:\Gt{S_1} \to\Gt{S_2},
  $$
  which is however not assumed to be of the form $\ch$, for any consonance $h$ from $S_1$ to $S_2$. Elementary examples
show that this situation is indeed possible, so the question arises as to what extent are $S_1$ and $S_2$ related as a
consequence of their tight groupoids being isomorphic.

Since the situation is now perfectly symmetric in terms of $S_1$ and $S_2$, we would rather not prioritize one over the
other, so we will choose a \otg/ $G$ in the isomorphism class of either $\Gt{S_1}$ or $\Gt{S_2}$, and we will let
  $$
  \pointfun_1:\Gt{S_1} \to G, \and\pointfun_2:\Gt{S_2} \to G, \eqno{(9.2)}
  $$
  be any two chosen isomorphisms. One possibility will be to take $G = \Gt{S_2}$, with $\pointfun_1 = \varphi$, and
$\pointfun_2$ being the identity map, but, even if we do so, we will still think of $G$ as an abstract object.

We will then consider, for each $i = 1,2$, the map
  $$
  h_i:S_ i \to G\op,
  $$
  given by
  $$
  h_i(s) = \pointfun_i(\Delta^i_s), \for s\in S_ i, \eqno{(9.3)}
  $$
  where $\Delta_s^i$ is as in (3.11) relative to $S_i$. Should we be concerned with just one of the $S_ i$'s, it might
be a good idea to tacitly identify $G$ and $\Gt{S_i}$, hence blurring $\pointfun_i$ out of the picture. So,
e.g.~$h_i(s)$ would be thought of as $\Delta^i_s$, in which case $h_i$ becomes an instance of the tight regular
representation introduced in (2.19), en passant proving that $h_i$ is an inverse semigroup homomorphism.

A situation in which the $\pointfun_i$ cannot possibly be ignored is as follows.

\state 9.4. Proposition. Let $\Tcal$ be the inverse sub-semigroup of $G\op$ generated by the union of the ranges of
$h_1$ and $h_2$, and, for every $i = 1,2$, let
  $$
  h_i':S_ i\to\Tcal
  $$
  be the same map as $h_i$, with codomain restricted to $\Tcal$. Then $h_i'$ is a consonance.

\Proof For every $i = 1, 2$, let us consider the map\fn{We are all used to the abuse of language of denoting by $f(A)$
the image of a set $A$ under a function $f$ whose domain \emph{contains} $A$, even though $A$ might not be an
\emph{element} of the domain of $f$. We could therefore also write $\setfun_i(U) = \pointfun_i(U)$, but, as we believe
this abuse of language might cause some confusion, we are introducing the extra notation $\setfun_i$ to refer to the
\emph{set-function} associated to the \emph{point-function} $\pointfun_i$.}
  $$
  \setfun_i:\Gt{S_i}\op\to G\op,
  $$
  given by
  $$
  \setfun_i(U) = \big\{\pointfun_i(\gamma) : \gamma\in U\big\}, \for U \in\Gt{S_i}\op.
  $$

Since $\pointfun_i$ is an isomorphism of \otgs/, we have that $\setfun_i$ is an up-set preserving isomorphism of inverse
semigroups of slices.

We may then redefine $h_i$ as
  $$
  h_i = \setfun_i\circ\Delta^i,
  $$
  and, besides reinforcing the fact that $h_i$ is an inverse semigroup homomorphism, we deduce from (8.6) that the range
of $h_i$ consists of compact up-slices.

Even though we haven't proved (8.7) for general \otgs/, this result does apply to $G$, because $G$ is in all respects
indistinguishable from $\Gt{S_i}$, where (8.7) holds. This said, we conclude that $\Tcal$ consists of compact up-slices
of $G$.

Observing that the situation is symmetric in terms of $S_1$ and $S_2$, we need only prove that $h_1'$ is a
consonance. Letting
  $$
  \Tcal_1 = \setfun_1^{-1}(\Tcal),
  $$
  we have that $\Tcal_1$ is an inverse sub-semigroup of $\Gt{S_1}\op$ consisting of compact up-slices, and clearly
containing the range of $\setfun_1^{-1}\circ h_1 = \Delta^1$. So we deduce that
  $$
  \Delta^1(S_1) \subseteq\Tcal_1\subseteq\Cpl{S_1} \subseteq\Gt{S_1}\op.
  $$
  Several conclusions follow: firstly, as $\Tcal_1$ contains the range of $\Delta^1$, we can restrict the codomain of
$\Delta^1$ to $\Tcal_1$, producing a homomorphism which we will call $\sigma_1$.

\begingroup\noindent\hfill\beginpicture\setcoordinatesystem units <0.025truecm, 0.025truecm> \setplotarea x from -30 to
210, y from 130 to -45 \put{\null} at -30 130 \put{\null} at -30 -45 \put{\null} at 210 130 \put{\null} at 210 -45
\put{$\Tcal_1$} at 0 0 \put{$\Cpl{S_1}$} at 0 50 \put{\rotatebox{90}{$\subseteq$}} at 0 25 \put{$\Gt{S_1}\op$} at 0 100
\put{\rotatebox{90}{$\subseteq$}} at 0 75 \put{$\Tcal$} at 180 0 \put{$G\op$} at 180 100
\put{\rotatebox{90}{$\subseteq$}} at 180 50 \arrow<0.11cm> [0.5, 1.8] from 37 100 to 158 100 \put{$\setfun_1$} at 97.5
114 \put{$S_1$} at 90 50 \arrow<0.11cm> [0.5, 1.8] from 70.769 60.684 to 19.231 89.316 \put{$\Delta^1$} at 35.624 64.193
\arrow<0.11cm> [0.5, 1.8] from 109.231 60.684 to 160.769 89.316 \put{$h_1$} at 144.376 64.193 \arrow<0.11cm> [0.5, 1.8]
from 109.231 39.316 to 160.769 10.684 \put{$h_1'$} at 143.405 34.058 \arrow<0.11cm> [0.5, 1.8] from 22 0 to 158 0
\put{$\setfun_1|_{\Tcal_1}$} at 96.8 -18 \setdashes<1.5pt> \arrow<0.11cm> [0.5, 1.8] from 70.769 39.316 to 19.231 10.684
\put{$\sigma_1$} at 36.595 34.058 \endpicture\hfill\null\endgroup

\noindent Secondly, observing that $\setfun_1\circ\sigma_1 = h_1'$, and that $\setfun_1$ restricts to an isomorphism
from $\Tcal_1$ to $\Tcal$, our goal of proving that $h_1'$ is a consonance will be achieved once we prove that
$\sigma_1$ has this property, but this follows immediately from (8.11). \endProof

This brings us to the main result of this section.

\state 9.5. Theorem. Let $S_1$ and $S_2$ be inverse semigroups (with zero). Then the following are equivalent:
  \item {(i)} $\Gt{S_1}$ and $\Gt{S_2}$ are isomorphic as \otgs/,
  \item {(ii)} $S_1$ and $S_2$ are consonant,
  \item {(iii)} there exists an inverse semigroup $T$, and consonances
  $$
  S_ 1 \labelarrow{h_1} T \ {\buildrel\textstyle h_2 \over{\hbox to 24pt {\leftarrowfill}}} \ S_2.
  $$

\Proof Implication (i$\Rightarrow$iii) is the content of (9.4), while (iii$\Rightarrow$ii) is obvious. Finally,
(ii$\Rightarrow$i) follows from a comment made shortly after (9.1). \endProof

In conclusion we see that, when two inverse semigroups are consonant, the zig-zag mentioned in the definition needs not
consist of more than one ``zig'' and one ``zag''.

\section 10 Tight envelope of inverse semigroups

\firstLine We would now like to conduct a careful study of the tight envelope of an inverse semigroup $S$, defined in
(8.8.b). The first relevant result is as follows:

\state 10.1. Theorem. Let $S_1$ and $S_2$ be inverse semigroups. Then the following are equivalent:
  \item {(i)} $S_1$ and $S_2$ are consonant,
  \item {(ii)} $\Cpl{S_1}$ and $\Cpl{S_2}$ are isomorphic as inverse semigroups.

\Proof Assuming that $S_1$ and $S_2$ are consonant, we have by (9.5) that there exists an isomorphism of \otgs/
  $$
  \pointfun:\Gt{S_1} \to\Gt{S_2}.
  $$
  The corresponding map on the respective inverse semigroups of slices, namely
  $$
  \setfun:U\in\Gt{S_1}\op\mapsto\pointfun(U)\in\Gt{S_2}\op
  $$
  is therefore clearly a semigroup isomorphism, which moreover sends compact up-slices on the left-hand-side to compact
up-slices in the right-hand-side and vice-versa. In other words, $\setfun$ restricts to an isomorphism between
$\Cpl{S_1}$ and $\Cpl{S_1}$.

Conversely, assuming (ii), it is obvious that $\Gt{\Cpl{S_1}}$ and $\Gt{\Cpl{S_2}}$ are isomorphic as \otgs/. By (8.11)
we then have that $\Gt{S_1}$ and $\Gt{S_2}$ are also isomorphic as \otgs/. Applying (9.5), we then get (i). \endProof

A useful consequence is as follows:

\state 10.2. Corollary. If $S_1$ and $S_2$ are consonant inverse semigroups, then there is a consonance
$h:S_2\to\Cpl{S_1}$. In other words, the inverse semigroup $T$ and the consonances $h_1$ and $h_2$ referred to in
(9.5.iii) can be taken to be
  $$
  S_ 1 \labelarrow{{\scriptstyle\rho_1}} \Cpl{S_1} \ {\buildrel h \over{\hbox to 24pt {\leftarrowfill}}} \ S_2,
  $$
  where $\rho_1$ is the fundamental consonance of $S_1$.

\Proof Using (10.1), let $\varphi:\Cpl{S_2} \to\Cpl{S_1}$ be any isomorphism, and let $\rho_2$ be the fundamental
consonance of $S_2$. Then $h = \varphi\circ\rho_2$ clearly satisfies the required conditions. \endProof

The next result highlights the fact that any consonance is a factor of the fundamental one.

\state 10.3. Theorem. Let $S_1$ and $S_2$ be inverse semigroups and let $h:S_1\to S_2$ be a consonance. Then the
fundamental consonance $\rho_1$ of $S_1$ factors as $\rho_1 = k\circ h$, for some homomorphism $k :S_2\to\Cpl{S_1} $
which is also a consonance.

\begingroup\noindent\hfill\beginpicture\setcoordinatesystem units <0.025truecm, -0.02truecm> \setplotarea x from -10 to
160, y from -10 to 100 \put{\null} at -10 -10 \put{\null} at -10 100 \put{\null} at 160 -10 \put{\null} at 160 100
\put{$S_1$} at 0 0 \put{$\Cpl{S_1}$} at 150 0 \arrow<0.11cm> [0.5, 1.8] from 15 0 to 135 0 \put{$\rho_1$} at 75 -18
\put{$S_2$} at 75 80 \arrow<0.11cm> [0.5, 1.8] from 13.5 14.4 to 61.5 65.6 \put{$h$} at 28.016 48.891 \arrow<0.11cm>
[0.5, 1.8] from 88.5 65.6 to 136.5 14.4 \put{$k$} at 121.984 48.891 \endpicture\hfill\null\endgroup

\Proof We have that $\ch$ is an isomorphism by (8.5), so the map $\hzao$ referred to in (6.4) sends compact up-slices to
compact up-slices, and so does $\hzao^{-1}$. Consequently $\hzao(\Cpl{S_1}) = \Cpl{S_2}$. Letting $\varphi$ be the map
obtained by restricting the domain of $\hzao$ to $\Cpl{S_1}$, as well as its codomain to $\Cpl S_2$, we get the
following updated version of Diagram (6.5):

\begingroup\noindent\hfill\beginpicture\setcoordinatesystem units <0.025truecm, 0.025truecm> \setplotarea x from -40 to
150, y from 180 to -40 \put{\null} at -40 180 \put{\null} at -40 -40 \put{\null} at 150 180 \put{\null} at 150 -40
\put{$S_1$} at 0 140 \put{$S_2$} at 110 140 \arrow<0.11cm> [0.5, 1.8] from 18 140 to 92 140 \put{$h$} at 55 155
\put{$\Cpl{S_1}$} at 0 56 \put{$\Cpl{S_2}$} at 110 56 \arrow<0.11cm> [0.5, 1.8] from 18 56 to 92 56 \put{$\varphi$} at
55 41 \put{$\Gt{S_1}\op$} at 0 0 \put{$\Gt{S_2}\op$} at 110 0 \arrow<0.11cm> [0.5, 1.8] from 35 0 to 75 0 \put{$\hzao$}
at 55 -15 \arrow<0.11cm> [0.5, 1.8] from 95.694 129.076 to 14.306 66.924 \put{$k$} at 45.896 109.922 \arrow<0.11cm>
[0.5, 1.8] from 0 124 to 0 72 \put{$\rho_1$} at -15 98 \arrow<0.11cm> [0.5, 1.8] from 110 124 to 110 72 \put{$\rho_2$}
at 125 98 \put{\rotatebox{270}{$\subseteq$}} at 0 28 \put{\rotatebox{270}{$\subseteq$}} at 110 28
\endpicture\hfill\null\endgroup

\noindent in which both $\hzao$ and $\varphi$ are isomorphisms of inverse semigroups. The desired homomorphism $k$ may
then be taken to be $\varphi^{-1}\circ\rho_2$. It is a consonance because so is $\rho_2$, and $\varphi$ is an
isomorphism. \endProof

From now on we would like to analyze two special properties of the tight envelope $\Cpl S$, with the goal of obtaining
an abstract characterization of it. The first one of these, as we already know by (8.10), is flatness. Regarding the
second property we first need some preparations.

Given an inverse semigroup $S$, recall from \cite{lawsonBook} that two elements $s$ and $t$ in $S$ are said to be
\emph{compatible} if $st^*$ and $s^*t$ are idempotent. The following gives other useful characterizations of
compatibility:

\state 10.4. Proposition. Given an inverse semigroup $S$, and given $s$ and $t$ in $S$, let $ e = s^*st^*t, \and f =
ss^*tt^*. $ Then the following are equivalent:
  \item {(i)} $s$ and $t$ are compatible,
  \item {(ii)} $se = te$, and $ft = fs$,
  \item {(iii)} $st^*t = ts^*s$, and $ss^*t = tt^*s$.

\Proof \proofImply {i}{iii} Assuming (i), let $p = st^*$. Then
  $$
  ts^*s = tt^*ts^*s = ts^*st^*t = p^*pt = pt = st^*t,
  $$
  and one may similarly prove that $tt^*s = ss^*t$.

\Medskip\proofImply{iii}{i} Assuming (iii) we have
  $$
  st^*st^* = st^*{st^*t}\, t^* = st^*ts^*st^* = st^*(st^*)^*st^* = st^*,
  $$
  and one may similarly prove that $s^*t$ is idempotent. \Medskip(ii$\Leftrightarrow$iii)\enspace Left to the reader.

\endProof

Perhaps the most intuitive characterization of compatibility is (10.4.ii), especially within the symmetric inverse
semigroup ${\cal I}(X)$ for a given set $X$. In this case, among other things, the element $e$ mentioned in (10.4.ii)
corresponds to the intersection of the domains of $s$ and $t$, and the condition $se = te$ says that $s$ and $t$
coincide on their common domain.

One of the easiest ways to prove that two elements $s$ and $t$ are compatible is to find a third element $r$ such that
$s\leq r$, and $t\leq r$. This is because
  $$
  st^*t = rs^*st^*t = rt^*ts^*s = ts^*s,
  $$
  and similarly $tt^*s = ss^*t$, verifying (10.4.ii). This is however not to say that compatible elements are always
dominated by a single element, as elementary examples show.

Given a finite family $\{s_1,\ldots,s_n\}$ of pairwise compatible elements in an inverse semigroup $S$, one says that
$t$ is their \emph{join}, if $t$ is the smallest upper bound of the $s_i$, relative to the usual order of $S$. In this
case $t$ is denoted by
  $$
  t = \bigvee_{i = 1}^n s_i.
  $$

\state 10.5. Definition. One says that the inverse semigroup $S$ \emph{has finite joins} if all finite, pairwise
compatible families admit a join. If moreover the multiplication in $S$ distributes over all existing finite joins, then
$S$ is said to be \emph{distributive}.

\state 10.6. Lemma. Let $S$ be a flat distributive inverse semigroup. Also suppose that we are given elements $s, s_1,
\ldots, s_n$ in $S$, such that each $s_i\leq s$. Then the following are equivalent:
  \item {(i)} $\{s_1^*s_1, \ldots, s_n^*s_n\}$ is a cover for $s^*s$,
  \item {(ii)} $ s = \bigvee_{i = 1}^n s_i. $

\Proof As observed above, the fact that the $s_i$ have a common upper bound implies that they are pairwise
compatible. So we let
  $$
  r = \bigvee_{i = 1}^n s_i,
  $$
  observing that necessarily $r\leq s$. Assuming (i) we then claim that $r^*r = s^*s$. In order to prove this, we notice
that there is no nonzero idempotent $g$ such that
  $$
  r^*r\perp g\leq s^*s,
  $$
  because $g$ would also be orthogonal to the $s_i^*s_i$, violating the properties of the cover in (i). This said, there
is no witness denying that $s^*s\leq r^*r$, so we conclude that $s^*s\menor r^*r$. Since $S$ is flat, it follows that
$s^*s\leq r^*r$. As it is obvious that $r^*r\leq s^*s$, the conclusion is that $r^*r = s^*s$. Therefore
  $$
  r = sr^*r = ss^*s = s,
  $$
  verifying (ii). Conversely, assuming (ii), let $g$ be a nonzero idempotent with $g\leq s^*s$. Then, using the
distributive property,
  $$
  0\neq g = gs^*s = gs^* \bigvee_{i = 1}^n s_i = \bigvee_{i = 1}^n gs^*s_i = \bigvee_{i = 1}^n gs_i^*s_i,
  $$
  where we have used that $s^*s_i = s^*s_is_i^*s_i = s_i^*s_i$. Therefore $g\Cap s_i^*s_i$, for some $i$, and we are
done. \endProof

\state 10.7. Proposition. The tight envelope of every inverse semigroup is distributive.

\Proof Given an inverse semigroup $S$, recall that $\Cpl S$ consists of all compact up-slices of the groupoid $\Gt
S$. One may then easily prove that a finite family $\{U_1,\ldots,U_n\}$ of elements of $\Cpl S$ is pairwise compatible
if and only if the union of all the $U_i$ is a slice, in which case said union is a compact up-slice simply because the
union of up-sets is always an up-set, and the $U_i$ are all compact. Moreover, as the order among slices is simply the
order of inclusion, we have that
  $$
  \bigvee_{i = 1}^n U_i = \bigcup_{i = 1}^n U_i,
  $$
  so $\Cpl S$ has finite joins. The fact that it is also distributive is now easily proved by inspection. \endProof

The universal characterization of the tight envelope is then as follows.

\state 10.8. Theorem. Let $S$ and $T$ be consonant inverse semigroups. Then the following hold:
  \item {(i)} if\/ $T$ is flat, then $T$ is isomorphic to an inverse sub-semigroup of $\Cpl S$,
  \item {(ii)} if\/ $T$ is flat and distributive, then $T$ is isomorphic to $\Cpl S$.

\Proof Consider the diagram
  $$
  S \labelarrow{{\scriptstyle\rho}} \Cpl{S} \ {\buildrel h \over{\hbox to 24pt {\leftarrowfill}}} \ T,
  $$
  arising from the application of (10.2), where $\rho$ is the fundamental consonance of $S$. Assuming that $T$ is flat,
we claim that $h$ is injective. To see this, observe that, for $t_1$ and $t_2$ in $ T$, one has that
  $$
  h(t_1) = h(t_2) \imply h(t_1) \equiva h(t_2) \mathrel{\buildrel(*) \over\Longrightarrow} t_1 \equiva t_2
\mathrel{\buildrel(**) \over\Longrightarrow} t_1 = t_2,
  $$
  where the implication marked with $(*)$ is due to the tight injectivity of $h$, while the one marked with $(**)$
follows from the flatness of $T$. We then have that $h$ is an isomorphism onto its image, proving (i).

\def\setin#1{\big\{\, #1\, \big\}_{i = 1}^n}

Assuming that $T$ is flat and distributive, we will prove that $h$ is also surjective. For this, pick any $U$ in $\Cpl
S$. Since $h$ is tightly surjective, there exists a finite set $\{t_1, \ldots, t_2\}\subseteq T$, such that
$h(t_i)\menor U$, for every $i$, and
  $$
  C := \setin{h(t_i^*t_i)\, U^*U}
  $$
  is a cover for $U^*U$. Since $\Cpl S$ is flat, we have that $h(t_i)\leq U$, for every $i$. This implies that
$h(t_i^*t_i)\leq U^*U$, so actually
  $$
  C = \setin{h(t_i^*t_i)},
  $$
  and we deduce from the fact that $\Cpl S$ is flat and distributive, and from (10.6) that
  $$
  U = \bigvee_{i = 1}^n h(t_i).
  $$
  This obviously also entails that the $h(t_i)$ are pairwise compatible. In particular, for all $i$ and $j$, one has
that $h(t_i^*t_j)$ is idempotent which means no more no less than
  $$
  h(t_i^*t_j) = h(t_i^*t_jt_i^*t_j).
  $$
  The injectivity of $h$, given by (i), then implies that $t_i^*t_j = t_i^*t_jt_i^*t_j$. By the same method one proves
that $t_it_j^*$ is also idempotent, hence showing that the $t_i$ are pairwise compatible. Since $T$ has finite joins by
hypothesis, we may let
  $$
  t = \bigvee_{i = 1}^n t_i,
  $$
  and the proof will be finished once we prove that $h(t) = U$.

In order to verify this claim we first notice that each $h(t_i)\leq h(t)$, and that $\setin{t_i^*t_i}$ is a cover for
$t^*t$ by (10.6). Therefore $\setin{h(t_i^*t_i)}$ is a cover for $h(t^*t)$ by (5.5), and then Lemma (10.6) once more
comes to our rescue giving
  $$
  h(t) = \bigvee_{i = 1}^n h(t_i) = U,
  $$
  proving the claim and finishing the proof. \endProof

Last but not least, we have the following main result:

\state 10.9. Corollary. Up to isomorphism, the tight envelope of an inverse semigroup $S$ is the only flat distributive
inverse semigroup that is consonant to $S$.

\section 11 Tight-like topological spaces

\def\cli#1{[\, #1, \infty)} \caldef E

\firstLine Up to this point, virtually all topological spaces and groupoids treated in this work arouse from inverse
semigroups, including, of course, semilattices. However, from now on we would like to study topological spaces and
groupoids from an abstract standpoint, but still under the aegis of inverse semigroups. By this we mean that we will
introduce abstract classes of topological spaces and groupoids by requiring axioms inspired by the properties satisfied
by tight spectra and tight groupoids. The present section will deal with topological spaces and the next one, with
groupoids. As before, order relations will be paramount.

\state 11.1. Definition. A \emph{tight-like space} is a locally compact, Hausdorff topological space $X$ equipped with a
partial order relation ``$\leq$'' such that:
  \item {(a)} for every $x$ in $X$, there exists some compact open up-set containing $x$,
  \item {(b)} for every $x$ and $y$ in $X$, such that $y\not\geq x$, there exists a compact open up-set $U$, such that
$x\in U$, but $y\notin U$,
  \item {(c)} the maximal elements are dense in $X$.

\bigskip The existence of an up-set $U$, as in (11.1.b), is clearly a \emph{sufficient} condition for $y\not\geq x$. On
the other hand, this condition requires that this be a \emph{necessary} condition.

If $x$ and $y$ are distinct points of $X$, assume without loss of generality that $y\not\geq x$. Then the set $U$
provided by (11.1.b) splits $X$ as a disjoint union of two open sets
  $$
  X = U\cup(X\setminus U)
  $$
  leaving $x$ and $y$ in different sides, so they cannot be in the same connected component of $X$. In other words, $X$
is totally disconnected. Being locally compact, the topology of $X$ necessarily admits a basis formed by compact open
sets.

The motivation for the terminology ``tight-like space'' is the following:

\state 11.2. Proposition. If $E$ is a semilattice, then $\tsp E$ is a tight-like space.

\Proof Given $\xi$ in $\tsp E$, and choosing any $e\in\xi$, we have that $\xi$ lies in the compact open up-set $D_e$,
proving (a). Now suppose that $\xi$ and $\eta$ are filters in $\tsp E$, such that $\eta\not\supseteq\xi$. Then there
exists some $e\in\xi$, such that $e\notin n$. Therefore $\xi\in D_e$, while $\eta\notin D_e$, proving (b). Finally (c)
holds due to (2.8). \endProof

\fix From now on we fix a tight-like space $X$.

\state 11.3. Lemma. Given $x$ in $X$, define
  $$
  \cli x = \{y\in X: y\geq x\},
  $$
  and
  $$
  \xi_x = \{V\subseteq X: \text{$V$ is a compact open up-set such that } x\in V\}.
  $$
  Then
  $$
  \bigcap_{V\in\xi_x}V = \cli x.
  $$

\Proof Given $y\in\bigcap_{V\in\xi_x}V$, suppose by contradiction that $y\notin\cli x$, meaning that $y\not\geq x$. By
(11.1.b) there exists a compact open up-set $U$, such that $x\in U$, but $y\notin U$. Such an $U$ is then a member of
$\xi_x$, so $y\in U$, a contradiction. This proves that $\bigcap_{V\in\xi_x}V \subseteq\cli x$, while the reverse
inclusion is obvious. \endProof

Our next goal is to prove the following crucial technical tool supporting Theorem (11.8), below, the main result of this
section.

\state 11.4. Lemma. Let $U$ and\/ $V$ be compact open up-sets, such that $V\varsubsetneq U$. Then there exists a
nonempty compact open up-set $W\subseteq U\setminus V$.

\Proof By (11.1.c) we may choose a maximal element $x$ in $U\setminus V$. With $\xi_x$ as in (11.3), we then claim that
there exists some $W$ in $\xi_x$ satisfying $V\cap W = \varnothing$. Arguing by contradiction, suppose that
  $$
  V\cap W\neq\varnothing,
  $$
  for all $W\in\xi_x$. Given $W_1$ and $W_2$ in $\xi_x$, observe that $W_1\cap W_2$ is also in $\xi_x$, so
  $$
  (V\cap W_1)\cap(V\cap W_2) = V\cap(W_1\cap W_2) \neq\varnothing,
  $$
  and we see that the family of sets of the form $V\cap W$, with $W$ in $\xi_x$, satisfies the finite intersection
property. Since these are all compact sets, the intersection of the whole family is nonempty, so we may choose some
  $$
  y\in\bigcap_{W\in\xi_x}V\cap W = V\cap\Big(\bigcap_{W\in\xi_x}W\Big) \explica{(11.3)}{=} V\cap\cli x.
  $$
  Therefore $y\in V$, and $y\geq x$, but since $x$ is maximal, we deduce that $y = x$, so
  $$
  V\ni y = x\in U\setminus V,
  $$
  a contradiction. This proves the claim, so we may pick some $W$ in $\xi_x$ such that $V\cap W = \varnothing$. The set
$W' = W\cap U$, is then nonempty because it contains $x$, and clearly $W'\subseteq U\setminus V$, concluding the
proof. \endProof

\state 11.5. Definition. We shall denote the set of all compact open up-sets of $X$ by $\E$.

Observing that $\E$ is closed under intersections, we see that $\E$ is a semilattice, with the empty set playing the
role of zero.

\state 11.6. Proposition. Let $U,U_1,U_2,\ldots,U_n\in\E$, be such that $U_i\subseteq U$, for every $i$. Then
  $$
  \{U_1,U_2,\ldots,U_n\}
  $$
  is a cover for $U$, in the semilattice sense, if and only if
  $$
  U = \bigcup_{i = 1}^n U_i.
  $$

\Proof The ``if'' part being trivial, we focus on the ``only if'' part. For this, suppose by contradiction that the
above is a cover in the semilattice sense and yet $\bigcup_{i = 1}^n U_i\varsubsetneq U$. It is then easy to see that
  $$
  V := \bigcup_{i = 1}^n U_i
  $$
  is a compact open up-set, so we may apply (11.4) to produce some nonzero (nonempty) $W$ in $\E$ such that $W\leq U$,
and $W\perp U_i$, for all $i$, contradicting the cover property. \endProof

\state 11.7. Proposition. For every $x$ in $X$, one has that $\xi_x$ is a tight filter.

\Proof We must prove that, whenever $\{U_1,U_2,\ldots,U_n\}$ is a cover for some $U$ in $\E$, and $U\in\xi_x$, then
there exists some $i$ such that $U_i\in\xi_x$. To say that $U\in\xi_x$ is to say that $x\in U$, so (11.6) implies that
$x\in U_i$, for some $i$, whence $U_i\in\xi_x$. \endProof

We have already seen that every tight spectrum is a tight-like space. Now we prove the converse.

\state 11.8. Theorem. Every tight-like space is order-homeomorphic to the tight spectrum of a semilattice. More
precisely, for every tight-like space $X$, the map
  $$
  \varphi:x\in X\mapsto\xi_x\in\tsp\E
  $$
  is an order-homeomorphism, where $\E$ is the semilattice formed by all compact open up-sets of $X$.

\Proof Let us first prove that $\varphi$ is continuous. For this it is enough to prove that, for each $U$ in $\E$, the
map
  $$
  x\in X\mapsto\bool{U\in\xi_x}\in\{0,1\}
  $$
  is continuous. Now, since
  $$
  \bool{U\in\xi_x} = \bool{x\in U},\for x\in X,
  $$
  we see that the above map is nothing but the characteristic function of $U$, which is clearly continuous.

We next prove that the range of $\varphi$ is closed. For this suppose that $\{x_i\}_i$ is a net in $X$, such that
$\{\varphi(x_i)\}_i$ converges to some $\xi$ in $\tsp\E$. Since filters are nonempty sets, we may pick some $U$ in
$\xi$. It then follows that
  $$
  1 = \bool{U\in\xi} = \lim_i \bool{U\in\varphi(x_i)} = \lim_i \bool{U\in\xi_{x_i}} = \lim_i \bool{x_i\in U},
  $$
  which means that $x_i$ lies in $U$ for all sufficiently large $i$. Since $U$ is compact, we may take a converging
subnet $\{x_{i_j}\}_j$, and then
  $$
  \xi= \lim_j\varphi(x_{i_j}) = \varphi(\lim_jx_{i_j}),
  $$
  proving that $\xi$ lies in the range of $\varphi$.

In order to prove that $\varphi$ is onto we argue by way of contradiction and suppose that $\tsp\E\setminus\varphi(X)$
is nonempty. This is an open set, and the ultra-filters are dense in $\tsp\E$, so we may find some ultra-filter $\xi$ in
$\tsp\E\setminus\varphi(X)$. By (2.13) we may then find $U$ in $\xi$ such that
  $$
  \xi\in D_U\subseteq\tsp\E\setminus\varphi(X).
  $$
  In particular this says that $D_U\cap\varphi(X) = \varnothing$, but this is certainly false because $U$ is nonempty
and $\varphi(U)\subseteq D_U$.

This shows that $\varphi$ is surjective and it is easy to see that it is also injective. Focusing now on proving that
the inverse of $\varphi$ is continuous, we have already seen that $\varphi(U)\subseteq D_U$, for any $U$ in $\E$, and we
claim that in fact $\varphi(U) = D_U$. This is proved as follows: given any $\xi$ in $D_U$, write $\xi= \varphi(x)$, for
some $x$ in $X$, by the surjectivity of $\varphi$. To say that $\xi\in D_U$ is to say that $U\in\xi$, so
  $$
  1 = \bool{U\in\xi} = \bool{U\in\varphi(x)} = \bool{U\in\xi_x} = \bool{x\in U},
  $$
  so $x\in U$, whence $\xi= \varphi(x) \in\varphi(U)$, as required.

This said we see that $\varphi$ restricts to a continuous and bijective map from $U$ to $D_U$. These being compact sets,
we deduce that $\varphi$ is a homeomorphism between them, and hence $\varphi^{-1}$ is continuous on $D_U$. As the $D_U$
form an open cover of $\tsp\E$, we deduce that $\varphi^{-1}$ is continuous on the whole of $\tsp\E$.

Finally, it is easy to prove that, for all $x$ and $y$ in $X$, one has that
  $$
  x\leq y \reliff\xi_x \subseteq\xi_y,
  $$
  so $\varphi$ is an order-isomorphism. \endProof

\section 12 Tight-like groupoids

\def\G{\Gt S}

\firstLine As already mentioned in section (8), our definition of ``\otg/'' requires no compatibility between the order
and the groupoid structures in (see 8.4 and the footnote there) but now we will change that.

\state 12.1. Definition. We shall say that a groupoid $G$ is a \emph{reverse Ehresmann ordered groupoid}, or
\emph{RE-groupoid} for short if, given any $\gamma, \gamma_1,\gamma_2,\delta, \delta_1,\delta_2\in G $, one has that
  \item {(a)} if $\gamma\leq\delta$, then $\gamma^{-1}\leq\delta^{-1}$,
  \item {(b)} if $\gamma_1\leq\gamma_2$, and $\delta_1\leq\delta_2$, then $\gamma_1\delta_1\leq\gamma_2\delta_2$,
provided all multiplications are defined,
  \item {(c)} given $x$ in the unit space $G\ex0 $, such that $d(\gamma)\leq x$, there exists a unique $\delta$ in $G$,
such that $\gamma\leq\delta$, and $d(\delta) = x$,
  \item {(d)} given $x$ in $G\ex0$, such that $r(\gamma)\leq x$, there exists a unique $\delta$ in $G$, such that
$\gamma\leq\delta$, and $r(\delta) = x$,

The reason for the expression ``reverse Ehresmann'' is the fact that the axioms above are exactly the same as in
Ehresmann's classical definition of ordered groupoids \cite[Section 4.1]{lawsonBook}, except that one must consider the
opposite order relation.

Even though the choice between an order relation and its opposite is a question of taste, there are strong reasons for
Ehresmann's choice to be what it is, but one could argue that our choice is not the best one and should therefore be
reversed to conform with Ehresmann's. The basis behind our choice is that we are aiming at the spectral order (see
Proposition 12.2, below) whose essence is to consider a filter $\xi$ \emph{smaller} than another filter $\eta$ when
$\xi$ is \emph{contained} in $\eta$. The devil's advocate point of view, on the other hand, would argue that a
\emph{bigger} topology is called \emph{finer}, and a \emph{bigger} filter may be seen as filtering a \emph{smaller}
region of space, a very strong argument indeed. Weighing the pros and cons we have decided to keep the above definition
after balking at the prospect of viewing an ultra-filter as a minimal element!

Should the reader disagree, implementing the change is straightforward by replacing ``$\leq$'' with``$\geq$''and vice
versa, and perhaps also substituting ``order ideal'' for ``up-set'', and ``down-slice for ``up-slice''.

The formal difference caused by this order reversal notwithstanding, it is surprising that Ehresmann's axioms, which
have been introduced with a very specific purpose of seeing inverse semigroups as groupoids, have turned out to be
absolutely central to the present work, which also aims to relate inverse semigroups to groupoids, but through a very
different optics.

\state 12.2. Proposition. The tight groupoid of any inverse semigroup $S$ is an RE-groupoid when equipped with its
spectral order.

\Proof Referring to the axioms in (12.1), we leave (a) and (c) as easy exercises and move on to (c). For this we write
$\gamma= \gt s\eta$, as usual, and take some
  $$
  \xi\geq d(\gamma) = \eta\ni s^*s,
  $$
  from where we see that $s^*s$ lies in $\xi$, meaning that the germ $\delta:= \gt s\xi$ is well formed, and it clearly
satisfies the requirements in (c). To see that $\delta$ is unique, suppose that $\beta\in\G$ is such that
$\gamma\leq\beta$, and $d(\beta) = \xi$. Then we may write $\gamma= \gt t\eta$, and $\beta= \gt t\xi$, with
$t^*t\in\eta\subseteq\xi$. From
  $$
  \gt s\eta= \gamma= \gt t\eta,
  $$
  we have that there exist $e\in\eta$, satisfying $se = te$. Therefore $e\in\xi$, and
  $$
  \delta= \gt s\xi= \gt{se}\xi= \gt{te}\xi= \gt{t}\xi= \beta,
  $$
  proving that $\delta$ is unique. Point (d) may be proved by a similar method or using (a) to switch between domains
and ranges. \endProof

An ordered groupoid in Ehresmann's sense is called \emph{inductive} (see \cite[Section 4.1]{lawsonBook}) if its unit
space is a semilattice under the restricted order. However the spectral order on tight groupoids is highly unlikely to
posses this property, regardless of whether or not one switches to the opposite order.

Observe that, in any RE-groupoid, the range map is order-preserving by axioms (i) and (ii), and because $r(\gamma) =
\gamma\gamma^{-1}$. By a similar reasoning one proves that the domain map is also order-preserving.

\state 12.3. Proposition. If\/ $G$ is an RE-groupoid, and $U\subseteq G$ is both a bisection and an up-set, then $d(U)$
and $r(U)$ are up-sets in $G\ex0$. Moreover both $r$ and $d$ are order-isomorphisms from $U$ to their ranges.

\Proof Assuming that $x$ and $y$ are elements in $G\ex0$ satisfying $y\geq x\in d(U)$, choose $\gamma\in U$ such that
$d(\gamma) = x$. By (12.1.c) we may find $\delta\geq\gamma$, such that $d(\delta) = y$. Since $U$ is an up-set, we
deduce that $\delta\in U$, so that $y = d(\delta) \in d(U)$, as needed. The proof that $r(U)$ is an up-set is similar.

We have already seen that $r$ and $d$ are order preserving. Suppose now that $\gamma,\delta\in U$ are such that
$d(\gamma)\leq d(\delta)$. Then, again by (12.1.c), there is some $\beta\geq\gamma$, such that $d(\beta) =
d(\delta)$. Since $U$ is an up-set, we have that $\beta\in U$, and since $U$ is a bisection, $\beta= \delta$. This
implies that $\delta\geq\gamma$, and we are done proving that $d$ is an order-isomorphism. The proof that $r$ is an
order-isomorphism is similar. \endProof

The next result will serve as similar purpose as Proposition (8.7).

\state 12.4. Proposition. If\/ $U_1$ and $U_2$ are, at the same time, bisections and up-sets in an RE-groupoid $G$, then
the same properties apply to $U_1U_2$ and $U_1^{-1}$.

\Proof It is well known that $U_1U_2$ and $U_1^{-1}$ are bisections, so we need only prove that they are up-sets. Take a
pair of elements
  $$
  (\gamma_1, \gamma_2)\in(U_1\times U_2)\cap G\ex2,
  $$
  so that $\gamma_1\gamma_2$ is a typical element of $U_1U_2$. Supposing that $\delta$ is another element of $G$ such
that $\delta\geq\gamma_1\gamma_2$, we must prove that $\delta\in U_1U_2$. Observing that
  $$
  r(\delta) \geq r(\gamma_1\gamma_2) = r(\gamma_1) \in r(U_1),
  $$
  we have that $r(\delta)\in r(U_1)$, because $r(U_1)$ is an up-set thanks to (12.3). So there exists some $\delta_1$ in
$U_1$, such that $r(\delta_1) = r(\delta)$, and hence $r(\delta_1) \geq r(\gamma_1)$. As $r$ is an order-isomorphism
from $U_1$ to $r(U_1)$, again by (12.3), we deduce that $\delta_1\geq\gamma_1$.

The fact that $r(\delta_1) = r(\delta)$ allows us to consider the element $\delta_2 := \delta_1^{-1}\delta$, which
satisfies
  $$
  \delta_2 = \delta_1^{-1}\delta\geq\gamma_1^{-1}\gamma_1\gamma_2 = \gamma_2\in U_2.
  $$
  Using that $U_2$ is an up-set, it follows that $\delta_2\in U_2$, whence $\delta= \delta_1\delta_2\in U_1U_2$, as
required. The proof that $U_1^{-1}$ is an up-set is trivial. \endProof

It is well known that the collection of all bisections in a groupoid is an inverse semigroup. The above result then says
that the subset of this semigroup formed by the bisections that are up-sets is an inverse sub-semigroup. This in turn
allows for the introduction of some important notions:

\state 12.5. Definition. Let $G$ be a (not necessarily Hausdorff) ample topological groupoid equipped with a partial
order relation making it a reverse Ehresmann ordered groupoid. We shall then say that:
  \item {(a)} a slice $U\subseteq G$ is an \emph{up-slice} if $U$ is an up-set,
  \item {(b)} the \emph{fundamental inverse semigroup} of $G$, denoted $\Up(G)$, is the collection of all compact
up-slices in $G$,
  \item {(c)} $G$ is a \emph{tight-like groupoid} if its unit space is a tight-like space (Definition 11.1), and
$\Up(G)$ covers $G$.
  \item {(d)} the tight-like groupoids $G$ and $H$ are \emph{isomorphic} if there exists a map $\varphi:G\to H$ that is
an isomorphism of groupoids, a homeomorphism of topological spaces, and an order-isomorphism.

The motivation for the expression ``tight-like groupoid'' is similar to that of ``tight-like spaces''.

\state 12.6. Proposition. The tight groupoid of every inverse semigroup $S$ is tight-like.

\Proof Denoting by $E$ the idempotent semilattice of $S$, the unit space of $\Gt S$ is $\tsp E$, which is tight-like by
(11.2). On the other hand, the fundamental slices $\Delta_s$ are compact up-slices covering $\Gt S$. \endProof

We next plan to use \cite[Proposition 5.4]{actions} to represent a given tight-like groupoid $G$ as the groupoid of
germs for the natural action of $\Up(G)$ on $G\ex0$, and this requires that we first prove the following:

\state 12.7. Lemma. Let $G$ be a tight-like groupoid, let $U$ and $V$ be compact up-slices, and let $\gamma\in U\cap
V$. Then there exists a compact up-slice $W$, such that $\gamma\in W\in U\cap V$.

\Proof Letting $W_1 = U\cap V$, observe that $W_1$ is open because the intersection of open sets is open, and $W_1$ is
an up-set because the intersection of up-sets is an up-set. As a conclusion, $W_1$ is an up-slice. Were $G$ assumed to
be Hausdorff, the intersection of compact sets would be compact and then $W_1$ would be a compact up-slice and the proof
would be finished with this choice of $W_1$.

Not knowing whether or not $G$ is Hausdorff, we need a different strategy. The new strategy will be based on (8.6),
except that this was only proved for tight groupoids, and not for tight-like ones, so we cannot invoke (8.6)
directly. The way out is to observe that $G\ex0$ is a tight-like space and hence it is order-homeomorphic to the tight
spectrum $\tsp E$, for some semilattice $E$, by (11.8). Viewing $E$ as an inverse semigroup, we've seen in the comments
following (2.19) that the tight groupoid of $E$ is $\tsp E$, so (8.6) duly applies to $\tsp E$.

An immediate consequence is then that, given any open up-set $Y\subseteq\tsp E$, and given any $\xi$ in $Y$, there
exists a compact open up-set $Z$, namely one of the $\Delta_{s_i}$'s mentioned in the statement of (8.6), such that
$\xi\in Z\subseteq Y$.

Considering that $d(W_1)$ is an open up-set by (12.3), and that
  $$
  d(\gamma)\in d(W_1)\subseteq G\ex0 \simeq\tsp E,
  $$
  we deduce that there exists a compact open upset $Z\subseteq d(W_1)$, such that $d(\gamma)\in Z$. This said, we claim
that
  $$
  W := d^{-1}(Z)\cap W_1
  $$
  is a compact up-slice. It is a bisection, because it is contained in the bisection $W_1$, it is compact and open
because it is the inverse image of the compact open set $Z$ under the homeomorphism between $W_1$ and $d(W_1)$ given as
the restriction of $d$. So we need only prove that $W$ is an up-set. For this, suppose that $\alpha\geq\beta\in W$, and
let us prove that $\alpha\in W$. Observing that $\alpha\geq\beta\in W\subseteq W_1$, and that $W_1$ is an up-set, we
have that $\alpha$ belongs to $W_1$. Also, since
  $$
  d(\alpha)\geq d(\beta)\in d(W) = Z,
  $$
  and since $Z$ is an up-set, we see that $d(\alpha)\in Z$, so $\alpha\in d^{-1}(Z)$. Therefore
  $$
  \alpha\in d^{-1}(Z)\cap W_1 = W.
  $$
  Summarizing we have proved that $W$ is a compact up-slice and that $\gamma\in W\subseteq U\cap V$, so we are
done. \endProof

With the technical issues out of the way, we have the following application of \cite[Proposition 5.4]{actions}:

\state 12.8. Theorem. Let $G$ be a tight-like groupoid and let $\Theta$ be the natural action of\/ $\Up(G)$ on $G\ex0$,
as described in \cite[5.3]{actions}. Then there is an isomorphism of topological groupoids
  $$
  \mu: \Up(G)\ltimes_\Theta G\ex0 \to G,
  $$
  such that
  $$
  \mu\big(\germ U{d(\gamma)}\Theta\big) = \gamma,
  $$
  whenever $U$ lies in $\Up(G)$, and $\gamma\in U$.

\Proof The hypotheses of \cite[5.4]{actions} hold due to (12.7) and the last condition in (12.5.c), hence leading up to
an isomorphism from $\Up(G)\ltimes_\Theta G\ex0$ to $G$, whose expression is precisely given in \cite[5.4.2]{actions},
which in turn coincides with the expression in the statement. \endProof

We thus arrive at the main result of this section.

\state 12.9. Theorem. Every tight-like groupoid is isomorphic to the tight groupoid of some inverse semigroup. More
precisely, given a tight-like groupoid $G$, its fundamental inverse semigroup $\Up(G)$ is the a unique (up to
isomorphism) flat distributive inverse semigroup whose tight groupoid is isomorphic to $G$, as \otgs/.

\Proof Letting $S = \Up(G)$, we will first prove that $G$ is isomorphic to $\Gt S$, as \otgs/. We begin by observing
that a compact up-slice $U\subseteq G$, seen as an element of $\Up(G)$, is idempotent if and only if $U\subseteq G\ex0$,
in which case it is a compact open up-set in $G\ex0$. Conversely, any compact open up-set $U\subseteq G\ex0$ is a
compact up-slice. In other words, the the idempotent semilattice of $\Up(G)$ coincides with the set of all compact open
up-sets in $G\ex0$, which incidentally was denoted by $\E$ in (11.5), provided we take $X = G\ex0$. Using (11.8) we then
conclude that the map
  $$
  \varphi:x\in G\ex0\mapsto\xi_x\in\tsp\E
  $$
  is an order-homeomorphism. Regarding the natural action $\Theta$ of $S$ on $G\ex0$, already mentioned in (12.8),
vis-\`a-vis the canonical action $\theta$ of $S$ on $\tsp\E$ (see 2.17), we claim that $\varphi$ is a covariant map or,
more precisely, $(id,\varphi)$ is a covariant pair in the sense of (2.21), where $id$ denotes the identity map on $S$.

In order to prove (2.21.1) it suffices to check that, for every $U$ in $S$, and every $x$ in $X$, one has that
  $$
  x\in d(U) \reliff\varphi(x)\in D_{U^*U}.
  $$
  Starting with the consequent, notice that
  $$
  \varphi(x)\in D_{U^*U} \explica{(2.10)}{\reliff} U^*U\in\varphi(x) = \xi_x \reliff x\in U^*U = d(U),
  $$
  hence proving (2.21.1). Assuming that $x$ is in $d(U)$, we must then prove that
  $$
  \varphi\big(\Theta_U(x)\big) = \theta_U\big(\varphi(x)\big). \eqno{(12.9.1)}
  $$
  For this, let $\gamma\in U$ be such that $d(\gamma) = x$. We then have that $ \Theta_U(x) =
\Theta_U\big(d(\gamma)\big) = r(\gamma), $ whence
  $$
  \varphi\big(\Theta_U(x)\big) = \varphi\big(r(\gamma)\big) = \xi_{r(\gamma)},
  $$
  and then, for any $V\in\E$, one has that
  $$
  V\in\varphi\big(\Theta_U(x)\big) \reliff r(\gamma)\in V.
  $$
  On the other hand,
  $$
  V\in\theta_U\big(\varphi(x)\big) \explica{(2.17)}{\reliff} U^*VU\in\varphi(x) = \xi_x \reliff x = d(\gamma)\in U^*VU.
  $$
  As it is clear that $d(\gamma)\in U^*VU$ if and only if $r(\gamma)\in V$, we have shown (12.9.1), proving the desired
covariance.

We then deduce from (2.22) that there exists a continuous homomorphisms
  $$
  id\times\varphi: S\ltimes_\Theta G\ex0 \to S\ltimes_\theta\tsp\E= \Gt S,
  $$
  such that
  $$
  (id\times\varphi)\big(\germ Ux\Theta\big) = \gt U{\varphi(x)},
  $$
  for every $U$ in $S$, and every $x$ in $d(U)$.

Recalling that $\varphi$ is a homeomorphism, we could instead have developed the above argument starting with
$\varphi^{-1}$, rather than $\varphi$, and this would lead to continuous groupoid homomorphism
  $$
  id\times\varphi^{-1}: S\ltimes_\theta\tsp\E\to S\ltimes_\Theta G\ex0,
  $$
  which is clearly the inverse of $id\times\varphi$. As a consequence, we have that $id\times\varphi$ is an isomorphism
of topological groupoids.

Combining this with the isomorphism $\mu$ provided by (12.8), we get the following commutative diagram, where $\psi=
(id\times\varphi)\circ\mu^{-1}$

\begingroup\noindent\hfill\beginpicture\setcoordinatesystem units <0.025truecm, 0.025truecm> \setplotarea x from -40 to
240, y from 120 to -20 \put{\null} at -40 120 \put{\null} at -40 -20 \put{\null} at 240 120 \put{\null} at 240 -20
\put{$S\ltimes_\Theta G\ex0$} at 0 80 \put{$S\ltimes_\theta\tsp E$} at 200 80 \arrow<0.11cm> [0.5, 1.8] from 45 80 to
160 80 \put{$id\times\varphi$} at 102.5 95 \put{$G$} at 100 0 \arrow<0.11cm> [0.5, 1.8] from 23.465 61.228 to 82.821
13.743 \put{$\mu$} at 62.513 49.199 \arrow<0.11cm> [0.5, 1.8] from 117.179 13.743 to 177.339 61.871 \put{$\psi$} at
134.881 47.114 \put{$ = \Gt S$} at 265 80 \endpicture\hfill\null\endgroup

\noindent and a quick computation shows that, for any $\gamma$ in $G$, one has that
  $$
  \psi(\gamma) = \gt U{\varphi(d(\gamma))}, \eqno{(12.9.2)}
  $$
  for any compact up-slice $U$ containing $\gamma$, which always exists by hypothesis.

In order to complete the proof it is now enough to prove that $\psi$ is an order-isomorphism. For this, suppose first
that $\gamma_1$ and $\gamma_2$ are elements of $G$, such that $\gamma_1\leq\gamma_2$. Choosing any compact up-slice $U$
containing $\gamma_1$, we have that $U$ necessarily contains $\gamma_2$. Therefore
  $$
  \psi(\gamma_1) = \gt U{\varphi(d(\gamma_1))} \leq\gt U{\varphi(d(\gamma_2))} = \psi(\gamma_2, )
  $$
  where we have used that $d$ and $\varphi$ are order-preserving.

Conversely, supposing that $\psi(\gamma_1) \leq\psi(\gamma_2)$, we must prove that $\gamma_1 \leq\gamma_2$. Using (8.2)
we may write
  $$
  \psi(\gamma_1) = \gt U{\xi_1}, \and\psi(\gamma_2) = \gt U{\xi_2}, \eqno{(12.9.3)}
  $$
  where $\xi_1$ and $\xi_2$ are in $\tsp\E$, and $U$ is an element of $S$ such that
$U^*U\in\xi_1\subseteq\xi_2$. Observing that $d\big(\psi(\gamma_i)\big) = \varphi\big(d(\gamma_i)\big)$, for $i = 1, 2$,
we must have
  $$
  \xi_i = \varphi\big(d(\gamma_i)\big).
  $$

Should the expressions for the $\psi(\gamma_i)$ in (12.9.3) have come as a result of applying (12.9.2), it would bring
with it the information that $\gamma_i$ lies in $U$, but unfortunately this was not the case. Nevertheless we may prove
that $\gamma_i\in U_i$ as follows: choose a compact up-slice $V_i$ containing $\gamma_i$, so that (12.9.2) implies that
  $$
  \psi(\gamma_i) = \gt{V_i}{\varphi(d(\gamma_i))} = \gt U{\varphi(d(\gamma_i))}.
  $$
  So there exists an idempotent compact up-slice, namely a compact open up-set $W\subseteq G\ex0$, such that $V_iW =
UW$, and
  $$
  \varphi(d(\gamma_i)) \in\dom(\theta_W) = \big\{\xi\in\tsp\E: W\in\xi\big\}.
  $$
  This implies that $ W\in\varphi(d(\gamma_i)) = \xi_{d(\gamma_i)}, $ meaning that $d(\gamma_i)\in W$, and then
  $$
  \gamma_i = \gamma_id(\gamma_i)\in V_iW = UW \subseteq U,
  $$
  showing that $\gamma_i \in U$, as desired. Recalling from (12.3) that $d$ is an order-isomorphism from $U$ to $d(U)$,
the fact that both $\gamma_1$ and $\gamma_2$ lie in $U$, gives
  $$
  \gamma_1\leq\gamma_2 \reliff d(\gamma_1)\leq d(\gamma_2),
  $$
  so it suffices to prove the consequent, and this follows from
  $$
  \varphi(d(\gamma_1)) = \xi_1 \leq\xi_2 = \varphi(d(\gamma_2)),
  $$
  and the fact that $\varphi$ is an order-isomorphism. This proves that $\psi$ is an isomorphism of \otgs/ from $G$ to
$\Gt S$.

The starting point for this proof was to take $S = \Up(G)$, namely the inverse semigroup of all compact up-slices of $G$
(Definition 12.5.b). This is now known to be isomorphic to the inverse semigroup of all compact up-slices of $\Gt S$,
a.k.a.~the tight envelope of $S$ (Definition 8.8.b). In other words, $S$ is isomorphic to its own tight envelope, and
hence $S$ is flat and distributive by (8.10.ii) and (10.7).

Finally, to show that $S$ is the unique flat distributive inverse semigroup whose tight groupoid is isomorphic to $G$,
as \otgs/, suppose that $T$ is another.

Then $\Gt T$ is isomorphic to $\Gt S$, so $T$ is consonant to $S$ by (9.5), and hence $T$ is isomorphic to $S$ by
(10.9). \endProof

The next result is to be interpreted as a duality result along the lines of \cite{lawsonV} and \cite{StoneDuality}.

\state 12.10. Theorem. The correspondence
  $$
  S\mapsto\Gt S
  $$
  is a bijection between the isomorphism classes of flat distributive inverse semigroups and the isomorphism classes of
tight-like groupoids. Its inverse is given by
  $$
  G \mapsto\Up(G),
  $$
  where $\Up(G)$ is the fundamental inverse semigroup of $G$.

\Proof This follows from Proposition (12.6) (the correspondence is well defined, i.e.~$\Gt S$ is tight-like), Theorems
(9.5) and (10.1) (the correspondence is injective up to isomorphism), and Theorem (12.9) (the correspondence is
surjective up to isomorphism, plus the expression for the inverse). \endProof

\section 13 Appendix -- Homeomorphisms between full spectra

\def\Dh{\widehat\D_h}

\firstLine This section is largely independent from the rest of the paper and its sole purpose is to contrast what we
have done so far with the attempt at solving the same fundamental problem that lead to the above development, namely
finding necessary and sufficient conditions under which a homomorphism between inverse semigroups induces an isomorphism
between their groupoids, with the vastly significant difference that we consider Paterson's universal groupoid
\cite[Section 4.3]{paterson} as opposed to tight groupoids.

As we shall see, big differences already appear at the level of semilattices, so we will not go beyond this case,
leaving a further development as a possible research project.

Paterson's construction doesn't care whether or not inverse semigroups possess a zero element, so we will break with our
previous tradition and adopt the same stance. Should inverse semigroups $S$ and $T$ possess a zero, we also will not
care if a homomorphism between them sends zero to zero.

If $E$ is a semilattice in such an inverse semigroup, then $E$ may also fail to have a zero. So, from now on, our
semilattices will no longer be required to possess a zero element, and, even if they do, we will ignore it. A
\emph{character} on such a semilattice $E$ is any \emph{nonzero} homomorphism
  $$
  \varphi:E\to\{0,1\},
  $$
  and the spectrum of $E$, denoted $\hatE$, is the space of all characters with the topology of pointwise convergence.

None of this is in fact such a big break with the case where semigroups and semilattices have zero since we may think of
every semilattice $E$ as a subset of $E\mathop{\dot\cup}\, \{0\}$, and then any character of $E$ may be viewed as a
zero-preserving character of $E\mathop{\dot\cup}\, \{0\}$ that has been restricted to $E$. In other words, zero does not
make a big difference. Tightness does.

\fix From now on we shall fix semilattices $E$ and $F$, as well as a homomorphism $h:E\to F$.

\medskip Given any character $\psi$ on $ F$, notice that the composition $\psi\circ h$ is a homomorphism from $E$ to
$\{0, 1\}$, so it will be a character on $E$ provided it is not identically zero. Setting
  $$
  \Dh= \big\{\psi\in\hatF: \psi\circ h \neq0\big\},
  $$
  we may then define
  $$
  \hath: \psi\in\Dh\mapsto\psi\circ h \in\hatE.
  $$

\state 13.1. Proposition. $\Dh$ is open and $\hath$ is continuous on $\Dh$.

\Proof Let $\psi\in\Dh$, so there exists some $e$ in $E$, such that $\psi(h(e))\neq0$. The set
  $$
  U = \big\{\rho\in\hatF: \rho(h(e))\neq0\big\}
  $$
  is then an open subset of $\hatF$, and clearly $\psi\in U\subseteq\Dh$. This shows that $\Dh$ is open. We leave the
rest of the proof as an easy exercise. \endProof

\state 13.2. Theorem. Given a semilattice homomorphism $h:E\to F$, the following conditions are equivalent:
  \item {(i)} $h$ is an isomorphism of semilattices,
  \item {(ii)} $\Dh= \hatF$, and $\hath$ is a homeomorphism from $\hatF$ onto $\hatE$,
  \item {(iii)} $\Dh= \hatF$, and $\hath$ is bijective.

\Proof \proofImply {i$\imply$ii}{iii} Obvious.

\Medskip\proofImply{iii}{i} Assuming (iii) we first claim that, for every $f_0$ in $F$, there exists some $e$ in $E$,
such that $h(e)\geq f_0$. Consider the (principal) character $\psi$ on $F$ given by
  $$
  \psi(f) = [f\geq f_0], \for f\in F, \eqno{(13.2.1)}
  $$
  where the brackets stand for Boolean value. Then $\psi\in\Dh$ by hypotheses, so $\psi\circ h$ is nonzero, and hence
there exists some $e$ in $E$ such that
  $$
  0 \neq\psi(h(e)) = [h(e)\geq f_0],
  $$
  proving the claim. In order to prove that $h$ is surjective, let $f_0\in F$. We then consider two characters $\psi$
and $\psi'$ on $F$, the first one of which is defined precisely as in (13.2.1), while $\psi'$ is defined for $f$ in $F$,
by
  $$
  \psi'(f) = \clauses{ \cl{1} {\text{if there exists some $e\in E$, such that $f\geq h(e)\geq f_0$, } } \cl{0}
{\text{otherwise.}} }
  $$
  We then claim that $\psi'$ is a character. To see that $\psi'$ respects multiplication it is enough to show that,
given $f_1$ and $f_2$ in $F$, one has that
  $$
  \big(\psi'(f_1) = 1\big) \wedge\big(\psi'(f_2) = 1\big) \iff\psi'(f_1f_2) = 1.
  $$
  Regarding the forward implication, suppose that $\psi'(f_1) = \psi'(f_2) = 1$. Then, for each $i = 1,2$, there exists
$e_i$ in $E$ such that $f_i\geq h(e_i)\geq f_0$. So
  $$
  f_1f_2 \geq h(e_1)h(e_2) = h(e_1e_2)\geq f_0,
  $$
  whence $\psi'(f_1f_2) = 1$. Conversely, assuming that $\psi'(f_1f_2) = 1$, there exist $e$ in $E$ such that
$f_1f_2\geq h(e)\geq f_0$, so, for each $i = 1,2$, we have that
  $$
  f_i \geq f_1f_2 \geq h(e)\geq f_0,
  $$
  whence $\psi'(f_i) = 1$. This shows that $\psi'$ respects multiplication. To see that $\psi'$ is not identically zero,
recall from the first part of this proof that there exists $e$ in $E$ such that $h(e)\geq f_0$, so clearly $\psi'(h(e))
= 1$. So $\psi'$ is indeed a character, as claimed.

Observing that, for every $e$ in $E$, one has that
  $$
  \psi'(h(e)) = 1 \iff h(e)\geq f_0,
  $$
  we have that $\psi'(h(e)) = \psi(h(e))$, or, equivalently, that $\hath(\psi') = \hath(\psi)$. As $\hath$ is supposed
to be one-to-one, we deduce that $\psi' = \psi$, so
  $$
  1 = \psi(f_0) = \psi'(f_0),
  $$
  meaning that there exists some $e$ in $E$, such that $f_0\geq h(e)\geq f_0$. This shows that $h$ is onto.

To show that $h$ is one-to-one, let $e_1, e_2\in E$, with $e_1\neq e_2$. Supposing without loss of generality that
$e_1\not\geq e_2$, let $\varphi$ be the character on $E$ given by
  $$
  \varphi(e) = [e\geq e_2], \for e\in E.
  $$
  Since $\hath$ is assumed to be onto, there exists some character $\psi$ on $F$ such that $\varphi= \hath(\psi)$,
whence
  $$
  \psi(h(e_1)) = \hath(\psi) (e_1) = \varphi(e_1) = 0 \neq1 = \varphi(e_2) = \hath(\psi) (e_2) = \psi(h(e_2)).
  $$
  Consequently $h(e_1) \neq h(e_2)$, as needed. \endProof

Like before, $\hatE$ is equipped with a spectral order relation, namely pointwise order of characters, or inclusion of
filters, in case we adopt this point of view. If a homomorphism $h$, as above, is an isomorphism, then clearly $\hath$
is an order-isomorphism, so the equivalent conditions of Theorem (13.2) could be enlarged by adding that $\Dh= \hatF$,
and $\hath$ is an order-isomorphism.

Should we take inspiration from Corollary (5.7) in defining a notion of consonance in the present context, we would
therefore simply add a synonym for isomorphism.

Of course we could go on and discuss a generalization of Theorem (13.2) to inverse semigroups vis-\`a-vis Paterson's
universal groupoid, but it is already clear that this situation is not nearly as deep as when one focus on tight
groupoids.

\overfullrule0pt

\references

\Bibitem LPA G. Abrams, P. Ara, M. Molina; Leavitt Path Algebras; Lecture Notes in Mathematics, vol. 2191, (2017),
Springer

\Article Boava G. Boava, G. G. de Castro, F. de L. Mortari; Inverse semigroups associated with labelled spaces and their
tight spectra; Semigroup Forum, 94 (2017), no. 3, 582-609

\Article bussExelMeyer A. Buss, R. Exel and R. Meyer; Inverse semigroup actions as groupoid actions; Semigroup Forum, 85
(2012), 227-243

\Article BiceStar T. Bice and C. Starling; Hausdorff tight groupoids generalised; Semigroup Forum, 100 (2020), 399-438

\Article Lisa L. O. Clark, C. Farthing, A. Sims, M. Tomforde; A groupoid generalisation of Leavitt path algebras;
Semigroup Forum, 89 (2014), no. 3, 501-517

\Article Cuntz J. Cuntz; Simple $C^*$-Algebras Generated by Isometries; Commun. Math. Phys., 57 (1977), 173-185

\Article CK J. Cuntz and W. Krieger; A Class of $C^*$-Algebras and Topological Markov Chains; Inventiones Math., 56
(1980), 251-268

\Article DonMil A. P. Donsig and D. Milan; Joins and covers in inverse semigroups and tight C*-algebras;
Bull. Aust. Math. Soc., 90 (2014), no. 1, 121-133

\Article DokExelShift M. Dokuchaev and R. Exel; Partial actions and subshifts; J. Funct. Analysis, 272 (2017), 5038-5106

\Article actions R. Exel; Inverse semigroups and combinatorial C*-algebras; Bull. Braz. Math. Soc. (N.S.), 39 (2008),
191-313

\Article tightrep R. Exel; Tight representations of semilattices and inverse semigroups; Semigroup Forum, 79 (2009),
159-182

\Article reconstru R. Exel; Reconstructing a totally disconnected groupoid from its ample semigroup;
Proc. Amer. Math. Soc., 138 (2010), 2991-3001

\Bibitem covertojoin R. Exel; Tight and cover-to-join representations of semilattices and inverse semigroups; Operator
theory, functional analysis and applications, 183-192, Oper. Theory Adv. Appl., 282, Birkh\umlaut auser/Springer, 2021

\Article ExelStar R. Exel, C. Starling; Self-similar graph C*-algebras and partial crossed products; J. Operator Theory,
75 (2016), no. 2, 299-317

\Article EGS R. Exel, D. Goncalves and C. Starling; The tiling C*-algebra viewed as a tight inverse semigroup algebra;
Semigroup Forum, 84 (2012), no. 2, 229-240

\Article exelLaca R. Exel and M. Laca; Cuntz-Krieger algebras for infinite matrices; J. reine angew. Math., 512 (1999),
119-172

\Article EP R. Exel, E. Pardo; Self-similar graphs, a unified treatment of Katsura and Nekrashevych C*-algebras;
Adv. Math., 306 (2017), 1046-1129

\Article EPdois R. Exel and E. Pardo; The tight groupoid of an inverse semigroup; Semigroup Forum, 92 (2016), no. 1,
274-303

\Bibitem ExelSteinbergOne R. Exel and B. Steinberg; Representations of the inverse hull of a 0-left cancellative
semigroup; arXiv:1802.06281 [math.OA], 2018

\Bibitem ExelSteinberg R. Exel and B. Steinberg; Subshift semigroups; arXiv:1908.08315 [math.OA], 2019

\Article LaLMil S. M. LaLonde and D. Milan; Amenability and uniqueness for groupoids associated with inverse semigroups;
Semigroup Forum, 95 (2017), no. 2, 321-344

\Bibitem lawsonBook M. V. Lawson; Inverse semigroups, the theory of partial symmetries; World Scientific, 1998

\Article lawsonC M. V. Lawson; Compactable semilattices; Semigroup Forum, 81 (2010), 187-199

\Article lawsonB M. V. Lawson; The Booleanization of an inverse semigroup; Semigroup Forum, 100 (2020), no. 1, 283-314

\Bibitem StoneDuality M. V. Lawson; Non-commutative Stone Duality; Semigroups, Algebras and Operator Theory. ICSAOT
2022. Springer Proceedings in Mathematics \& Statistics, vol 436

\Article LawsonLenz M. V. Lawson and D. H. Lenz; Pseudogroups and their \'etale groupoids; Adv. Math., 244 (2013),
117-170

\Article lawsonV M. V. Lawson and A. Vdovina; The universal Boolean inverse semigroup presented by the abstract
Cuntz--Krieger relations; J. Noncommut. Geom., 15 (2021), 279-304

\Article MilanStein D. Milan and B. Steinberg; On inverse semigroup {C}$^*$-algebras and crossed products; Transactions
of the American Mathematical Society, 367 (2015), no. 6, 4043-4071

\Article VNOne F. J. Murray e J. von Neumann; On rings of operators; Ann. of Math., 37 (1936), 116-229

\Article VNTwo F. J. Murray e J. von Neumann; On rings of operators, II; Trans. Amer. Math. Soc., 41 (1937), 208-248

\Article VNFour F. J. Murray e J. von Neumann; On rings of operators, IV; Ann. of Math., 44 (1943), 716-808

\Article OrtegaPardo E. Ortega and E. Pardo; The tight groupoid of the inverse semigroups of left cancellative small
categories; Transactions of the American Mathematical Society, 373 (2020), no. 7, 5199-5234

\Bibitem paterson A. L. T. Paterson; Groupoids, inverse semigroups, and their operator algebras; Birkhauser, 1999

\Bibitem Raeburn I. Raeburn; Graph Algebras; CBMS Reg. Conf. Ser. Math., vol. 103, Amer. Math. Soc., Providence, RI,
2005

\Bibitem Renault J. Renault; A groupoid approach to C*-algebras; Lecture Notes in Mathematics vol.~793, Springer, 1980

\Bibitem Sims A. Sims; Hausdorff \'etale groupoids and their C*-algebras; in A. Sims, G. Szab\'{o}, and D. Williams,
{\it Operator algebras and dynamics: groupoids, crossed products, and {R}okhlin dimension} (F. Perera, ed), pp. 58--120,
Advanced Courses in Mathematics - CRM Barcelona, Springer, (2020)

\Article StarlingTwo C. Starling; Boundary quotients of $\rm C^*$-algebras of right LCM semigroups; J. Funct. Anal., 268
(2015), no. 11, 3326-3356

\Article StarlingThree C. Starling; Inverse semigroups associated to subshifts; Journal of Algebra, 463 (2016), 211-233

\Article Steinberg B. Steinberg; A groupoid approach to discrete inverse semigroup algebras; Adv. Math., 223 (2010),
no. 2, 689-727

\Article VNThree J. von Neumann; On rings of operators, III; Ann. of Math., 41 (1940), 96-161

\Article ZM G. Zeller-Meier; Produits crois\'es d'une C*-alg\`ebre par un group d'automorphismes; J. Math Pures Appl.,
47 (1968), 101-239

\endgroup

\vskip1cm

\def\Address#1#2#3{{\bigskip{\tensc#2} \par\it E-mail address: \tt#3 \par}}

\Address{R. Exel} {Departamento de Matem\'atica, Centro de Ci\^encias F\'{\i}sicas e Mate\-m\'a\-ticas, Universidade
Federal de Santa Catarina, Florian\'opolis, SC, 88040-900, Brazil} {ruyexel@gmail.com}

\close